\documentclass[leqno]{amsart}

\usepackage{amssymb}
\usepackage{amsmath, amsfonts, vmargin, enumerate}
\usepackage{color}
\usepackage{amsthm}
\usepackage{bbm}
\usepackage{verbatim}
\usepackage{mathrsfs}
\usepackage{txfonts}
\usepackage[all]{xy}
\usepackage{mathtools}
\usepackage{mleftright,xparse}
\usepackage{cite}
\usepackage{xspace}
\usepackage{hyperref}

\begin{document}
	\newtheorem{theorem}{Theorem}
	\theoremstyle{definition}
	\newtheorem{proposition}[theorem]{Proposition}
	\newtheorem{lemma}[theorem]{Lemma}
	\newtheorem{remark}[theorem]{Remark}
	\newtheorem{definition}[theorem]{Definition}
	\newtheorem{notation}[theorem]{Notation}
	\newtheorem{conjecture}[theorem]{Conjecture}
	\newtheorem{corollary}[theorem]{Corollary}
	\newtheorem{question}[theorem]{Question}
	\newtheorem{example}[theorem]{Example}
	\newtheorem{problem}[theorem]{Problem}
	\newtheorem{exercise}[theorem]{Exercise}
	\newtheorem{fact}[theorem]{Fact}
	\numberwithin{theorem}{section}
	\numberwithin{equation}{section}

	\def\CC{\mathbb{C}}
	\def\EE{\mathbb{E}}
	\def\GG{\mathbb{G}}
	\def\HH{\mathbb{H}}
	\def\II{\mathbb I}
	\def\MM{\mathbb{M}}
	\def\N{\mathbb{N}}
	\def\NN{\mathbb{N}}
	\def\QQ{\mathbb{Q}}
	\def\RR{\mathbb{R}}
	\def\TT{\mathbb{T}}
	\def\VV{\mathbf{V}}
	\def\Z{\mathbb{Z}}
	\def\ZZ{\mathbb{Z}}
	
	\def\A{{\mathcal A}}
	\def\B{{\mathcal B}}
	\def\C{{\mathcal C}}
	\def\D{{\mathcal D}}
	\def\E{{\mathcal E}}
	\def\F{{\mathcal F}}
	\def\G{{\mathcal G}}
	\def\H{{\mathcal H}}
	\def\I{{\mathcal I}}
	\def\K{{\mathcal K}}
	\def\L{{\mathcal L}}
	\def\M{{\mathcal M}}
	\def\P{{\mathcal P}}
	\def\Q{{\mathcal Q}}
	\def\R{{ R}}
	\def\S{{\mathcal S}}
	\def\T{{\mathcal T}}
	\def\U{{\mathcal U}}
	\def\V{{\mathcal V}}
	\def\W{{\mathcal W}}
	\def\X{{\mathcal X}}
	
	\def\a{\mathfrak{a}}
	\def\b{\mathfrak{B}}
	\def\d{\mathfrak{d}}
	\def\g{\mathfrak{g}}
	\def\h{\mathfrak{h}}
	\def\n{\mathfrak{n}}
	\def\s{\mathfrak{s}}
	\def\z{\mathfrak{z}}

	\newcommand{\abr}[1]{ \langle  #1 \rangle}
	\newcommand{\Norm}[1]{ \left\|  #1 \right\| }
	\newcommand{\set}[1]{ \left\{ #1 \right\} }
	\newcommand{\inner}[1]{ \langle  #1 \rangle}
	\newcommand{\lp}[1]{L_{#1}(\GG)}
	\newcommand{\wlp}[1]{\L^{#1,\infty}}
	\newcommand{\dl}[1]{\delta_{#1}}
	\newcommand{\sg}[1]{\sigma_{#1}}
	\newcommand{\lgm}[1]{L_{\infty}(\GG;\mathbb{M}_{#1})}
	\newcommand{\lm}[1]{L_{\infty}(\mathbb{M}_{#1})}
	\newcommand{\wkp}[1]{\dot{W}^{1,#1}(\GG;\MM_{n})}
	\newcommand{\mlp}[1]{L_{#1}(\GG;\MM_{n})}
	\newcommand{\sects}[1]{\Gamma(#1)}

\newcommand{\quasiriesz}[2]{\R^{#1}_{#2}}
\newcommand{\co}{\mathrm}
\newcommand{\gra}[2]{\smash{\g^{#2}_{#1}}}
\newcommand{\named}{Jacobian\xspace}
\newcommand{\conamed}{Euclidean\xspace}
\newcommand{\itemcase}[1]{\uppercase{#1}}
\newcommand{\Jdet}[1]{\mathrm J_{#1}}

\def\Span{\operatorname{Span}\,}
\def\Lie{\operatorname{Lie}\,}
\def\Ord{\operatorname{Ord}\,}
\def\Id{\operatorname{Id}\,}
	\def\dim{\operatorname{dim}\,}
	\def\p{\partial}
	\def\rp{ ^{-1} }
	\def\rn{\mathbb{R}^{d}}
	\def\sn{\mathbb{S}^{d-1}}
	\def\Exp{\operatorname{\rm Exp}}
	\def\one{\mathbf 1}
	\def\unit{\mathbbm 1}
	\def\id{{\rm id}}
	\def\supp{{\rm supp}}
	\def\ri{{\rm i}}
	
	\def\po{\pi_{1}}
	\def\pt{\pi_{2}}
	\def\vn{{\rm VN}(\GG)}
	\def\vh{{\rm VN_{hom}}(\GG)}
	\def\sym{{\rm sym}}
	\def\aut{\rm Aut}
	\def\cci{C_{c}^{\infty}(\GG)}
	\def\gl{{\rm GL}(d,\RR)}
	\def\glo{{\rm GL}(n_{1},\RR)}
	\def\at{{\rm Aut^{\delta}}(\GG)}
	\def\diag{{\rm diag}}
	\def\cca{C^{\infty}(\GG;\at)}
	\def\BMO{{\rm BMO}(\GG)}
	\def\sgm{\dot{W}^{1,Q}(\GG;\MM_{n_{1}})}
	\def\gln{{\rm GL}(n,\RR)}
	\def\sgn{\dot{W}^{1,Q}(\GG;\MM_{n})}
	\def\ov{\overline}
	\def\Ext{{\rm Ext}}
	\def\Rest{{\rm Rest}}
	
	\def\homo{{\rm hom}}
	\def\Homo{{\rm Hom}}
	\newcommand{\Borel}{\mathfrak M}
	\newcommand{\norm}[1]{\left\lVert#1\right\rVert}


	\title{A $C^{*}$-Algebraic Approach To Principal Symbol Calculus On Filtered Manifolds}
	
	\thanks{{\it 2020 Mathematics Subject Classification:} Primary: 58C50,46L87  Secondary: 58A50 ,47L80 }
	
	\thanks{{\it Key words:} Principal symbol, filtered manifold, $C^{*}$-algebra, stratified Lie group, sub-Riemannian manifold}

	\author{David Farrell}
	\address{School of Mathematics and Statistics, UNSW, Kensington, NSW 2052, Australia}
	\email{david.farrell@unsw.edu.au}

	\author{Fedor Sukochev}
	\address{School of Mathematics and Statistics, UNSW, Kensington, NSW 2052, Australia}
	\email{f.sukochev@unsw.edu.au}

	\author{Fulin Yang}
	\address{Institute for Advanced Study in Mathematics, Harbin Institute of Technology, Harbin 150001, China}
	\email{fulinyoung@hit.edu.cn}

	\author{Dmitriy Zanin}
	\address{School of Mathematics and Statistics, UNSW, Kensington, NSW 2052, Australia}
	\email{d.zanin@unsw.edu.au}

	%

	\begin{abstract} 
		From the viewpoint of $*$-homomorphism on $C^{*}$-algebras, we establish the principal symbol mapping for filtered manifolds which are locally isomorphic to stratified Lie groups.
		Let $\GG$ be a stratified Lie group, and let $M$ be a filtered manifold with a $\GG$-atlas and a smooth positive density $\nu$. 
		For the $C^{*}$-algebra bundle $E_{hom}$ of $M$ constructed from quasi-Riesz transforms on $\GG$, we show that there exists a surjective $*$-homomorphism $$\sym_{M}:\Pi_{M}\to C_{b}(E_{hom})$$ such that 
		$$\ker(\sym_{M})=\K(L_{2}(M,\nu))\subset \Pi_{M}$$
		where the domain $\Pi_{M}\subset\B(L_{2}(M,\nu))$ is a $C^{*}$-algebra and $C_{b}(E_{hom})$ is the $C^{*}$-algebra of bounded continuous sections of $E_{hom}$. Especially, we do not make any assumptions on the lattice of the osculating group of $M$ or the assumption of compactness on manifolds in \cite{DAO3,DAO4}.
		
	\end{abstract}

	\date{}
	\maketitle

	\section{Introduction}
	\subsection{Background}
	
	The theory of pseudodifferential operators (PSDOs) plays an important role in partial differential equations, harmonic analysis and non-commutative geometry.
	A central notion of PSDOs is that of principal symbol, which is roughly a homomorphism from the algebra of pseudo-differential operators into an algebra of functions, \cite[Lemma 5.1]{Kohn-Nirenberg}, \cite[Theorem 5.5]{Joshi}, \cite[pp.54-55]{Treves-book}. Usually, it is defined in a manner inhospitable for operator theorists.  An approach to pseudodifferential calculi based on $C^*$-algebras theory was first suggested by Cordes \cite{Cordes87} (see \cite{melo05} for the case of a compact manifold). Recently, in \cite{DAO1}, a new and hospitable approach (in the setting of PSDOs on $\mathbb{R}^n$) to a principal symbol on a certain $C^{\ast}$-subalgebra $\Pi_{\RR^{n}}$ in $B(L_2(\mathbb{R}^n))$ is proposed; this mapping turns out to be a $\ast$-homomorphism from $\Pi_{\RR^{n}}$ into a commutative $C^{\ast}$-algebra. This $C^{\ast}$-algebra $\Pi_{\RR^{n}}$ introduced in \cite{DAO1} is much wider than the class of all classical compactly based PSDOs of order 0. The key point is that $\Pi_{\RR^{n}}$ is the closure (in the uniform norm) of the $\ast$-algebra of all compactly supported classical PSDOs of order $0$ on $\mathbb{R}^d$. The idea to consider this closure appears in \cite{AS1986-1} (see Proposition 5.2 on page 512) and more recently in \cite{melo05}. This opens an avenue to investigate the theory of PSDOs by using $C^{\ast}$-algebraic approach, which is arguably more ``user friendly".

	In \cite{DAO1}, the authors consider two commutative $C^{*}$-algebras $\A_{1}(\RR^{n})=\CC+C_{0}(\RR^{n})$ and $\A_{2}(\RR^{n})=C(\mathbb{S}^{n-1})$. These algebras are respectively represented on Hilbert space $L^{2}(\RR^{n})$ via representations
	\begin{align*}
		\pi_{1}(f)=M_{f},\quad f\in \A_{1}(\RR^{n}),
	\end{align*}
	and 
	\begin{align*}
		\pi_{2}(g)=g\,\biggl(\frac{\ri\nabla}{\sqrt{\Delta}}\biggr),\quad g\in \A_{2}(\RR^{n}),
	\end{align*}
	where $\nabla$ denotes the gradient operator and $\Delta$ denotes the Laplacian on $\RR^{n}$. The $C^{*}$-algebra $\Pi_{\RR^{n}}$ is generated by $\pi_{1}(\A_{1}(\RR^{n}))$ and $\pi_{2}(\A_{2}(\RR^{n}))$, and $\Pi_{\RR^{n}}$ is called the domain of the principal symbol mapping. 
	They then show that there is a unique norm-continuous $*$-homomorphism
	\begin{align*}
		\sym:\Pi_{\RR^{n}}\rightarrow \A_{1}(\RR^{n})\otimes_{\rm min}\A_{2}(\RR^{n})
		\quad\mbox{with}\quad \ker(\sym)=\K(L_{2}(\RR^{n}))\subset\Pi_{\RR^{n}}
	\end{align*}
	such that 
	\begin{align*}
		\sym(\pi_{1}(f))=f\otimes\one \quad\mbox{and}\quad 
		\sym(\pi_{2}(g))=\one\otimes g.
	\end{align*}
	The $*$-homomorphism $\sym$ is called the principal symbol mapping on $\RR^{n}$.
	Later on, this idea is extended to the settings of non-commutative tori and non-commutative Euclidean space \cite{DAO2}, smooth compact manifolds \cite{DAO3}, Heisenberg groups and compact contact manifolds \cite{DAO4}.
	
	In sub-Riemannian manifold settings, the principal symbol of a Heisenberg pseudodifferential operator in  \cite{BG88,Taylor84} is defined in local coordinates only, so the definition a priori depends on the choice of these coordinates. A coordinate-free definition of the principal symbol in this setting was given first by Ponge in \cite{Po08}. The principal symbol of an operator is a polyhomogeneous function for the appropriate dilations on the duals of the osculating groups.  Another invariant definition was given in \cite{EMbook,E04} as a section over a bundle of jets of vector fields representing the bundle of osculating groups. 
	In \cite{V2010b}, van Erp developed the $C^{\ast}$-algebraic approach to the principal symbol map in the Heisenberg calculus. This $C^{\ast}$-algebraic approach is convenient for the purposes of index theory, but it is not useful when one is interested in regularity properties of operators.  
	Therefore, a new $C^{\ast}$-algebraic approach is suggested on compact contact manifold in \cite{DAO4}. 
	Let $\HH^{n}$ be a Heisenberg group and let $\{R_{k}\}_{k=1}^{2n}$ be the left-translation invariant Riesz transforms. Consider two $C^{*}$-algebras 
	$$\A_{1}(\HH^{n})=\CC+C_{0}(\HH^{n})\quad \mbox{and}\quad \A_{2}(\HH^{n})=C^{*}(\{R_{k}\}_{k=1}^{2n}),$$ 
	where $\A_{2}(\HH^{n})$ represents a $C^{*}$-algebra generated by $\{R_{k}\}_{k=1}^{2n}$.
	Obviously, $\A_{1}(\HH^{n})$ is commutative but $\A_{2}(\HH^{n})$ is non-commutative. The construction of $\A_{2}(\HH^{n})$ differs from that on $\RR^{n}$ as given in the prior paper \cite{DAO1} by the same authors, however from translation bi-invariance of Riesz transforms together with functional calculus one can show that the $\A_{2}(\RR^{n})$ constructed therein is generated by the Riesz transforms and so the constructions coincide in the Euclidean case. 
	On $\HH^{n}$, the two algebras are separately represented on Hilbert space $L_{2}(\HH^{n})$ via the following representations
	\begin{align*}
		\pi_{1}(f)=M_{f},\quad f\in \A_{1}(\HH^{n}),
	\end{align*}
	and
	\begin{align*}
		\pi_{2}(g)=g,\quad g\in \A_{2}(\HH^{n}).
	\end{align*}
	The domain of the principal symbol mapping $\Pi_{\HH^{n}}$ is the $C^{*}$-algebra generated by $\pi_{1}(\A_{1}(\HH^{n}))$ and $\pi_{2}(\A_{2}(\HH^{n}))$.
	In \cite{DAO4}, they show that there is a unique norm-continuous $*$-homomorphism
	\begin{align*}
		\sym:\Pi_{\HH^{n}}\rightarrow \A_{1}(\HH^{n})\otimes_{\rm min}\A_{2}(\HH^{n})
		\quad\mbox{with}\quad \ker(\sym)=\K(L_{2}(\HH^{n}))\subset\Pi_{\HH^{n}}
	\end{align*}
	such that 
	\begin{align*}
		\sym(\pi_{1}(f))=f\otimes\one \quad\mbox{and}\quad 
		\sym(\pi_{2}(g))=\one\otimes g.
	\end{align*} 
	This construction is also based on the $C^{*}$-algebraic point but doesn’t involve operator kernels. It directly treats the model operators as operators in the Hilbert space $L_{2}(\HH^{n})$ and elements of
	the von Neumann algebras on $\HH^{n}$. 
	After then, by applying schr\"{o}dinger representations of Heisenberg group, the paper \cite{DAO4} establishes invariance of $\Pi_{\HH^{n}}$ and equivariance of $\mathrm{sym}$ under local Heisenberg diffeomorphisms \cite[definition 1.3.5]{DAO4} on compact contact manifolds. 
	This gives a good definition of the domain of principal symbol mapping on contact manifolds:
	given a compact contact manifold $M$ with a Heisenberg atlas $\{(\U_{i},h_{i})\}_{i\in\II}$ and with a smooth positive density $\nu$, the domain $\Pi_{M}$ is the set of operators $T\in\B(L_{2}(M,\nu))$ satisfying
	\begin{enumerate}[\rm(i)]
		\item \itemcase{f}or all $i\in\II$ and $\phi\in C_{c}(\U_{i})$, the operator $M_{\phi}TM_{\phi}$ regarded in $\B(L_{2}(\HH^{n}))$ belongs to $\Pi_{\HH^{n}}$, and
		\item \itemcase{f}or each $\phi\in C_{c}(M)$, the operator $[T,M_{\phi}]$ is compact on $L_{2}(M,\nu)$.
	\end{enumerate}
	This $\Pi_{M}$ is a $C^{*}$-algebra. They prove that
	there exists a surjective $*$-homomorphism $$\sym_{M}:\Pi_{M}\to C(E_{M})$$ such that $$\ker(\sym_{M})=\K(L_{2}(M,\nu)),$$
	where $E_{M}=(M,\A_{2}(\HH^{n}),\overline{\omega})$ is the $C^{*}$-algebra bundle with $\overline{\omega}$ a family of continuous mappings associated to a Heisenberg atlas \cite[definition 1.3.6]{DAO4} and $C(E_{M})$ is the $C^{*}$-algebra of bounded continuous sections of $E_{M}$. They then establish the Connes' trace theorem and Connes' integral formula on the compact contact manifold. This $*$-homomorphism significantly promotes the development of non-commutative geometry on contact manifolds.

	Motivated by the results and applications on $\RR^{n}$, compact manifolds and compact contact manifolds, we aim to investigate the principal symbol mapping on a class of manifolds generalising that of compact contact manifolds; namely $\GG$-filtered manifolds. These are those (not necessarily compact!) manifolds which are equipped with a filtration of their tangent bundles which is locally isomorphic to the natural filtration on some stratified Lie group $\GG$. We will show that:
	\begin{enumerate}[\rm(i)]
		\item If $\GG$ is a stratified Lie group,
		then there is a surjective $*$-homomorphism $\sym:\Pi\rightarrow\A_{1}\otimes_{\rm min}\A_{2}$ with $\ker(\sym)=\K(L_{2}(\GG))\subset\Pi$ such that, for all $f\in\A_{1}$ and $g\in\A_{2}$,
		\begin{align*}
			\sym(\po(f))=f\otimes1\quad\text{and}\quad\sym(\pt(g))=1\otimes g,
		\end{align*}
		where $\A_{1}=\CC+C_0(\GG)$, $\A_{2}=C^{*}\Big(\{\quasiriesz{A}{k}:A\in\at\}_{k=1}^{n_{1}}\Big)$ and $\Pi=C^{*}\Big(\{\po(\A_{1}),\pt(\A_{2})\}\Big)$.

		\item Given a $\GG$-filtered manifold $M$ (see Definition \ref{G-filtered manifold}), 
		let $\nu$ be a smooth positive density on $M$.  Let the domain $\Pi_{M}$ be as before with $L_{2}(\GG)$ and $\Pi$ in place of $L_{2}(\HH^{n})$ and $\Pi_{\HH^{n}}$, and with an additional condition that there exists a sequence $(\phi_k)_{k\in\NN}$ in $C_c(M)$ such that $M_{\phi_k}T\to T$ uniformly in $\B(L_{2}(M,\nu))$. 
		Then there exists a surjective $*$-homomorphism $\sym_{M}:\Pi_{M}\to C_{b}(E_{hom})$ such that $$\ker(\sym_{M})=\K(L_{2}(M,\nu))\subset\Pi_{M},$$
		where $E_{hom}=(M,\A_{2},\omega)$ is the $C^{*}$-algebra bundle with $\omega$ a family of continuous mappings from the $\GG$-atlas (see Definition \ref{G-atlas}) and $C_{b}(E_{hom})$  is the $C^{*}$-algebra of bounded continuous sections of $E_{hom}$ (more details, see subsection \ref{symbol on manifold-section}).
	\end{enumerate}
	It is worth noting that we do not ask the compactness condition on manifold $M$ and we do not make any lattice assumption on the osculating group of $M$. So our results are suitable for more settings. When $\GG$ is $\HH^{n}$, by schr\"{o}dinger representations on Heisenberg group, our $C^{*}$-algebra $\A_{2}$ is identical to the $C^{*}$-algebra $\A_{2}(\HH^{n})$. Therefore, we recover the previous result if $M$ is a compact contact manifold.
	
	Comparing with the methods on Heisenberg groups and compact manifolds, the main difficulties here are that the unitary representation from the orbit theory is obscure and there is no lattice sub-group for more general stratified Lie groups, and the globalisation theorem for the existence of a $*$-homomorphism introduced in \cite{DAO3} is unsuitable for our settings. On compact contact manifold, the unitary representation and lattice of Heisenberg group guarantee the invariance of $\Pi_{\HH^{n}}$ and equivariance of $\mathrm{sym}$ so that they combine the (compact) globalisation theorem in \cite{DAO3} to  produce the existence of the principal symbol. 
	This requires us to rethink the key purpose of the unitary representation and lattice, and revisit the globalisation theorem in \cite{DAO3}.

	To overcome the difficulties, we bypass the use of unitary representation and lattice group, and provide a new globalisation theorem. We suggest the $C^{*}$-algebra constructed from  quasi-Riesz transforms of the strata-preserving automorphisms on stratified Lie groups. Combining functional calculus on \cite{MSZ-Cwikel} and the property of stratification on stratified Lie group, such a $C^{*}$-algebra allow us to avoid the use of the unitary representation. To avoid lattice group, we use  a covering lemma (see Lemma \ref{bd multiplicity covering}) and then use the technique of partition of unity arising from the covering. 
	After carefully investigating the globalisation theorem in \cite{DAO3}, 
	we add the condition in $\Pi_{M}$: there exists a sequence $(\phi_k)_{k\in\NN}$ in $C_c(M)$ such that $M_{\phi_k}T\to T$ uniformly in $\B(L_{2}(M,\nu))$. This allow us to establish a new globalisation theorem for the existence of a $*$-homomorphism without the compactness assumption on the manifold.
	Otherwise, for our purpose, we also use the notions of $\GG$-diffeomorphism on stratified Lie group and $\GG$-atlas on filtered manifold. 
	Our approach, throughout the paper, is based on the $C^{\ast}$-algebraic point of view and was developed in different settings in previous works \cite{DAO1,DAO2,DAO3,DAO4}. It is closer to that of \cite{Taylor84,V2010b}. 
	Similar to \cite{DAO3,DAO4} but unlike \cite{V2010b}, we don't use operator kernels while directly apply tools of operator theory on stratified Lie groups \cite{MSZ-Cwikel}. 
	As an expectation, we hope that our principal symbol mapping can benefit to the development of the theory of PSDOs \cite{GV2022,JE-arxiv2023,RS1976,DGD2007,Po08,VY19,CY2022,GV2023}, the index theory \cite{AS1986-1,AS1986-3,MK2022,GV2022,MK2022,V2010,V2010b} and non-commutative geometry \cite{Connes1994,Po07,CY2023} on sub-Riemannian manifolds.

	To state our main results clearly, let us give a brief introduction to basic structures and analysis on stratified Lie groups and filtered manifolds.
	
	\subsection{Stratified Lie Groups and Filtered Manifolds}
	
	\subsubsection{Stratified Lie Groups}
	Let $\g$ be a real Lie algebra. Then $\g$ is said to be nilpotent if there exist $l>0$ such that $\g^{(l)}=\{0\}$ where $\g^{(i)}=[\g^{(i-1)},\g]$ $1\leq i\leq l$ and $\g^{(0)}=\g$. A \textit{grading} of $\g$ is a sequence $\g_1, \g_2, \dotsc$ of linear subspaces of $\g$ 
	such that $[\g_{k},\g_{m}]\subset\g_{k+m}$ for all $k,m\geq 1$ and such that $\g$ decomposes as a direct sum
	$$\g=\bigoplus_{k=1}^{\infty}\g_{k}.$$
	If $\g$ is equipped with a grading, it is said to be \textit{graded}. Further, the grading is said to be a \textit{stratification} if, in addition, $\g_{1}$ generates $\g$ as a Lie algebra \cite[Section 3.1]{FR-Book}. A Lie algebra equipped with a stratification is said to be stratified.
	The homogeneous dimension $Q$ of a graded Lie algebra $\g$ is defined to be
	$$Q=\sum_{k=1}^{\infty}k\cdot\dim(\g_{k}).$$
	Obviously, when $\g$ is graded and nilpotent, the homogeneous dimension is finite.
	A connected and simply connected Lie group $\GG$ is called stratified if its Lie algebra $\g$ is stratified. A stratified Lie group $\GG$ is diffeomorphic to $\g$ via the exponential mapping $\Exp$ \cite[Proposition 1.2]{FS1982}. 
	The Lie algebra $\g$ may also be identified with the set of all vector fields on $\GG$ commuting with the action of $\GG$ on itself by left translation.
	
	Throughout the remainder of the paper, let $\GG$ be a $d$-dimensional stratified Lie group with Lie algebra $\g$, stratification
	$$\g=\bigoplus_{n=1}^\iota\g_n,$$
	where $\iota$ is the greatest integer such that $\g_\iota\neq\{0\}$. We will denote the identity of $\GG$ by $o$ or $0$. Further, let $\smash{\{Z_{k}^{(i)}:k=1,2,\,\dotsc,n_{i},i=1,2,\,\dotsc,\iota\}}$ be a basis for $\g$ adapted to the stratification in the sense that, for each $1\leq i\leq \iota$, the set $\smash{\{Z_k^{(i)}:k=1,2,\,\dotsc,n_{i}\}}$ is a basis for $\g_i$. We will use the following notation,
	\begin{notation}
		For all $1\leq k\leq n_1$, we will denote $Z^{(1)}_k$ by $X_k$.
	\end{notation}
	
	\begin{notation}\label{lindexing}
		Also, we will take the liberty to index the basis $\smash{(Z^{(i)}_k)}$ using only one index; that is, if $1\leq i\leq \iota$ and $1\leq k\leq n_i$ we will write $Z_j$ for $\smash{Z^{(i)}_k}$, where $j=n_1+\cdots+n_{i+1}+k$. We will also write so that, if $1\leq j\neq d$,
		$$Z_j=Z^{(\Ord(Z_j))}_k$$
		where $\displaystyle{k=j-\sum_{i=1}^{\Ord(Z_j)-1}}$. Similarly, we will write $\partial/\partial x_j$ for $\smash{\partial/\partial x^{(i)}_k}$ (see Remark \ref{canoncoords} below). Conversely, the double-indexing scheme for $\smash{Z^{(i)}_k}$ is copied over to the standard basis $(e_j)$ of $\RR^d$. That is, we write $\smash{e^{(i)}_k}$ for $e_j$ when $j=n_1+\cdots+n_{i-1}+k$.
	\end{notation}
	
	\begin{remark}\label{canoncoords}
		\par We use the diffeomorphic exponential mapping $\g\to \GG$ to identify the underlying manifold of $\GG$ with that of $\g$. Further, we use the given basis to identify the underlying vector space of $\g$ with $\RR^d$. Thereby we may think of $\GG$ both $\g$ and $\GG$ having underlying set $\RR^d$, although this identification relies on a fixed choice of basis $\smash{\{Z^{(i)}_k\}}$. We will refer to this coordinate system on both $\g$ and $\GG$ as \textit{\named coordinates}, and denote them by $\smash{x^{(i)}_k}$. The induced coordinates on tangent spaces will be denoted by $\smash{\partial/\partial x^{(i)}_k}$ or simply $\smash{\partial^{(i)}_k}$ and called \textit{\conamed coordinates}. The distinction between \named and \conamed coordinates on $\g$ will feature prominently, starting in Remark \ref{basis-rep remark} below.\par
	\end{remark}
	
	\begin{remark}
		The grading on $\g$ induces a natural semigroup $(\delta_r)_{r>0}$ of linear maps on $\g$, given by
		$$\delta_r(Z^{(i)}_k)=r^iZ^{(i)}_k.$$
		The $\delta_r$ are called the \textit{dilations}, and are Lie algebra automorphisms on $\g$ and Lie group automorphisms on $\GG$ (using the identification $\Exp:\g\to\GG$, see Remark \ref{canoncoords} above). These mappings play the same role in analysis on stratified groups that scalar multiplication plays in the Euclidean case, and allow us to consider homogeneity of functions, vector fields and operators, defined as usual, suited to the group structure at hand.
	\end{remark}
	
	\begin{remark}\label{basis-rep remark}
		Each $\smash{Z_k^{(i)}}$ can be regarded as a left-invariant vector field on $\GG$, and by \cite[Remark 1.4.6]{BLU2007}, these vector fields are $\delta_{r}$-homogeneous of degree $i$ and may be written, with respect to \conamed coordinates, as
		\begin{align}\label{basis-rep}
			Z_{j}^{(i)}=\frac{\p}{\p x_{j}^{(i)}}+\sum_{h=i+1}^{\iota}\sum_{k=1}^{n_{h}}a_{j,k}^{(i,h)}(x^{(1)},\,\dotsc,x^{(h-i)})\frac{\p}{\p x_{k}^{(h)}},
		\end{align}
		where the $\smash{a^{(i,j)}_{j,k}}$ are $\delta_r$-homogeneous polynomial functions of degree $h-i$. In the linear indexing scheme (see Notation \ref{lindexing}), we write (\ref{basis-rep}) as
		\begin{equation}\label{basis-rep-linear}
			Z_i=\frac{\partial}{\partial x_i}-\sum_{j=n_1+1}^d \alpha_{ij}\frac{\partial}{\partial x_j}.
		\end{equation}
	\end{remark}
	
	The Haar measure on $\GG$ is identical to the pushforward of the Lebesgue measure on $\g$ under $\Exp$,
	and write $\lp{p}$ as the $L_{p}$-space that are always defined with the Haar measure. We omit the measure when writing integrals of $f\in\lp{1}$. 
	That is,
	\begin{align*}
		\int_{\GG}f:=\int_{\g}f(\xi)\,d\xi.
	\end{align*}
	There are many equivalent symmetric homogeneous norms on $\GG$ and denote, by $B(x,r)$, $r>0$, the open ball induced by the Koranyi gauge $\rho$ (see \cite{BLU2007}) i.e. a symmetric homogeneous norm satisfying that there is a constant $A_{0}\geq1$ such that 
	\begin{align}\label{triangle inequality}
		\rho(xy)\leq A_{0}(\rho(x)+\rho(y))\quad\mbox{for}\quad x,y\in\GG.
	\end{align}
	The measure of the ball $B(x,r)$ is a $Q$-homogeneous function with respect to $r>0$.
	
	\begin{remark}\label{linearity}
		We recall that Lie group endomorphisms of $\GG$ are linear with respect to the vector space structure on $\GG$ induced by the identification with $\g$ via the exponential. This follows from the Campbell-Baker-Hausdorff formula, since, if $x\in\GG$,
		$$x^2=x+x+[x, x]+\frac{1}{2}\left([x, [x, x]]+[x, [x, x]]\right)+\cdots=x+x$$
		and similarly $x^{-1}=-x$ hence applying induction and continuity we get that $x^n=nx$, $x^{p/q}=(p/q)x$ and hence $x^t=tx$ for all $t\in\RR$, thus if $A:\GG\to\GG$ is a homomorphism we get $A(tx)=A(x^t)=(A(x))^t=tA(x)$, and of course, any differentiable map satisfying this agrees with its derivative and therefore is linear. Thus, automorphisms of $\aut(\GG)$ can be represented by matrices (relative to the \named coordinates on $\GG$). This fact will be used heavily throughout the paper.
	\end{remark}

	\subsubsection{Filtered Manifolds}\label{Carnot manifold-section}\
	Recall that, if $E$ is a smooth vector bundle over a smooth manifold $M$, with projection $\pi:E\to M$, then a \textit{smooth section} of $E$ is a smooth map $X:M\to E$ such that $\pi \circ X=\Id_M$, i.e., $X_x=X(x)\in \pi^{-1}(\{x\})=E_x$ for all $x\in M$. The set of smooth sections of $E$ is denoted by $\Gamma(E)$. Note that, if $F$ is a smooth vector sub-bundle of $E$, then $\Gamma(F)=\{X\in\Gamma(E) : X_x\in F_x\,\text{for all}\,x\in M\}$.
	
	\begin{definition}
		Let $M$ be a smooth manifold. A \textit{horizontal distribution} or, simply, \textit{distribution} on $M$ is a smooth vector sub-bundle $H\subset TM$ of the tangent bundle $TM$ of $M$. If $H$ is a distribution on $M$, a vector field $X\in\Gamma(TM)$ is said to be \textit{horizontal} if $X\in\Gamma(H)$. Also, a piecewise smooth curve $\gamma:I\to M$, where $I\subset\RR$ is an interval, is said to be \textit{horizontal} with respect to $H$ if, for each $t\in I$ for which $\gamma$ is differentiable, we have $\dot\gamma(t)\in H_{\gamma(t)}$.
	\end{definition}

	\begin{definition}
		Let $M$ be a smooth manifold. A horizontal distribution $H$ on $M$ is said to be \textit{bracket generating} if, for each $x\in M$,
		$$T_xM=\Span\Bigl\{[V_1, V_2,\,\hdots, V_k]_x : k\geq 1, V_1,\,\hdots, V_k\in\Gamma(H)\Bigr\}$$
		where we use the shorthand notation $[V_1,\,\hdots, V_k]$ to mean $[V_1, [V_2,\,\hdots, [V_{k-1}, V_k]\hdots]]$.
	\end{definition}
	
	\begin{definition}\label{subriemannianmanifold}
		A \textit{sub-Riemannian manifold} is a smooth manifold $M$ equipped with a bracket generating horizontal distribution $H$ and a field $g\in\Gamma(H^*\otimes H^*)$ of inner products.
	\end{definition}
	
	A connected sub-Riemannian manifold carries a natural metric space structure given by the \textit{Carnot-Carath\'eodory} metric, defined as follows. First, if $\gamma:[0,1]\to M$ is a horizontal path, define
	$$L_{cc}(\gamma):=\int_0^1 g_{\gamma(t)}(\dot\gamma(t),\dot\gamma(t))^{1/2}\,dt.$$
	Now, for all $x,y\in M$, define
	$$d_{cc}(x,y):=\inf\set{L_{cc}(\gamma):\gamma:[0,1]\to X\text{ horizontal}, \gamma(0)=x, \gamma(1)=1}.$$
	The metric space axioms follow by the usual proof, but it is not a-priori guaranteed that the infimum will always be finite. This is the content of the celebrated \textit{Chow-Rashevskii Theorem}.
	
	\begin{theorem}[Chow-Raveskii Theorem]\cite{chow, rashevskii}
		Let $M$ be a smooth manifold equipped with a bracket generating horizontal distribution $H$. Then, for all $x,y\in M$, if there is a path from $x$ to $y$, there is a horizontal path from $x$ to $y$.
	\end{theorem}
	
	\begin{definition}\label{filtereddef}
		A \textit{filtration} on a smooth manifold $M$ is an increasing sequence $H^1\subset H^2\subset H^3\subset \hdots$ of distributions on $M$ such that
		\begin{enumerate}
			\item\label{step} \itemcase{t}here is $n\in\NN$ such that $H^n=TM$, and
			\item\label{filtration} \itemcase{f}or each $V\in\Gamma(H^i)$ and $W\in\Gamma(H^j)$, one has $[V,W]\in\Gamma(H^{i+j})$.
		\end{enumerate}
		A manifold $M$ equipped with a filtration $(H^n)_{n\in\NN}$ is called a \textit{filtered manifold}.
	\end{definition}
	
	\begin{remark}
		Filtered manifolds are also called \textit{Carnot manifolds} in, e.g., \cite{Po19, MK2022,S2013}.
	\end{remark}
	
	\begin{definition}
		Filtered manifold $(M, H)$ is called \textit{equiregular} if
		\begin{equation}\label{brackgenerationfilt}
			\Gamma(H^{n+1})=\Gamma(H^n)+{\rm Span}([\Gamma(H^n),\Gamma(H^1)]),\quad n\geq 1.
		\end{equation}
		Equiregular filtered manifolds are also called \textit{equiregular Carnot-Carath\`eodory} manifolds, or ECC manifolds for short.
	\end{definition}
	
	If $(H^j)_{j\in\NN}$ is a filtration on $M$ and $\iota\in\NN$ is the least integer satisfying (\ref{step}), we say that $(H^j)_{j\in\NN}$ is an \textit{$\iota$-step filtration}, and that $M$ is an \textit{$\iota$-step filtered manifold} or a \textit{filtered manifold of step $\iota$.}
	
	If $M$ is a filtered manifold, with filtration $(H^j)_{j\in\NN}$, then the distribution $H^1$ is bracket generating, and the bundle $H^n$ is determined by $H^1$. Therefore, if a filtered manifold is equipped with a field of inner products on $H^1$, it forms a sub-Riemannian manifold. The converse is not true: that is, if $H$ is a bracket generating distribution on a manifold $M$, it is not true that formula (\ref{brackgenerationfilt}) defines a bracket generated filtration, as the $H^n$ may fail to even be vector bundle. A well known example of this failing is seen in the \textit{Martinet distribution}, which is the bracket generating distribution on $\RR^3$ given as the kernel of the $1$-form $dz-x^2dy$. 
	
	\begin{remark}\label{gradingfiltration}
		If the tangent bundle $TM$ of a smooth manifold $M$ carries a grading, i.e.,
		$$TM=\bigoplus_{i=1}^\iota V_i$$
		for distributions $V_1, V_2, \hdots, V_\iota$ on $M$ satisfying $[X, Y]\in \Gamma(V_{i+j})$ for all $X\in\Gamma(V_i)$ and $Y\in\Gamma(V_j)$, then we obtain a sequence of distributions satisfying (\ref{step}) of Definition \ref{filtereddef} via the formula $$H^n:=\sum_{i=1}^n V_i.$$
	\end{remark}
	\begin{definition}\label{stratifiedfiltration}
		Let $\GG$ be a graded Lie group. We obtain a grading on $T\GG$ by left-translating the grading of $\g$; i.e.,
		$$\gra{n}{x}:=((L_x)_*)_o\g_n$$
		where $L_x$ is left-multiplication map $L_x(y)=xy$, and $((L_x)_*)_o$ is its pushforward at $o$; see Section \ref{pushforwardsubsub}. Therefore, we obtain a filtration on $\GG$ via the construction in Remark \ref{gradingfiltration}. This filtration is called \textit{the standard filtration on $\GG$}. 
	\end{definition}
	\begin{remark}\label{spanner}
		In the context of Definition \ref{stratifiedfiltration}, if $(Z^{(i)}_k)$ is a basis for $\g$ adapted to the grading, i.e., $\g_i=\Span\{Z^{(i)}_k : 1\leq k\leq n_i\}$, then, by left-invariance,
		\begin{align*}
			\gra{n}{x} &= ((L_x)_*)_o\big(\Span\{Z^{(i)}_k : 1\leq k\leq n_i\}\big) = \Span\{((L_x)_*)_o(Z^{(i)}_k) : 1\leq k\leq n_i\}=\Span\{Z^{(i)}_k|_x : 1\leq k\leq n_i\}
		\end{align*}
		in particular,
		$$\gra 1 x = \Span\{X_1|_x,\hdots X_{n_1}|_x\},$$
		where $X_j=Z^{(1)}_j$ as usual.
	\end{remark}
	\begin{remark}
		In the context of Definition \ref{stratifiedfiltration}, the statement that $\GG$ is stratified, is equivalent to $\GG$ being equiregular as a filtered manifold when equipped with its standard filtration.
	\end{remark}
	
	\begin{remark}
		Let $(M, H)$ be a filtered manifold and $\Omega\subset M$ an open subset. By the usual identification of $T_x\Omega$ with $T_xM$, for all $x\in \Omega$, we can ``restrict'' $H$ to a filtration $H'$ on $\Omega$ by setting $H'_x=H_x$ for all $x\in\Omega$, making $\Omega$ a filtered manifold in its own right. We will implicitly supply open subsets of filtered manifolds with the filtration inherited in this way from their ambient space throughout the paper.
	\end{remark}
	
	\subsubsection{Pushforwards Of Differential Operators}\label{pushforwardsubsub}
	
	We will use some notation from differential geometry. Suppose $\Phi:M\to N$ is a smooth mapping between smooth manifolds, such as for instance open subsets of $\GG$. If $x\in M$, the \textit{pushforward} of $\Phi$ at $x$ is the linear map $(\Phi_*)_x:T_xM\to T_{\Phi(x)}N$ given by $(\Phi_*)_xX(f)=X(f\circ\Phi)$ for all $X\in T_xM$, where $X$ is realised as a linear operator $C^\infty(M)\to \RR$.
	
	As is well known, if $\varphi$ and $\psi$ are coordinate charts around $x\in M$ and $\Phi(x)\in N$, respectively, then
	\begin{equation}\label{pfptcoord}
		(\Phi_*)_x\Bigl(\sum_{i=1}^d a_i\frac{\partial}{\partial \varphi_i}\Bigr)=\sum_{i,j=1}^d \biggl(a_i \frac{\partial}{\partial\varphi_i}\Phi_j\biggr)(x)\frac{\partial}{\partial \psi_j}\Bigg|_{\Phi(x)}
	\end{equation}
	where $\Phi_j=\pi_j\circ\psi\circ\Phi$. Therefore the pushforward, in local coordinates, is multiplication by the Jacobian.
	
	If, in addition, $\Phi$ is a diffeomorphism, then we can extend the pushforward $\Phi_*$ to act on vector fields, and, indeed, operators $T:C^\infty(M)\to C^\infty(M)$ in general. Suppose that $T$ is a mapping $:C^\infty(M)\to C^\infty(M)$. Then the \textit{pushforward} of $T$ along $\Phi$ is the mapping $\Phi_*(T):C^\infty(M)\to C^\infty(M)$ defined by
	\begin{equation}\label{pushforward}
		\Phi_*(T)f=T(f\circ\Phi)\circ\Phi^{-1}
	\end{equation}
	for all $f\in C^\infty(M)$.
	
	\begin{remark}\label{pushforwardoperators}
		Suppose now that $\Omega\subset\GG$ is an open set. Then, due to the Haar measure on $\GG$, we can consider the space $\D'(\Omega)=C_c^\infty(\Omega)^*$ of distributions on $\Omega$, defined on compactly supported smooth functions. A vector $X\in T_x\Omega$ ($x\in\Omega$), considered as a linear functional on $C^\infty(\Omega)$, extends (e.g., by continuity) from $C^\infty_c(\Omega)$ to $\D'(\Omega)$. Similarly, a smooth vector field $X\in\Gamma(T\Omega)$ gives rise to an operator, also denoted $X$, taking $\D'(\Omega)$ to itself.
	\end{remark}
	\begin{theorem}\label{Vdef}
		Let $\Omega,\Omega'\subset\GG$ be open sets and let $\Phi:\Omega\to\Omega'$ be a diffeomorphism. Then there is a unique continuous operator $V_\Phi:\D'(\Omega')\to\D'(\Omega)$ such that
		$$V_\Phi(f)=f\circ\Phi$$
		for all $f\in C^\infty_c(\Omega')$.
	\end{theorem}
	\begin{proof}
		See, e.g., \cite[Theorem 6.1.2]{hormander}.
	\end{proof}
	We will need to consider pushforwards of operators which are defined on subspaces of $\D'$, e.g. $L_2$, or $W^{1,p}$, which do not extend continuously to $\D'$. Therefore we allow the domain of definition of operators to be some subset $\Xi\subset\D'$. The appropriate domain for their pushforwards is then $$\Lambda=\{f\in\D'(\Omega'):\exists\; g\in \Xi, f=V_\Phi(g)\}=V_{\Phi}^{-1}(\Xi)$$
	\begin{definition}\label{pushforwarddef}
		Let $\Omega,\Omega'\subset\GG$ be open, let $\Phi:\Omega\to\Omega'$ be a diffeomorphism, let $\Xi\subset\D'(\Omega)$ be some set of distributions, let $\Lambda=V_{\Phi}^{-1}(\Xi)$ and let $T:\Xi\to \Xi$ be a mapping. Then \textit{the pushforward of $T$ along $\Phi$} is the mapping $\Phi_*(T): \Lambda\to \Lambda$ defined by
		\begin{equation}\label{simple}
			\Phi_*(T)=V_\Phi^{-1} T V_\Phi.
		\end{equation}
	\end{definition}
	Let us now collect some useful facts about pushforwards.
	\begin{proposition}\label{pffacts}
		Let $\Omega,\Omega'\subset\GG$ be open, let $\Phi:\Omega\to\Omega'$ be a diffeomorphism, let $\Xi\subset\D'(\Omega)$ be some set of distributions and let $\Lambda=V_{\Phi}^{-1}(\Xi)$. Let $T,S:\Xi\to \Xi$, let $X\in\Gamma(T\Omega)$, let $g\in C^\infty(\Omega)$ and let $f:\Omega\to\GG$ be a function such that $M_f$ is defined on $\Xi$. Then
		\begin{enumerate}
			\item\label{pfmultiplier} $\Phi_*(M_f)=M_{f\circ\Phi^{-1}}$;
			\item\label{pfhomo} $\Phi_*(TS)=\Phi_*(T)\Phi_*(S)$;
			\item\label{pffunc} $\Phi_*(gX)=(g\circ\Phi^{-1})\Phi_*(X)$;
			\item\label{pfbracket} $\Phi_*([T, S])=[\Phi_*(T), \Phi_*(S)]$;
			\item\label{pfvectorcoord} \itemcase{i}f $\varphi$ and $\psi$ are coordinate charts on $\Omega$ and $\Omega'$, respectively, and $X=\sum_i\alpha_i\partial/\partial\varphi_i$,
			$$\Phi_*(X)=\sum_{i,j=1}^d \biggl(\biggl(\alpha_i\frac{\partial}{\partial\varphi_i}\Phi_j\biggr)\circ \Phi^{-1}\biggr)\frac{\partial}{\partial \psi_j}=\sum_{j=1}^d \Bigl((X\Phi_j)\circ \Phi^{-1}\Bigr)\frac{\partial}{\partial \psi_j};$$
			\item\label{pfleftinvariant} \itemcase{i}f $\Omega=\GG$, then $X$ is left-invariant if and only if $(L_x)_*X=X$ for all $x\in\GG$;
			\item\label{pfleftinvariantauto} \itemcase{i}f $\Omega=\GG$ and $\Phi$ is an automorphism of $\GG$, then $T$ is left-invariant, in the sense that $(L_x)_*T=T$ for all $x\in\GG$, if and only if $\Phi_*(T)$ is.
		\end{enumerate}
	\end{proposition}
	\begin{proof}
		(\ref{pfhomo}) is immediate from the definition, and (\ref{pfbracket}) is immediate from (\ref{pfhomo}). For (\ref{pfmultiplier}), suppose that $M_f$ is defined on $\Xi$. That is, for every $g\in \Xi$, the mapping $M_fg:C_c^\infty(\Omega)\to\CC$ defined by $M_fg(h)=g(fh)$ is a continuous linear functional and is contained in $\Xi$. Let $g\in \Lambda$ and $h\in C^\infty_c(\Omega')$. Then
		$$(V_\Phi^{-1}M_fV_\Phi g)(h)=(M_fV_\Phi g)(h\circ \Phi)=(V_\Phi g)(f(h\circ \Phi))=g(f(h\circ \Phi))\circ\Phi^{-1}=M_{V_\Phi^{-1}f}g(h).$$
		so that $V_\Phi^{-1}M_fV_\Phi g=M_{V_\Phi^{-1}f}$ as desired.
		(\ref{pffunc}) is the combination of (\ref{pfhomo}) and \ref{pfmultiplier} applied to $M_f$ and $X$, with $\Xi=C^\infty(\Omega)$, and (\ref{pfvectorcoord}) results from the application of (\ref{pfptcoord}) at all points, plus linearity and (\ref{pffunc}). For (\ref{pfleftinvariant}), by definition $X$ is left-invariant if and only if $X_x=((L_x)_*)_xX_o$ for all $x\in\GG$. This is equivalent to that $X_x=X_o(f\circ L_x)$ for all $x\in\GG$ and $f\in C^\infty(\GG)$, or equivalently that $X(f)(x)=X(f\circ L_x)(L_x\rp(x))$ for all $f\in C^\infty(\GG)$, which is of course equivalent to that $(L_x)_*(X)=X$. For (\ref{pfleftinvariantauto}), observe that, if $x\in \GG$, $y=\Phi^{-1}(x)$ and $(L_y)_*(T)=T$, then
		$$\Phi_*(T)=\Phi_*((L_y)_*(T))=(\Phi_*\circ (L_y)_*)(T)=((L_{y})_*\circ \Phi_*)(T)=\Phi_*(T)=(L_{x})_*(\Phi_*T)$$
		hence $\Phi_*T$ is left-invariant. To obtain the converse, apply the forward implication with $\Phi^{-1}$ and $\Phi_*(T)$.
		
	\end{proof}
	
	\begin{proposition}\label{automorphism push-formawrd}
		Let $A\in\aut(\GG)$ be an automorphism, with matrix $(A_{ij})$ with respect to \named coordinates, and let $1\leq i\leq d$. Then
		$$A_*(Z_i)=\sum_{j=1}^d A_{ij}Z_j.$$
	\end{proposition}
	\begin{proof}
		By Proposition \ref{pffacts} (\ref{pfleftinvariantauto}), the vector field $A_*(Z_i)$ is left-invariant; write $\displaystyle A_*(Z_i)=\sum_j a_{ij}Z_j$. Then
		$$A_*(Z_i)|_o=\sum_{j=1}^d a_{ij}Z_j|_o=\sum_{j=1}^d a_{ij}\frac{\partial}{\partial x_j}\bigg|_o$$
		However, by (\ref{pfptcoord}), 
		$$A_*(Z_i)|_o=(A_*)_o(Z_i|_o)=(A_*)_o\biggl(\frac{\partial}{\partial x_i}\bigg|_o\biggr)=\sum_{j=1}^d\frac{\partial}{\partial x_i} A_j(o)\frac{\partial}{\partial x_j}\bigg|_{A(o)}=\sum_{j=1}^dA_{ji}\frac{\partial}{\partial x_j}\bigg|_{o}.$$
		Therefore, $a_{ij}=A_{ij}$ for all $1\leq j\leq d$.
	\end{proof}

	\subsubsection{$\GG$-Diffeomorphisms and $\GG$-Atlases.}
	
	\begin{definition}
		Let $M$ and $N$ be filtered manifolds, with filtrations $(H^j)_{j\in\NN}$ and $(K^j)_{j\in\NN}$, respectively. A smooth map $\Phi:M\to N$ is said to be \textit{filtered} if it respects the filtrations on $M$ and $N$, in the sense that
		$$\Phi_*(H^j)\subset K^j$$
		for all $j\in\N$.
	\end{definition}
	
	\begin{remark}\label{firstlayer}
		If $M$ is an equiregular filtered manifold and $\Phi$ is a diffeomorphism, then for $\Phi$ to be filtered, it is sufficient that $\Phi_*(H_1)\subset K^1$, since the pushforward of vector fields along a diffeomorphism respects Lie bracket (\ref{pfhomo}), hence if $x\in M$ and $X\in H^{n+1}_x$, then $\displaystyle X=V+\sum_j[Y_j, Z_j]_x$ for some $Y_i\in H^1_x$ and $Z_i\in\Gamma(H^n)$ and so
		$$(\Phi_*)_x(X)=(\Phi_*)_x(V)+\sum_{j=1}^k((\Phi_*)([Y_i, Z_i]))_x=(\Phi_*)_x(V)+\sum_{j=1}^k[\Phi_*(Y_i), \Phi_*(Z_i)]_x$$
		thus if $(\Phi_*)_x(Z_x)\in\Gamma(K^n)$ for every $Z\in\Gamma(H^n)$, we get $(\Phi_*)_x(X)\in K^{n+1}_x$, and hence $(\Phi_*)(Z')\in\Gamma(K^{n+1})$ for every $Z'\in\Gamma(H^{n+1})$, so the result follows by induction.
	\end{remark}
	
	\begin{theorem}\label{vphiv}
		Let $\Omega, \Omega'\subset\GG$ be open and $\Phi:\Omega\to\Omega'$ a diffeomorphism. The following are equivalent. 
		\begin{enumerate}[\rm(1)]
			\item $\Phi$ is filtered for the standard filtration on $\GG$ (restricted to $\Omega$ and $\Omega'$);\vspace{0.1cm}
			\item\label{David's G-diffeomorphism} \itemcase{f}or all $x\in\Omega$, one has
			$$(\Phi_x)_*(\gra{1}{x})\subset\gra{1}{\Phi(x)};$$
			\item \itemcase{f}or all $1\leq k\leq n_1$.
			$$\Phi_*(X_k)\in\Gamma(\gra{1}{{}});$$
			\item \itemcase{t}here are smooth functions $p_{ij}:\Omega'\to\RR$, for $1\leq i,j\leq n_1$, such that, for all $1\leq k\leq n_1$,
			$$\Phi_*(X_k)=\sum_{j=1}^{n_1}p_{kj}X_j;$$
			\item\label{Dima's G-diffeomorphism} \itemcase{f}or all $1\leq k\leq n_1$,
			$$\Phi_*(X_k)=\sum_{j=1}^{n_1} (X_k\Phi_j\circ \Phi^{-1})X_j;$$
			\item\label{Fulin's G-diffeomorphism}\itemcase{f}or all $x\in\Omega$,
			$$J_{\Phi}(x)(F(x))\subset F(\Phi(x)),$$
			where, for all $x\in\GG$,
			$$F(x):=\bigcap_{i=n_1+1}^d \biggl\{e_i-\sum_{j=1}^{n_1} \alpha_{ij}(x)e_j\biggr\}^\perp\subset\RR^d,$$
			and the $\smash{\alpha_{ij}}$ are the functions appearing in Remark \ref{basis-rep}.
		\end{enumerate}
	\end{theorem}
	The proof of Theorem \ref{vphiv} is  deferred to Section \ref{vphivsec}. This theorem motivates the following definition, and provides a number of approaches to working with these maps in \named coordinates on $\GG$ and in terms of the standard filtration on $\GG$.

	\begin{definition}\label{G-diffeomorphism}
		Let $\Omega,\Omega'\subset\GG$ be open sets and $\Phi:\Omega\to\Omega'$ a diffeomorphism. Then $\Phi$ is said to be a \textit{$\GG$-diffeomorphism} if it satisfies the equivalent conditions of Theorem \ref{vphiv}.
	\end{definition}
	
	In the case $\GG=\HH^d$, by Theorem \ref{vphiv}, $F(x)$ coincides with the $F(x)^\perp$ defined in \cite[Definition 1.3.5]{DAO4}. Therefore, the notion of $\HH^d$-diffeomorphism coincides with that of Heisenberg diffeomorphism, which plays an central role in establishing principal symbol and Connes' trace theorem on compact contact manifolds.\
	
	We will also have need for the related concept of \textit{admissibility}.
	
	\begin{definition}\label{admissibledef}
		Let $\Omega, \Omega'\subset\GG$ be open and let $\Phi:\Omega\to\Omega'$ be a diffeomorphism. Then $\Phi$ is said to be \textit{admissible} if, for each $x\in\GG$, there is an automorphism $H^\Phi(x)\in\at$ whose matrix $\smash{(H^\Phi_{ij}(x))_{ij=1}^d}$ in \named coordinates on $\GG$ satisfies
		$$H^\Phi_{ij}(x)=X_j\Phi_i(x)$$
		for all $1\leq i,j\leq n_1$. 
	\end{definition}
	\begin{remark}
		Since a strata preserving automorphism of $\GG$ (equivalently $\g$) is determined by its action on the first stratum $\g_1$, the $H^\Phi$ appearing in Definition \ref{admissibledef} above is uniquely determined, should it exist.
	\end{remark}
	
	\begin{theorem}\label{admissibility}
		Let $\Omega, \Omega'\subset\GG$ be open, let $x\in\Omega$ and let $\Phi:\Omega\to\Omega'$ be a $\GG$-diffeomorphism. Then $\Phi$ is admissible; indeed, the mapping on $\GG$ given by the matrix $(Z_j\Phi_i(x))_{ij=1}^d$ is a strata preserving automorphism.
	\end{theorem}
	
	\begin{proof}
		In \cite{Po19} it was shown that any filtered map $\Phi$ between filtered manifolds defines, at each point $x$, a Lie algebra homomorphism between the \textit{tangent Lie algebra} at $x$ and that at $\Phi(x)$, whose construction can be found in the same paper, and that these homomorphisms are isomorphisms if $\Phi$ is a diffeomorphism. In the case that the manifolds are open subsets of $\GG$, these tangent Lie algebras are canonically isomorphic to $\g$. It is an easy exercise to show that the Lie algebra automorphisms defined in this way are given, with respect to the \named coordinates, by the matrices $(Z_j\Phi_i(x))_{i,j}^d$. Finally, since any Lie algebra automorphism of $\g$ is identically a Lie group automorphism of $\GG$, the result follows.
	\end{proof}

	\begin{definition}\label{G-atlas}
		Let $M$ be a smooth manifold. A \textit{$\GG$-atlas on $M$} is a collection $\{(\U_{i},h_{i})\}_{i\in\II}$ of pairs $(\U_i, h_i)$ whose first entry $\U_i$ is an open subset of $M$ and whose second entry $h_i$ is a map $h_i:\U_i\to\GG$, and
		\begin{enumerate}
			\item \itemcase{t}he $\{\U_i\}_{i\in\II}$ cover $M$;
			\item \itemcase{f}or each $i\in\II$, $h_i$ is a diffeomorphism onto its image, and
			\item \itemcase{t}he transition mappings $h_j\circ h_i^{-1}:h_i(\U_i\cap \U_j)\to h_j(\U_i\cap \U_j)$ are $\GG$-diffeomorphisms.
		\end{enumerate}
	\end{definition}
	
	\begin{proposition}\label{inducedfiltration}
		Let $M$ be a smooth manifold and let $\{(\U_{i},h_{i})\}_{i\in\II}$ be a $\GG$-atlas on $M$. Then there is a unique bracket generated filtration $H$ on $M$ making $M$ into an equiregular filtered manifold such that, for each $i\in\II$, $h_i$ is a filtered diffeomorphism.
	\end{proposition}
	\begin{proof}
		Clear
	\end{proof}
	
	If a $\GG$-atlas induces the filtration of a filtered manifold $(M, H)$ in the sense of Proposition \ref{inducedfiltration}, we say that the atlas is \textit{compatible} with the filtration. Clearly, a filtered manifold with a compatible $\GG$-atlas is equiregular. It is natural to ask whether every equiregular filtered manifold admits a compatible $\GG$-atlas. By a famous theorem of Darboux, it is true for the special case of contact manifolds.
	
	\begin{theorem}[Darboux Theorem]
		Let $M$ be a filtered manifold of odd dimension such that $H^1$ has codimension 1, and $H^2=TM$. Then $M$ has a compatible $\HH^n$-atlas, where $\HH^n$ is the Heisenberg group with the same dimension as $M$.
	\end{theorem}
	
	A strengthening of the Darboux Theorem was given in \cite[Corollary 3.6.1]{morimoto1993geometric}, however, in full generality of equiregular filtered manifolds, it is false.
	
	\begin{definition}\label{G-filtered manifold}
		Let $M$ be a filtered manifold. Then, following \cite{MK2022}, we say that $M$ has the \textit{Darboux property} if there is a stratified group $\GG$ such that $M$ admits a compatible $\GG$-atlas. In this case we will also say that $M$ is a \textit{$\GG$-filtered manifold}.
	\end{definition}
	
	Cartan's famous ``two rolling balls'' example \cite{balls} is, to the best of our knowledge, the only known example of equiregular filtered manifolds for which the Darboux property fails (unfortunately there is no accessible proof of this failing in the literature at this time), and so a typical equiregular filtered manifold one comes across in practice generally satisfies the Darboux property.

	\subsection{Main results}\label{main results}
	
	\subsubsection{Principal Symbol On Groups}
	
	
	\begin{definition}
		An automorphism $A\in\aut(\GG)$ is said to be \textit{strata preserving} if it respects the filtration on $\GG$. The set of strata-preserving automorphisms on $\GG$ is denoted by $\at$.
	\end{definition}
	
	\begin{definition}\label{XALARA}
		Let $A\in\at$ and $1\leq k\leq n_1$. Define
		$$X_k^A:=A_*(X_k),\quad L_A:=A_*(\Delta)\quad\text{and}\quad \quasiriesz A k :=A_*(\R_k),$$
		where $\Delta$ is the \textit{sub-Laplacian} $\displaystyle{\Delta=-\sum_{k=1}^{n_{1}}X_{k}^{2}}$ and $R_k$ is the \textit{Riesz transform} $\R_k=X_j\Delta^{-1/2}$. The operators $\quasiriesz A k$ are called the \textit{quasi-Riesz transforms}.
	\end{definition}

%
	
	By definition of $A_*$, we have $X^A_k=V_A^{-1}X_k V_A$, $L_{A}=V_{A}\rp \Delta V_{A}$ and $R^{A}_{k}=V_{A}\rp R_{k}V_{A}$. By Proposition \ref{automorphism push-formawrd}, since $A$ is strata preserving, we have the identity
	$$X^A_k=\sum_{j=1}^{n_1} A_{kj} X_j$$
	Additionally, $\displaystyle L_A=-\sum_{k=1}^{n_{1}}(X^A_k)^{2}$ and $\quasiriesz{A}{k}=X^A_k L_A^{-1/2}$.
	
	\begin{remark}\label{left invariance}
		Since the vector fields $X_k$ are left-invariant, by Proposition \ref{pffacts} (\ref{pfleftinvariantauto}), each of the operators defined in Definition \ref{XALARA} are also left-invariant, and hence so is, e.g., any strong limit thereof.
	\end{remark}

	Now we introduce three $C^{*}$-algebras on $\GG$, to act as generalisations of their counterparts in the Euclidean case (\!\cite[Definition 1.1]{DAO1}, \cite[Definition 1.1]{DAO3}) and the Heisenberg case (\!\cite[Definition 1.3.1]{DAO4}).
	\begin{definition}
		Let 
		\begin{align*}
			\A_{1}=\CC+C_0(\GG),
			\quad
			\po(f)=M_{f}\quad\mbox{for}\quad f\in\A_{1};
		\end{align*}
		and 
		\begin{align*}
			\A_{2}=C^{*}\Big(\{\quasiriesz A k:1\leq k\leq n_1, A\in \at\}\Big),
			\quad
			\pt(g)=g\quad\mbox{for}\quad g\in\A_{2}.
		\end{align*}
		Let $\Pi$ be the $C^{*}$-subalgebra in $\B(\lp{2})$ generated by $\po(\A_{1})$ and $\pt(\A_{2})$, i.e.
		$\Pi=C^{*}\Big(\{\po(\A_{1}),\pt(\A_{2})\}\Big)$.
		
	\end{definition}
 \begin{remark}
	Note that different bases  on $\g$ are isomorphic. This results $\A_2$, as defined above, can also be defined by other basis, i.e. $\A_2$ is independent of the choice of the basis $\set{X_k}_{k=1}^{n_{1}}$.
 \end{remark}

	Let $\A_{1}\otimes_{\rm min}\A_{2}$ be the minimal tensor product which is defined as the completion of the algebraic tensor product $\A_{1}\otimes\A_{2}$ with respect to the norm induced by $\B(\lp{2}\otimes_{2}\lp{2})$ via the embedding
	$$\A_{1}\otimes\A_{2}\subset\B(\lp{2}\otimes_{2}\lp{2}).$$
	For details on minimal tensor products, see \cite{T-book2002}.
	In the sense of isometric isomorphism, we have
	$$\A_{1}\otimes_{\rm min}\A_{2}\cong\A_{2}\otimes_{\rm min}\A_{1}.$$
	The following theorem is the existence of principal symbol on $\GG$ which arises as an application of Theorem 3.3 in \cite{DAO2} (see also \cite[Theorem 5.2.6]{LMSZv2}).
	\begin{theorem}\label{symbol}
		Let $\GG$ be a stratified Lie group.
		There is a surjective $*$-homomorphism $\sym:\Pi\rightarrow\A_{1}\otimes_{\rm min}\A_{2}$ such that
		\begin{align*}
			\sym(\po(f))=f\otimes1,\quad\sym(\pt(g))=1\otimes g,\quad f\in\A_{1},g\in\A_{2}.
		\end{align*}
		Furthermore,
		\begin{align*}
			\K(\lp{2})=\ker(\sym)\subset\Pi.
		\end{align*}
	\end{theorem}
	The $*$-homomorphism $\sym$ is called the principal symbol mapping. It maps a non-commutative algebra to a commutative algebra. In (non-commutative) Euclidean space and Heisenberg groups, it plays important role in the theory of pseudodifferential operators and Connes' trace theorem. In this note, it serves as the grounding of principal symbol mapping on filtered manifolds.

	\subsubsection{Principal Symbol on Manifolds}\label{symbol on manifold-section}
	\begin{definition}
		Let $M$ be a smooth manifold and let $\Borel$ denote the Borel $\sigma$-algebra of $M$. We say that a measure $\nu$ on $\Borel$ is a \textit{smooth positive density} if, for each smooth chart $(\U,h)$ on $M$, the measure $\nu\circ h^{-1}$ is absolutely continuous with respect to Lebesgue measure on $h(\U)$, and the Radon-Nikodym derivative of $\nu\circ h^{-1}$ with respect to Lebesgue measure on $h(\U)$ is smooth and strictly positive.
	\end{definition}
	
	For a $d$-dimensional manifold $M$ with a smooth positive density $\nu$ and with an atlas $\{(\U_{i},h_{i})\}_{i\in\II}$, if an operator $T\in\B(L_{2}(M,\nu))$ is compactly supported in $\U_{i}$ for some $i\in\II$ (i.e. there exists $\psi\in C_{c}^{\infty}(\U_{i})$ such that $T=TM_{\psi}=M_{\psi}T$), then, by pre-composing with $h_{i}$, we can transfer $T$ to an operator on $L_{2}(\RR^{d})$. We do this with a $\GG$-atlas, obtaining operators on $L_2(\GG)$ from compactly supported operators on a $\GG$-filtered manifold.
	
	Let us give the definition of the $C^{*}$-algebra $\Pi_{M}$, the domain of our principal symbol mapping.
	\begin{definition}\label{Pi-X}
		Let $M$ be a $\GG$-filtered manifold, with $\GG$-atlas $\{(\U_{i},h_{i})\}_{i\in\II}$ and let $\nu$ be a smooth positive density on $M$. Then $\Pi_{M}$ is defined to be the set of operators $T\in\B(L_{2}(M,\nu))$ satisfying
		\begin{enumerate}[\rm(1)]
			\item\label{Pi tranfer to group} \itemcase{f}or each $i\in\II$ and each $\phi\in C_{c}(\U_{i})$, the operator $M_{\phi}TM_{\phi}\in L_{2}(\GG)$ belongs to $\Pi$;
			\item \itemcase{f}or each $\phi\in C_{c}(M)$, the operator $[T,M_{\phi}]$ is compact on $L_{2}(M,\nu)$;
			\item \itemcase{t}here exists a sequence $(\phi_k)_{k\in\NN}$ in $C_{c}(M)$ such that $M_{\phi_k}T\to T$ uniformly in $\B(L_{2}(M,\nu))$.
		\end{enumerate}
	\end{definition}

	
	Now we are ready to describe the codomain of our principal symbol mapping on $\GG$-filtered manifolds.
	For manifold $M$ with an atlas $\{(\U_{i},h_{i})\}_{i\in\II}$ and a unital $C^{*}$-algebra $\M$, let $\omega=\{\pi_{i,j}\}_{i,j\in\II}$ be a family of continuous mappings $\pi_{i,j}:\U_{i}\cap\U_{j}\rightarrow\aut(\M)$ defined on each nonempty overlap $\U_{i}\cap\U_{j}$ such that the following cocycle conditions hold:
	\begin{align}\label{cocycle-two} 
		\pi_{i,j}(x)\pi_{j,i}(x)=\id\quad \mbox{for}\quad x\in\U_{i}\cap\U_{j};
	\end{align}
	\begin{align}\label{cocycle-three} 
		\pi_{i,j}(x)\pi_{j,k}(x)\pi_{k,i}(x)=\id\quad \mbox{for}\quad x\in\U_{i}\cap\U_{j}\cap\U_{k}.
	\end{align}
	The triple $(M,\M,\omega)$ is called a bundle of $C^{*}$-algebras. More precisely, it is a bundle with fiber $\M$ and structure group $\aut(\M)$ that is trivializable over $\{\U_{i}\}_{i\in\II}$ with transition functions $\{\pi_{i,j}\}_{i,j\in\II}$.

	When $E=(M,\M,\omega)$ is a bundle of $C^{*}$-algebras, the space $C_{b}(E)$ of bounded continuous sections of this bundle is the set of all families $F=\{F_{i}\}_{i\in\II}$ satisfying the following three laws
	\begin{enumerate}[\rm(i)]
		\item \itemcase{f}or each $i\in\II$, we have $F_{i}\in C(\U_{i},\M)$;
		\item \itemcase{f}or each $i,j\in\II$, we have $F_{j}=\pi_{i,j}(F_{j})$ on $\U_{i}\cap\U_{j}$;
		\item \itemcase{t}he norm $\displaystyle\Norm{F}_{C_{b}(E)}=\sup_{i\in\II}\Norm{F_{i}}_{C(\U_{i},\M)}$ is finite.
	\end{enumerate}
	By the argument in \cite{DAO4}, if $E$ be a bundle of $C^{*}$-algebras, then $C_{b}(E)$ is a $C^{*}$-algebra.
	In addition, for $C^{*}$-algebra bundle $E=(M,\M,\omega)$ and each $i\in\II$, there is a natural embedding 
	\begin{align}\label{natural embedding}
		\Theta_{i}:C_{c}(\U_{i},\M)\rightarrow C_{b}(E).
	\end{align} 
	Given $f\in C_{c}(\U_{i},\M)$, the corresponding family $F=\{F_{j}\}_{j\in\II}$ is defined as follows. Let $F_{j}\in C(\U_{j},\M)$ be the continuous extension of $\pi_{j,i}(f|_{\U_{i}\cap\U_{j}})$ to $\U_{j}$ by zero if $\U_{i}\cap\U_{j}\neq\emptyset$ and be 0 if $\U_{i}\cap\U_{j}=\emptyset$. The each $F_{j}$ is well-defined and continuous, and the family $F=\{F_{j}\}_{j\in\II}$ satisfies the compatibility condition. 

	\begin{definition}
		Let $\Phi$ be an admissible diffeomorphism on $\GG$. Define
		\begin{align*}
			\pi_{H^{\Phi}}(\omega)(x)=V_{H^{\Phi}(x)}\rp\omega(x)V_{H^{\Phi}(x)}
			,\quad\omega\in \lp{\infty}\overline{\otimes}\B(\lp{2}).
		\end{align*}
	\end{definition}
	It is easy to check that
	\begin{align}\label{transfer}
		\pi_{H^{\Phi}}(f\otimes\one)=f\otimes\one,\quad
		\pi_{H^{\Phi}}(\one\otimes R_{k})\circ \Phi\rp=\one\otimes \quasiriesz{H^{\Phi}}{k}.
	\end{align}
	For our purpose, we consider a particular $C^{*}$-algebra bundle $E_{hom}=(M,\A_{2},\omega)$ with 
	$\omega=\{\pi_{i,j}\}_{i,j\in\II}$  given as
	\begin{align*}
		\pi_{i,j}(x)=\pi_{H^{\Phi_{i,j}}(h_{i}(x))},\quad i,j\in\II \quad such\quad that\quad x\in\U_{i}\cap\U_{j}\neq\emptyset.
	\end{align*}
	\begin{definition}
		Let $M$ be a $\GG$-filtered manifold with atlas $\{(\U_{i},h_{i})\}_{i\in\II}$. The codomain of the principal mapping is the algebra  $C_{b}(E_{hom})$ of bounded continuous sections of $E_{hom}$. 
	\end{definition}

	The main result of this paper is the following theorem.
	\begin{theorem}\label{symbol on manifold}
		Let $M$ be a $\GG$-filtered manifold, with $\GG$-atlas $\{(\U_{i},h_{i})\}_{i\in\II}$ and let $\nu$ be a smooth positive density on $M$. Then
		\begin{enumerate}[\rm(1)]
			\item $\Pi_{M}$ is a $C^{*}$-algebra which contains $\K(L_{2}(M,\nu))$;
			\item \itemcase{t}here exists a surjective $*$-homomorphism $\sym_{M}:\Pi_{M}\to C_{b}(E_{hom})$ such that $$\ker(\sym_{M})=\K(L_{2}(M,\nu)).$$
		\end{enumerate}
	\end{theorem}
	The global principal symbol mapping $\sym_{M}$ is compatible with the local one introduced in Theorem \ref{symbol} in the following sense. If $T\in\Pi_{M}$ is compactly supported in $\U_{i}$ for some $i\in\II$, then the operator $T$ transferred to an operator on $L_{2}(\GG)$ belongs to $\Pi$ and $$\sym_{M}(T)=\Theta_{i}(\sym(T)).$$

	This paper is organized as follows. 
	In Section \ref{principal-symbol}, we show Theorem \ref{symbol} the existence of principal symbol mapping on $\GG$. 
	In Section \ref{quasi-Riesz-transform} and Section \ref{Equivariance-of-PS}, we investigate the properties of quasi-Riesz transforms and the equivariance of principal symbol mapping on $\GG$.
	For our purpose, a globalisation theorem for the existence of $*$-homomorphism between two $C^{*}$-algebras on manifolds is established in Section \ref{Globalisation theorem}. 
	Lastly, in Section \ref{symbol on Carnot manifolds}, we give the proof of Theorem \ref{symbol on manifold} the existence of principal symbol mapping on filtered manifold.

	Throughout the paper, we use $A\lesssim B$ to denote the statement that $A\leq CB$ for some positive constant $C$ independent of $A$ and $B$, and $A\simeq B$ to denote the statement that $A\lesssim B$ and $B\lesssim A$.

	\section{Principal symbol on stratified Lie groups}\label{principal-symbol}
	
	In this section, we prove Theorem \ref{symbol}: the existence of the principal symbol mapping on stratified Lie groups.
	Recall that the Riesz transforms on $\GG$ are defined as 
	$$\R_{j}=X_{j}\Delta^{-1/2},\quad j=1,2,\,\dotsc,n_{1},$$ 
	where $\Delta$ is the sub-Laplacian $\displaystyle{\Delta=-\sum_{k=1}^{n_{1}}X_{k}^{2}}.$
	The operator $\Delta$ is a 2-homogeneous linear differential operator with polynomial coefficients and it commutes with left translations. 
	The Riesz transforms are 0-homogeneous and left-translation invariant on $\GG$. 
	Moreover, the $j$-th Riesz transform $\R_{j}$ is a left convolutional singular integral operator, i.e.
	\begin{align}\label{Riesz-convolution}
		\R_{j}f(x)=\int_{\GG}K_{j}(y\rp x)f(y)\,dy\quad{\rm for}\quad f\in C_{c}^{\infty}(\GG),
	\end{align}
	where
	\begin{align*}
		K_{j}(x)=\frac{1}{\sqrt{\pi}}\int_{0}^{\infty}t^{-Q/2-1}(X_{j}h_{1})(\delta_{t^{-1/2}} (x))\,dt.
	\end{align*} 
	Here $\delta_{r}$ $(r>0)$ is the dilation on $\GG$ and $h_{t}(\cdot)=u(t,\cdot)$ is the normalised $C^{\infty}$ solution of the heat equation $(\p/\p t+\Delta)u=0$ with initial condition $u(0,\cdot)=\delta$ (see  \cite[chapter 1]{FS1982}, \cite[IV.4]{VSC1992}). 
	By \cite{FS1982}, $K_{j}$ is smooth away from the origin $o$ and is homogeneous of degree $-Q$, i.e.
	\begin{align}\label{Riesz-homogeneous}
		K_{j}(rx)=r^{-Q}K_{j}(x) \quad \mbox{for} \quad x\neq o.
	\end{align}

	\subsection{Verification of conditions}
	
	Let $\lambda:\GG\rightarrow\B(\lp{2})$ be the left regular representation which is defined by setting
	\begin{align*}
		\lambda_{x}f(y)=f(x\rp y),\quad f\in\lp{2},\quad x, y\in\GG.
	\end{align*}
	It is easy to see that $\lambda_{x}$ is a unitary operator on $\lp{2}$ and $\lambda_{x}\rp=\lambda_{x\rp}$.

	
	By Remark \ref{left invariance}, every element of $\A_2$ is left-invariant. The dilations $\sg{r}$ ($r>0$) on $\lp{2}$ are defined by
	\begin{align*}
		\sg{r}(f):=r^{Q/2}(f\circ \dl{r})
	\end{align*}
	Basically, $\sg{r}$ is unitary on $\lp{2}$ and $\sg{r}\rp=\sg{r\rp}$. 
	Additionally, all $\R_{j}$ ($j=1,2,\,\dotsc,n_{1}$) are dilation-invariant. 
	Indeed, 
	by \eqref{Riesz-homogeneous}, we obtain
	\begin{align*}
		\sg{r} \R_{j} \sg{r\rp}(f)(x)=\int_{\GG}f(\delta_{r\rp}(y))K_{j}(y\rp\delta_{r}(x))\,dy
		=\R_{j}(f)(x),\quad f\in\lp{2}.
	\end{align*}
	This implies that each element in $\A_{2}$ is dilation-invariant.

	\medskip
	

	Recall some notations in \cite{Dixmier-book}.
	Denote by $(\pi,\H)$ a representation of $\A$ on a Hilbert space $\H$. A closed subspace $E$ of $\H$ is called an invariant subspace of representation $(\pi,\H)$ if $\pi(x)E\subset E$ for all $x\in\A$. If $(\pi,\H)$ has no invariant subspace other than $\H$ and $\{0\}$, then it is said to be irreducible.
	
	Let the set $S\subset\B(\H)$ and $S'$ be the commutant of $S$ in $\B(\H)$. 
	The following lemmas are well-known. 
	\begin{lemma}\label{id-irreducible}
		\cite[Proposition II.6.1.8]{Blackadar-book}
		Let $\pi$ be a representation of a $C^{*}$-algebra $\A$ on a Hilbert space $\H$. Then
		$\pi$ is irreducible if and only if $\pi(\A)'=\CC\one$.
	\end{lemma}

	\begin{lemma}\label{irreducible}
		\cite[Corollary 4.1.10]{Dixmier-book}
		Let $\A$ be a $C^{*}$-algebra and its representation $\pi:\A\rightarrow\B(\H)$ be irreducible. 
		One of the following mutually exclusive options holds:
		\begin{enumerate}[\rm(1)]
			\item $\pi(\A)$ does not contain any compact operator (except for 0);
			
			\item $\pi(\A)$ contains all compact operators. 
		\end{enumerate}
	\end{lemma}

	Let us continue the process now.
	\begin{lemma}\label{containalg}
		The algebra $\K(\lp{2})$ is contained in $\Pi$.
	\end{lemma}
	\begin{proof}
		Embed $\lp{\infty}$ in $\B(\lp{2})$ as a commutative $C^{*}$-subalgebra. It follows that
		\begin{align*}
			\po(\A_{1})'=\po(\lp{\infty})'=\po(\lp{\infty}).
		\end{align*}
		Since $\po(\A_{1})\subset\Pi$, it follows that $\Pi'\subset\po(\A_{1})'$. 
		Therefore, for $A\in\Pi'$, there is $f\in\lp{\infty}$ such that $A=\po(f)$. 
		By \cite[Thm 1.2]{DLLW2019}, it concludes that
		\begin{align*}
			\Norm{f}_{\BMO}\lesssim\Norm{[\po(f),\R_{k}]}_{\B(\lp{2})}=\Norm{[A,\R_{k}]}_{\B(\lp{2})}=0.
		\end{align*}
		Here $\Norm{\cdot}_{\BMO}$ denotes the BMO norm on $\GG$.
		This implies that $f$ is a constant and therefore $\Pi'=\CC\one$. By Lemma \ref{id-irreducible}, the representation $\id:\Pi\rightarrow\B(\lp{2})$ is irreducible.
		
		Pick a non-zero function $f\in\cci$ such that $[\po(f),\R_{k}]$ is non-zero. By \cite{CDLW2019}, one finds that $[\po(f),\R_{k}]$ is compact. Note that $[\po(f),\R_{k}]\in\Pi$. Applying Lemma \ref{irreducible} to $\H=\lp{2}$, one has $\K(\lp{2})\subset\Pi$. This is the desired result.
	\end{proof}
	
	\begin{lemma}\label{compactness}
		For $f\in\A_{1}$ and $g\in\A_{2}$, then $[\po(f),\pt(g)]\in\K(\lp{2})$.
	\end{lemma}
	\begin{proof}
		Firstly, let $f\in\cci$.
		By \cite{CDLW2019}, one finds that
		$$[\po(f),\R_{k}],\; [\po(f),\R_{k}^{*}]\in\K(\lp{2}),\quad 1\leq k\leq n_{1}.$$
		Note that 
		\begin{align*}
			V_{A}M_{f} V_{A}\rp=M_{f\circ A},\quad A\in\at.
		\end{align*}
		Then
		\begin{align*}
			[\po(f),\quasiriesz{A}{k}]=V_{A}\rp[\po(f\circ A),\R_{k}]V_{A}\in\K(\lp{2}).
		\end{align*}
		Similarly, $$[\po(f),(\quasiriesz{A}{k})^{*}]\in\K(\lp{2}).$$
		
		For $g\in\A_{2}$, pick polynomials $\{g_{m}\}_{m>0}$, generated by $\{\quasiriesz{A}{k},(\quasiriesz{A}{k})^{*},1\leq k\leq n_{1},A\in\at\}$, such that $g_{m}$ tends to $g$ in the uniform norm. Observe that $[\po(f),\pt(g_{m})]$ is a compact operator. Since
		\begin{align*}
			\Norm{[\po(f),\pt(g_{m}-g)]}_{\infty}\leq2\Norm{f}_{\infty}\Norm{g_{m}-g}_{\infty},
		\end{align*}
		the commutator $[\po(f),\pt(g)]$ is the limitation of $\set{[\po(f),\pt(g_{m})]}_{m>0}$ in the uniform norm. Hence $[\po(f),\pt(g)]$ is a compact operator.
		
		Next, if $f\in C_{0}(\GG)$ and $g\in\A_{2}$, then choose $f_{m}\in\cci$ such that $f_{m}$ tends to $f$ in the uniform norm. Since
		\begin{align*}
			\Norm{[\po(f_{m}-f),\pt(g)]}_{\infty}\leq2\Norm{f_{m}-f}_{\infty}\Norm{g}_{\infty},
		\end{align*}
		the compact operator $[\po(f_{m}),\pt(g)]$ tends to $[\po(f),\pt(g)]$ in the uniform norm. Therefore, $$[\po(f),\pt(g)]\in\K(\lp{2}).$$
		
		Finally, when $f$ is a constant and $g\in\A_{2}$, commutator $[\po(f),\pt(g)]$ is zero and the argument is always true. This completes the proof.
	\end{proof}

	\medskip

	In order to continue the investigation of the existence of principal symbol mapping, we need the following two elementary results.
	The first one is the translation property of multiplication operator,
	\begin{align}\label{multiplication-conjugate}
		\lambda_{x}\po(f)\lambda_{x}\rp=\po(\lambda_{x}f),\; f\in\lp{\infty}.
	\end{align}
	Indeed,
	\begin{align*}
		\lambda_{x}(M_{f}\lambda_{x}\rp h)(y)
		=(M_{f}\lambda_{x}\rp h)(x\rp y)
		=f(x\rp y) (\lambda_{x}\rp h)(x\rp y)
		=(M_{\lambda_{x}f}h)(y),\; h\in\lp{2}.
	\end{align*}
	In particular, when $f=\chi_{B(x,r)}$, one finds that
	\begin{align}\label{chibox}
		\lambda_{x}\rp\po(\chi_{B(x,r)})\lambda_{x}=\po(\chi_{B(0,r)}).
	\end{align}
	The second one is the dilation property of multiplication operator with respect to the characteristic function, i.e.
	\begin{align}\label{sgmaboo}
		\sg{r}\po(\chi_{B(0,r)})\sg{r\rp}=\po(\chi_{B(0,1)}).
	\end{align}
	Indeed, 
	\begin{align*}
		\sg{r}(M_{\chi_{B(0,r)}}\sg{r\rp}f)(x)
		=r^{Q/2}(M_{\chi_{B(0,r)}}\sg{r\rp}f)(\dl{r}(x))
		=\chi_{B(0,r)}(\dl{r}(x))f(x)
		=\chi_{B(0,1)}(x)f(x)
	\end{align*}

	\begin{lemma}\label{projection-compact}
		\cite[Lemma 4.1]{DAO2}
		Let $A\in\K(\H)$ and $\{P_{k}\}\subset\B(\H)$ be a sequence of pairwise orthogonal projections. Then $$\Norm{P_{k}A}_{\infty}\rightarrow0,\quad{\rm as}\quad k\rightarrow\infty.$$
	\end{lemma}
	
	\begin{lemma}\label{compact-zero}
		If $g\in\A_{2}$ satisfies $\po(\chi_{B(0,1)})\pt(g)\in\K(\lp{2})$, then $g=0$.
	\end{lemma}
	\begin{proof}
		Since $\po(\chi_{B(0,1)})\pt(g)$ is compact, by Lemma \ref{projection-compact}, select  pairwise disjoint balls  $\{B(x_{k},r_{k})\}_{k>0}$ containing in $B(0,1)$ so that
		\begin{align}\label{chibxko}
			\po(\chi_{B(x_{k},r_{k})})\pt(g)=\po(\chi_{B(x_{k},r_{k})})\po(\chi_{B(0,1)})\pt(g)\rightarrow0
		\end{align}
		in the uniform norm as $k\rightarrow\infty$.
		
		Recall that each element in $\A_{2}$ is left-invariant. By formula \eqref{chibox},
		\begin{align*}
			\lambda_{x_{k}}\rp \po(\chi_{B(x_{k},r_{k})})\pt(g) \lambda_{x_{k}}=\po(\chi_{B(0,r_{k})})\pt(g).
		\end{align*}
		Since $\lambda_{x_{k}}$ is unitary on $\lp{2}$, by \eqref{chibxko}, we have $$\po(\chi_{B(0,r_{k})})\pt(g)\rightarrow0,\quad\mbox{as}\quad k\rightarrow\infty,$$
		in the uniform norm.
		
		Observe that each element in $\A_{2}$ is dilation-invariant and $\sigma_{r_{k}}$ is unitary on $\lp{2}$. Thus, by formula \eqref{sgmaboo},
		\begin{align*}
			\po(\chi_{B(0,1)})\pt(g)=\sg{r_{k}}\po(\chi_{B(0,r_{k})})\pt(g)\sg{r_{k}\rp}\rightarrow0,
		\end{align*}
		in the uniform norm as $k\rightarrow\infty$. 
		Since $\po(\chi_{B(0,1)})\pt(g)$ is independent of $k$, it follows that
		\begin{align*}
			\po(\chi_{B(0,1)})\pt(g)=0.
		\end{align*}
		It immediately yields $g=0$. This completes the proof.
	\end{proof}

	\medskip

	If $f\in \mathbb{C}+C_0(\mathbb{G})$, then
	\begin{align*}
		\sg{r}(M_{\lambda_{x}\rp f}\sg{r\rp}h)(y)
		=f(x(\dl{r}(y)))h(y),\quad h\in\lp{2}.
	\end{align*}
	This yields that
	\begin{align*}
		\Norm{\sg{r}M_{\lambda_{x}\rp f}\sg{r\rp}h-f(x)h}_{\lp{2}}^{2}=\int_{\GG}|f(x(\dl{r}(y)))-f(x)|^{2}|h(y)|^{2}\,dy.
	\end{align*}
	By dominated convergence theorem, 
	\begin{align}\label{convergencef(x)}
		\sg{r}\po(\lambda_{x}\rp f)\sg{r\rp}\rightarrow f(x)\quad\mbox{as}\quad r\rightarrow0^{+},
	\end{align}
	in the sense of uniform norm.
	
	\begin{lemma}\label{primeness}
		Let $\set{f_{k}}_{k=1}^{m}\subset\A_{1}$ and $\set{g_{k}}_{k=1}^{m}\subset\A_{2}$. If
		$\smash{\displaystyle\sum_{k=1}^{m}\po(f_{k})\pt(g_{k})\in\K(\lp{2})}$,
		then $\smash{\displaystyle\sum_{k=1}^{m}f_{k}\otimes g_{k}=0}$.
	\end{lemma}
	\begin{proof} Note that each element in $\A_{2}$ is left-invariant. By the assumption of compactness on $\displaystyle\sum_{k=1}^{m}\po(f_{k})\pt(g_{k})$, formula \eqref{multiplication-conjugate} yields that
		\begin{align*}
			\sum_{k=1}^{m}\po(\lambda_{x}\rp f_{k})\pt(g_{k})=\lambda_{x}\rp\sum_{k=1}^{m}\po(f_{k})\pt(g_{k})\lambda_{x}\in\K(\lp{2}).
		\end{align*}
		Since $\sg{r}$ is unitary and $\chi_{B(0,r)}\in\lp{\infty}$, one has
		\begin{align*}
			y_{r}:=\sg{r}\po(\chi_{B(0,r)})\sum_{k=1}^{m}\po(\lambda_{x}\rp f_{k})\pt(g_{k})\sg{r\rp}\in\K(\lp{2}).
		\end{align*}
		
		Recall that each element in $\A_{2}$ is dilation-invariant. By \eqref{sgmaboo} and \eqref{convergencef(x)},
		\begin{align*}
			y_{r}=\po(\chi_{B(0,1)})\sum_{k=1}^{m}\sg{r}\po(\lambda_{x}\rp f_{k})\sg{r\rp}\pt(g_{k})
			\rightarrow\po(\chi_{B(0,1)})\pt(\sum_{k=1}^{m}f_{k}(x)g_{k}),\quad{\rm as}\quad r\rightarrow0^{+}
		\end{align*}
		where the limitation is taken in the uniform norm. Because the limit of compact operators is a compact operator, we obtain
		$$\displaystyle\po(\chi_{B(0,1)})\pt(\sum_{k=1}^{m}f_{k}(x)g_{k})\in\K(\lp{2}).$$
		By Lemma \ref{compact-zero},
		$$\sum_{k=1}^{m}f_{k}(x)g_{k}=0$$
		for {\it every} $x\in\mathbb{G}.$ This completes the proof.
	\end{proof}

	\medskip

	\subsection{Proof of Theorem \ref{symbol}}
	\ 
	\newline
	
	The following assertion can be found in \cite[Theorem 3.3]{DAO2}.
	\begin{lemma}\label{symbol-existence}
		Let $\C_{1}$ and $\C_{2}$ be $C^{*}$-algebras, and let $\alpha_{1}:\C_{1}\rightarrow\B(\H)$ and $\alpha_{2}:\C_{2}\rightarrow\B(\H)$ be representations. Denote $\Pi(\C_{1},\C_{2})$ the $C^{*}$-algebra generated by $\alpha_{1}(\C_{1})$ and $\alpha_{2}(\C_{2})$. Suppose that
		
		\begin{enumerate}[\rm(1)]
			\item $\C_{1}$, $\C_{2}$ are unital and $\C_{1}$ is abelian.
			
			\item \itemcase{T}he commutator $[\alpha_{1}(T_{1}),\alpha_{2}(T_{2})]$ is compact for all $T_{1}\in\C_{1}$, $T_{2}\in\C_{2}$.
			
			\item \itemcase{I}f $T_{1,k}\in\C_{1}$, $T_{2,k}\in\C_{2}$, $1\leq k\leq m$, then
			\begin{align*}
				\sum_{k=1}^{m}\alpha_{1}(T_{1,k})\alpha_{2}(T_{2,k})\in\K(\H)\Rightarrow\sum_{k=1}^{m}T_{1,k}\otimes T_{2,k}=0.
			\end{align*}
		\end{enumerate}
		Then there is a unique (in the uniform norm) continuous $*$-homomorphism $\sym:C^{\ast}(\alpha_1(\C_{1}),\alpha_2(\C_{2}))\rightarrow \C_{1}\otimes_{\rm min}\C_{2}$ such that
		\begin{align*}
			\sym(\alpha_{1}(T_{1}))=T_{1}\otimes1,\quad
			\sym(\alpha_{2}(T_{2}))=1\otimes T_{2},\quad\forall T_{1}\in\C_{1},T_{2}\in\C_{2}.
		\end{align*}
	\end{lemma}

	Let us now prove Theorem \ref{symbol}.
	\begin{proof}[Proof of Theorem \ref{symbol}]
		
		Letting $\C_{1}=\A_{1}$ and $\C_{2}=\A_{2}$, the condition (1) in Lemma \ref{symbol-existence} is obvious.
		The Lemma \ref{compactness} and Lemma \ref{primeness} yield the conditions (2) and (3) separately. Then Lemma \ref{symbol-existence} claims that there is a $*$-homomorphism $$\sym:\Pi\rightarrow\A_{1}\otimes_{\rm min}\A_{2}$$ such that
		\begin{align*}
			\sym(\po(f))=f\otimes1,\quad\sym(\pt(g))=1\otimes g,\quad f\in\A_{1},g\in\A_{2}.
		\end{align*}
		Recalling the proof of Theorem 3.3 in \cite{DAO2}, the $*$-homomorphism $\sym$ is constructed as a composition
		\begin{align*}
			\sym=\theta\circ q,
		\end{align*}
		where $\theta$ is some linear isomorphism and $$q:\B(\lp{2})\rightarrow\B(\lp{2})/\K(\lp{2})$$ is the quotient map. The quotient map $q$ kills all compact operators.
		By Lemma \ref{containalg}, we obtain $$\K(\lp{2})=\ker(\sym).$$
		
		It is left to prove that $\sym$ is surjective. 
		For $\displaystyle F=\sum_{k=1}^{m}f_{k}\otimes g_{k}\in\A_{1}\otimes\A_{2}$, take
		\begin{align*}
			T_{F}=\sum_{k=1}^{m}\po(f_{k})\pt(g_{k})\in\Pi.
		\end{align*}
		Since $\sym$ is a $*$-homomorphism, we have
		\begin{align*}
			\sym(T_{F})=\sum_{k=1}^{m}\sym(\po(f_{k}))\sym(\pt(g_{k}))=F.
		\end{align*}
		It follows that $\A_{1}\otimes\A_{2}\subset\sym(\Pi)$.
		By Theorem 4.1.9 in \cite{Kadison-Ringrose-book}, $*$-homomorphism image of a $C^{*}$-algebra is again a $C^{*}$-algebra. Then $\sym(\Pi)$ is a $C^{*}$-algebra. Because $\A_{1}\otimes\A_{2}$ is dense in $\A_{1}\otimes_{\rm min}\A_{2}$, this implies that $\sym$ is surjective.
	\end{proof}

	\section{Quasi-Riesz transforms}\label{quasi-Riesz-transform}
	
	In this section, we investigate the properties of quasi-Riesz transforms under $*$-homomorphism $\sym$. 
	Assume that $a=(a_{jk})_{j,k=1}^{d}:\GG\rightarrow\at$ is a smooth function. Set
	\begin{align*}
		X_{k}^{a}=\sum_{j=1}^{n_{1}}M_{a_{jk}}X_{j},\quad1\leq k\leq n_{1}\quad{\rm and}\quad L_{a}=\sum_{k=1}^{n_{1}}|X_{k}^{a}|^{2}.
	\end{align*}
	Denote
	\begin{align}\label{quasi-Riesz}
		\quasiriesz{a}{k}=X_{k}^{a}L_{a}^{-1/2},\quad1\leq k\leq n_{1}.
	\end{align}
	
	\begin{theorem}\label{symbolcalcu}
		Let $\psi\in\cci$ and let $a: \GG\rightarrow\at$ be a smooth function constant outside of a ball. For $1\leq k\leq n_{1}$, we have
		\begin{align*}
			M_{\psi} \quasiriesz{a}{k} M_{\psi}\in\Pi\quad{\rm and}\quad
			\sym(M_{\psi} \quasiriesz{a}{k} M_{\psi})(x)=\psi(x)^{2} \quasiriesz{a(x)}{k},\quad x\in\GG.
		\end{align*}
	\end{theorem}
	The proof of this theorem is arranged in subsection \ref{Proof of Theorem3.1} while we make some preparations in subsection \ref{Compactness section1}, subsection \ref{Estimates on quasi-Riesz transforms} and subsection \ref{Approximation property}.

	\subsection{Schatten classes}\label{Spdef}
	The following material is standard; for more details we refer to \cite{LSZv1,LMSZv2}.
	We recall the definitions of the Schatten classes $\mathcal{L}_{p}(H)$ and $\mathcal{L}_{p,\infty}(H)$ for a Hilbert space $H$. Let $\B(H)$ be the set of all bounded operators on $H.$ 
	The operator norm on $\B(H)$ is also denoted $\Norm{\cdot}_{\infty}$.
	Note that if $A$ is any compact operator on $H$, then $|A|$ is compact, positive and therefore diagonalizable. We define the singular values $\{\mu(k;A)\}_{k=0}^{\infty}$ to be the sequence of eigenvalues of $|A|$ (counted according to multiplicity). 

	For $0<p<\infty$, a compact operator $A$ on $H$ is said to belong to the Schatten class $\mathcal{L}_p(H)$ if $\{\mu(k,A)\}_{k=0}^{\infty}$ is $p$-summable, i.e. in the sequence space $l_p.$ 
	If $p\geq 1$, then the $\mathcal{L}_p(H)$ norm is defined as
	$$\|A\|_{\mathcal{L}_p(H)}=\left(\sum_{k=0}^{\infty}\mu(k,A)^p\right)^{1/p}.$$
	With this norm $\mathcal{L}_p(H)$ is a Banach space and, an ideal of $\B(H)$. 
	When $\frac{1}{r}=\frac{1}{p}+\frac{1}{q}$ with $p,q,r\geq1$, there is a constant $c_{p,q,r}$ such that below H\"{o}lder's inequality holds,
	\begin{align*}
		\Norm{AB}_{\L_{r}(H)}\leq c_{p,q,r}\Norm{A}_{\L_{p}(H)}\Norm{A}_{\L_{q}(H)}.
	\end{align*}
	
	For $0<p<\infty$, the weak Schatten class $\mathcal{L}_{p,\infty}(H)$ is the set consisting of operators $A$ on $H$ such that $\{\mu(k,A)\}_{k=0}^{\infty}$ is in $l_{p,\infty},$ with quasi-norm
	$$\|A\|_{\mathcal{L}_{p,\infty}(H)}=\sup_{k\geq0}(k+1)^{1/p}\mu(k,A).$$
	We have the H\"{o}lder-type inequality
	\begin{align*}
		\Norm{AB}_{\L_{r,\infty}(H)}\leq c_{p,q,r}\Norm{A}_{\L_{p,\infty}(H)}\Norm{A}_{\L_{q,\infty}(H)}
	\end{align*}
	for some constant $c_{p,q,r}$ with $\frac{1}{r}=\frac{1}{p}+\frac{1}{q}$ $(p,q,r\geq1)$.
	For convenience, we use the abbreviations $\mathcal{L}_p$ and $\mathcal{L}_{p,\infty}$ to denote $\mathcal{L}_p(H)$ and $\mathcal{L}_{p,\infty}(H)$ respectively, whenever the Hilbert space $H$ is clear from context.

	To continue, for $1\leq p<\infty$, we define an ideal $(\L_{p,\infty})_{0}$ by setting
	\begin{align*}
		(\L_{p,\infty})_{0}=\set{A\in\L_{p,\infty}:\lim_{k\rightarrow\infty}k\mu(k,A)^{p}=0}.
	\end{align*}
	This is a closed subspace of $\L_{p,\infty}$ and coincides with the closure of ideal of all finite rank operators in $\L_{p,\infty}$ and is commonly called the separable part of $\L_{p,\infty}$. 
	There are some elementary results:
	\begin{enumerate}[\rm(i)]
		\item $(\L_{p,\infty})_{0}\subset(\L_{p+\epsilon,\infty})_{0}$ for any $\epsilon>0$;
		\item $\L_{p,\infty}\cdot(\L_{q,\infty})_{0}\subset(\L_{r,\infty})_{0}$ and $(\L_{q,\infty})_{0}\cdot\L_{p,\infty}\subset(\L_{r,\infty})_{0}$ for $\frac{1}{r}=\frac{1}{p}+\frac{1}{q}$ with $1\leq p,q,r<\infty$;
		\item $\B(\H)\cdot(\L_{p,\infty})_{0}\subset(\L_{p,\infty})_{0}$ and $(\L_{p,\infty})_{0}\cdot\B(\H)\subset(\L_{p,\infty})_{0}$ (i.e. $(\L_{p,\infty})_{0}$ is a bi-ideal of $\B(\H)$).   
	\end{enumerate}
	Those ideal properties of $(\L_{p,\infty})_{0}$ are easy consequences of H\"{o}lder's inequality. 
	
	In addition, let us list two basic facts that are easy to check but we will use them frequently in this note.
	\begin{enumerate}[\rm (i)]
		\item If $T_{1}$ and $T_{2}$ are linear operators that are well-defined in a Hilbert space, then
		\begin{align}\label{commutator-inverse}
			[T_{1}\rp,T_{2}]=-T_{1}\rp[T_{1},T_{2}]T_{1}\rp.
		\end{align}
		
		\item The space of compact operators $\K(\lp{2})$ is a bi-ideal of $\B(\lp{2})$.
		
	\end{enumerate}

	\subsection{Compactness \expandafter{\romannumeral1}}\label{Compactness section1}
	\ 
	\newline
	For our purpose, we recall some function spaces at first. 
	Given $n\geq1$ and $1\leq p<\infty$, the $L_{p}$-space from $\GG$ to $\MM_{n}$ is defined to be 
	\begin{align*} 
		\mlp{p}=
		\set{w=(w_{kj})_{k,j=1}^{n}:\GG\rightarrow\MM_{n}\Big| \Norm{w_{kj}}_{\lp{p}}<\infty,\quad 1\leq k,j\leq n}
	\end{align*}
	equipped with the norm
	\begin{align*}
		\Norm{w}_{\mlp{p}}
		=\Big(\sum_{k,j=1}^{n}\Norm{w_{kj}}_{\lp{p}}^{p}\Big)^{1/p}.
	\end{align*}
	When $p=\infty$, the $L_{\infty}$-space from $\GG$ to $\MM_{n}$ is defined as
	\begin{align*} 
		\mlp{\infty}=
		\set{w=(w_{kj})_{k,j=1}^{n}:\GG\rightarrow\MM_{n}\Big| \Norm{w_{kj}}_{\lp{\infty}}<\infty,\quad 1\leq k,j\leq n}
	\end{align*}
	equipped with the norm
	\begin{align*}
		\Norm{w}_{\lgm{n}}=\sup_{k,j}\Norm{w_{kj}}_{\infty}.
	\end{align*}
	As in convention, $L_{p}(\GG;\MM_{1})=\lp{p}$ $(1\leq p<\infty)$ and $L_{\infty}(\GG;\MM_{1})=\lp{\infty}$.
	They are the usual $L^{p}$ and $L^{\infty}$ spaces on $\GG$.

	For $1\leq p<\infty$, the (homogeneous) Sobolev space $\wkp{p}$ is defined to be
	\begin{align*}
		\wkp{p}=
		\set{w=(w_{kj})_{k,j=1}^{n}:\GG\rightarrow\MM_{n}\Big| \Norm{X_{l}w_{kj}}_{\lp{p}}<\infty,\quad 1\leq l\leq n_{1}, 1\leq k,j\leq n}
	\end{align*}
	equipped with norm
	\begin{align*}
		\Norm{w}_{\wkp{p}}
		=\Big(\sum_{k,j=1}^{n}\sum_{l=1}^{n_{1}}\Norm{X_{l}w_{kj}}_{\lp{p}}^{p}\Big)^{1/p}.
	\end{align*}
	For $p=\infty$, the (homogeneous) Sobolev space $\wkp{\infty}$ is defined to be
	\begin{align*}
		\wkp{\infty}=
		\set{w=(w_{kj})_{k,j=1}^{n}:\GG\rightarrow\MM_{n}\Big| \Norm{X_{l}w_{kj}}_{\lp{\infty}}<\infty,\quad 1\leq l\leq n_{1}, 1\leq k,j\leq n}
	\end{align*}
	equipped with norm
	\begin{align*}
		\Norm{w}_{\wkp{\infty}}
		=\sum_{k,j=1}^{n}\sum_{l=1}^{n_{1}}\Norm{X_{l}w_{kj}}_{\lp{\infty}}.
	\end{align*}
	As in convention, $\dot{W}^{1,p}(\GG;\MM_{1})=\dot{W}^{1,p}(\GG)$ and $\dot{W}^{1,\infty}(\GG;\MM_{1})=\dot{W}^{1,\infty}(\GG)$. They are the usual Sobolev spaces on $\GG$.

	\begin{definition}\label{vqr}
		Assume that $w=(w_{jk})_{j,k=1}^{n_{1}}:\GG\rightarrow\glo$ is a differentiable function. Let
		\begin{align*}
			X_{k}^{w}=\sum_{j=1}^{n_{1}}M_{w_{jk}}X_{j},\quad1\leq k\leq n_{1}\quad{\rm and}\quad L_{w}=\sum_{k=1}^{n_{1}}|X_{k}^{w}|^{2}.
		\end{align*}
		Similar to \eqref{quasi-Riesz}, the quasi-Riesz transforms are defined as
		\begin{align*}
			\quasiriesz{w}{k}=X_{k}^{w}L_{w}^{-1/2},\quad1\leq k\leq n_{1}.
		\end{align*}
		When $w$ is the unit in $\glo$, it follows that $\quasiriesz{w}{k}$ is the usual $k$-th Riesz transform on $\GG$.
	\end{definition}

	\begin{lemma}\label{laplacinbd}
		Let $w:\GG\rightarrow\glo$ be a smooth function with bounded inverse. Then
		$$\Norm{\Delta^{1/2}L_{w}^{-1/2}}_{\infty}\leq\Norm{w\rp}_{\lgm{n_{1}}}.$$
	\end{lemma}
	
	\begin{proof}
		
		Let $\omega$ be the transpose of $w$. Writing $X$ as the row vector $(X_{1},X_{2},\,\dotsc,X_{n_{1}})$, we have
		\begin{align*}
			\Delta=X^{*}X
			=(\omega X)^{*} (\omega\rp)^{*}\omega\rp (\omega X)\leq\Norm{(\omega\rp)^{*}\omega\rp}_{\infty}(\omega X)^{*}\omega X.
		\end{align*}
		Therefore,
		\begin{align*}
			\Delta\leq\Norm{\omega\rp}_{\lgm{n_{1}}}^{2}L_{w}=\Norm{w\rp}_{\lgm{n_{1}}}^{2}L_{w}.
		\end{align*}
		Moreover,
		\begin{align*}
			L_{w}^{-1/2}\Delta L_{w}^{-1/2}\leq\Norm{w\rp}_{\lgm{n_{1}}}^{2}.
		\end{align*}
		Since 
		\begin{align*}
			|\Delta^{1/2}L_{w}^{-1/2}|^{2}=L_{w}^{-1/2}\Delta L_{w}^{-1/2}
			\quad\mbox{and}\quad 
			\Norm{|\Delta^{1/2}L_{w}^{-1/2}|^{2}}_{\infty}
			=\Norm{\Delta^{1/2}L_{w}^{-1/2}}_{\infty}^{2},
		\end{align*}
		it implies the desired result.
	\end{proof}

	\begin{lemma}\label{laplacian-bd-14}
		Let $w:\GG\rightarrow\glo$ be given as in Lemma \ref{laplacinbd}. Then
		$$\Norm{\Delta^{1/4}L_{w}^{-1/4}}_{\infty}\leq\Norm{w\rp}_{\lgm{n_{1}}}^{1/2}.$$
	\end{lemma}
	\begin{proof}
		This assertion follows from Lemma \ref{laplacinbd} and Hadamard's three lines lemma \cite[Lem 1.3.5]{G2014c}.
	\end{proof}

	\begin{lemma}\label{quasiriszbd}
		Let $\psi\in\cci$ and let $w$ be given as in Lemma \ref{laplacinbd}. Then $M_{\psi}\quasiriesz{w}{k}$ is bounded on $\lp{2}$ for $1\leq k\leq n_{1}$.
	\end{lemma}
	\begin{proof}
		Let $w=(w_{jk})_{j,k=1}^{n_{1}}$ and write
		\begin{align*}
			M_{\psi}\quasiriesz{w}{k}
			=\sum_{j=1}^{n_{1}}M_{\psi w_{jk}}X_{j}\Delta^{-1/2}\cdot\Delta^{1/2}L_{w}^{-1/2}.
		\end{align*}
		Since $\psi\in\cci$, it induces that $M_{\psi w_{kj}}$ is bounded on $\lp{2}$. Thus, the boundedness of Riesz transforms $\{ R_{j}\}_{j=1}^{n_{1}}$ and Lemma \ref{laplacinbd} yield the boundedness of $M_{\psi}\quasiriesz{w}{k}$. This completes the proof.
	\end{proof}

	\begin{lemma}\label{rieszcompactt}
		Let $\psi_{1},\psi_{2}\in\cci$ and let $w$ be given as in Lemma \ref{laplacinbd}. Then
		$$M_{\psi_{1}}\Big(\quasiriesz{w}{k}-X_{k}^{w}(1+L_{w})^{-1/2}\Big)M_{\psi_{2}}\in\K(\lp{2}).$$
	\end{lemma}
	\begin{proof}
		Denote ${\rm LHS}=M_{\psi_{1}}\Big(\quasiriesz{w}{k}-X_{k}^{w}(1+L_{w})^{-1/2}\Big)M_{\psi_{2}}$. Note that
		\begin{align*}
			{\rm LHS}=M_{\psi_{1}}\quasiriesz{w}{k}\Big(1-L_{w}^{1/2}(1+L_{w})^{-1/2}\Big)M_{\psi_{2}}.
		\end{align*}
		Since
		\begin{align*}
			1-L_{w}^{1/2}(1+L_{w})^{-1/2}
			=\frac{1}{(1+L_{w})^{1/2}((1+L_{w})^{1/2}+L_{w}^{1/2})},
		\end{align*}
		it follows that
		\begin{align*}
			{\rm LHS}
			=M_{\psi_{1}}\quasiriesz{w}{k}\cdot\frac{L_{w}^{1/2}}{(1+L_{w})^{1/2}((1+L_{w})^{1/2}+L_{w}^{1/2})}
			\cdot L_{w}^{-1/2}\Delta^{1/2}\cdot\Delta^{-1/2}M_{\psi_{2}}.
		\end{align*}
		On the right hand side,  applying Lemma \ref{quasiriszbd} and Lemma \ref{laplacinbd} separately, the first and third factors are bounded. The second factor is bounded by functional calculus. The last factor is compact by Cwikel estimate in \cite{MSZ-Cwikel}. Therefore, ${\rm LHS}$ is a product of bounded operators and compact operator, and thus it is a compact operator. 
	\end{proof}

	The following lemma is a version of \cite[Lemma 3.5.2]{DAO4}.
	\begin{lemma}\label{deltadecomo}
		Assume that $T\in\B(\H)$ is positive self-adjoint and bounded from below by a strictly positive constant. Then
		\begin{align*}
			T^{-1/2}=\frac{2}{\pi}\int_{0}^{\infty}\frac{dt}{T+t^{2}}
		\end{align*}
		where this integral is understood in the sense of strong operator topology.
	\end{lemma}

	\begin{lemma}\label{tail2-compact}
		Let $\psi,\phi\in\cci$ and let $w$ be given as in Lemma \ref{laplacinbd}. Then
		\begin{align*}
			[M_{\psi},(1+L_{w})^{1/2}](1+L_{w})^{-1/2}M_{\phi}\in\K(\lp{2}).
		\end{align*}
	\end{lemma}
	\begin{proof}
		Denote ${\rm Term}=[M_{\psi},(1+L_{w})^{1/2}](1+L_{w})M_{\phi}$.
		
		It follows from Lemma \ref{deltadecomo} that
		\begin{align*}
			(1+L_{w})^{1/2}=\frac{2}{\pi}\int_{0}^{\infty}\frac{1+L_{w}}{1+t^{2}+L_{w}}\,dt,
		\end{align*}
		in the strong operator topology on the domain of $L_w.$ Taking this equality into consideration, one has
		\begin{align*}
			[M_{\psi},(1+L_{w})^{1/2}]
			=-\frac{2}{\pi}\int_{0}^{\infty}[M_{\psi},\frac{1}{1+t^{2}+L_{w}}]t^{2}\,dt.
		\end{align*}
		
		By \eqref{commutator-inverse}, ${\rm Term}$ can be described in the following way
		\begin{align*}
			{\rm Term}
			=\frac{2}{\pi}\int_{0}^{\infty}\frac{t^{2}}{1+t^{2}+L_{w}}[M_{\psi},L_{w}](1+L_{w})^{-1/2}
			\frac{1}{1+t^{2}+L_{w}}M_{\phi}dt.
		\end{align*}
		Note that
		\begin{align}\label{ML-commutator}
			[M_{\psi},L_{w}]=\sum_{k,j,l=1}^{n_{1}}\Big(M_{X_{j}(\overline{w_{jk}}w_{lk}X_{l}\psi)}
			+M_{(\overline{w_{jk}}w_{lk}+\overline{w_{lk}}w_{jk})X_{l}\psi}X_{j}\Big).
		\end{align}
		Let
		\begin{align*}
			\psi_{kjl}^{(1)}=X_{j}(\overline{w_{jk}}w_{lk}X_{l}\psi)\quad\mbox{and}\quad
			\psi_{kjl}^{(2)}=(\overline{w_{jk}}w_{lk}+\overline{w_{lk}}w_{jk})X_{l}\psi.
		\end{align*}
		Then
		\begin{align*}
			{\rm Term}
			&=\frac{2}{\pi}\sum_{k,j,l=1}^{n_{1}}\int_{0}^{\infty}\frac{t^{2}}{1+t^{2}+L_{w}}M_{\psi_{kjl}^{(1)}}(1+L_{w})^{-1/2}
			\frac{1}{1+t^{2}+L_{w}}M_{\phi}dt+\\
			&+\frac{2}{\pi}\sum_{k,j,l=1}^{n_{1}}\int_{0}^{\infty}\frac{t^{2}}{1+t^{2}+L_{w}}M_{\psi_{kjl}^{(2)}}X_{j}(1+L_{w})^{-1/2}
			\frac{L_{w}^{1/4}}{1+t^{2}+L_{w}}\cdot L_{w}^{-1/4}M_{\phi}dt\\
			&=:{\rm Term1}+{\rm Term2}.
		\end{align*}
		
		Write
		\begin{align*}
			M_{\psi_{kjl}^{(1)}}(1+L_{w})^{-1/2}
			=M_{\psi_{kjl}^{(1)}}\Delta^{-1/2}\cdot \Delta^{1/2}L_{w}^{-1/2}\cdot \frac{L_{w}^{1/2}}{(1+L_{w})^{1/2}}.
		\end{align*}
		By Cwikel estimate in \cite{MSZ-Cwikel}, the above first factor on the right hand side is compact.
		By Lemma \ref{laplacinbd} and functional calculus, the second and third factors are bounded separately.
		Thus $M_{\psi_{kjl}^{(1)}}(1+L_{w})^{-1/2}$ is a compact operator.
		Moreover,
		\begin{align*}
			\Norm{{\rm Term1}}_{\infty}
			\leq\frac{2}{\pi}\Norm{w\rp}_{\lgm{n_{1}}}\sum_{k,j,l=1}^{n_{1}}\Norm{M_{\psi_{kjl}^{(1)}}\Delta^{-1/2}}_{\infty}
			\int_{0}^{\infty}\frac{dt}{1+t^{2}}<\infty.
		\end{align*}
		It follows that ${\rm Term1}$ is convergent in the uniform norm and therefore it is compact.
		
		Write
		\begin{align*}
			M_{\psi_{kjl}^{(2)}}X_{j}(1+L_{w})^{-1/2}
			=M_{\psi_{kjl}^{(2)}}X_{j}\Delta^{-1/2}\cdot\Delta^{1/2}L_{w}^{-1/2}\cdot \frac{L_{w}^{1/2}}{(1+L_{w})^{1/2}}.
		\end{align*}
		It follows from Lemma \ref{laplacinbd} and functional calculus that $M_{\psi_{kjl}^{(2)}}X_{j}(1+L_{w})^{-1/2}$ is a bounded operator.
		Write
		\begin{align}\label{L14-compact}
			L_{w}^{-1/4}M_{\phi}=L_{w}^{-1/4}\Delta^{1/4}\cdot \Delta^{-1/4}M_{\phi}.
		\end{align}
		Applying Lemma \ref{laplacian-bd-14}, the above first factor on the right hand side is bounded,
		and using Cwikel estimate in \cite{MSZ-Cwikel}, the second factor is compact.
		By spectral theorem,
		\begin{align}\label{LL-spectral}
			\Norm{\frac{L_{w}^{1/4}}{1+t^{2}+L_{w}}}_{\infty}
			\leq\sup_{s\geq0}\frac{s}{1+t^{2}+s^{4}}
			=\frac{3}{4}3^{-1/4}(1+t^{2})^{-3/4}.
		\end{align}
		Therefore,
		\begin{align*}
			\Norm{{\rm Term2}}_{\infty}
			\leq\frac{3}{2\pi}3^{-1/4}\Norm{w\rp}_{\lgm{n_{1}}}^{1/2}
			\sum_{k,j,l=1}^{n_{1}}\Norm{M_{\psi_{kjl}^{(2)}}X_{j}\Delta^{-1/2}}_{\infty}
			\int_{0}^{\infty}\frac{dt}{(1+t^{2})^{3/4}}<\infty.
		\end{align*}
		This implies that ${\rm Term2}$ is convergent in the uniform norm. Moreover it is a compact operator by compactness of $L_{w}^{-1/4}M_{\phi}$.
		This completes the proof.
		
	\end{proof}

	\begin{proposition}\label{compactrieszcommu}
		Let $\psi_{1},\psi_{2},\psi_{3}\in\cci$ and let $w$ be given as in Lemma \ref{laplacinbd}. Then
		\begin{align*}
			M_{\psi_{1}}[M_{\psi_{2}},\quasiriesz{w}{k}]M_{\psi_{3}}\in\K(\lp{2}).
		\end{align*}
	\end{proposition}
	\begin{proof}
		Let $w=(w_{jk})_{j,k=1}^{n_{1}}$ and ${\rm LHS}=M_{\psi_{1}}[M_{\psi_{2}},\quasiriesz{w}{k}]M_{\psi_{3}}$.
		Write
		\begin{align*}
			{\rm LHS}
			=M_{\psi_{1}}[M_{\psi_{2}},\quasiriesz{w}{k}-X_{k}^{w}(1+L_{w})^{-1/2}]M_{\psi_{3}}
			+M_{\psi_{1}}[M_{\psi_{2}},X_{k}^{w}(1+L_{w})^{-1/2}]M_{\psi_{3}}.
		\end{align*}
		Applying Lemma \ref{rieszcompactt}, the first commutator on the right hand side is compact. It suffices to show the compactness of the second commutator on the right hand side.
		
		By Leibniz rule, we can split the second commutator into a sum with the following form
		\begin{align*}
			M_{\psi_{1}}[M_{\psi_{2}},X_{k}^{w}](1+L_{w})^{-1/2}M_{\psi_{3}}
			+M_{\psi_{1}}X_{k}^{w}[M_{\psi_{2}},(1+L_{w})^{-1/2}]M_{\psi_{3}}=:{\rm Term1}-{\rm Term2}.
		\end{align*}
		
		Note that
		\begin{align*}
			[M_{\psi_{2}},X_{k}^{w}]=-\sum_{j=1}^{n_{1}}M_{w_{jk}}M_{X_{j}\psi_{2}}.
		\end{align*}
		Rewrite
		\begin{align*}
			{\rm Term1}
			=-\sum_{j=1}^{n_{1}}M_{\psi_{1}w_{jk}}M_{X_{j}\psi_{2}}\cdot\frac{L_{w}^{1/2}}{(1+L_{w})^{1/2}}
			\cdot L_{w}^{-1/2}\Delta^{1/2}\cdot\Delta^{-1/2}M_{\psi_{3}}.
		\end{align*}
		The first factor in ${\rm Term1}$ on the right hand side is bounded. 
		The second and third factors are bounded by applying functional calculus and Lemma \ref{laplacinbd} separately. Therefore, ${\rm Term1}$ is compact by Cwikel estimate of the last factor in \cite{MSZ-Cwikel}.
		
		By \eqref{commutator-inverse},
		\begin{align*}
			[M_{\psi_{2}},(1+L_{w})^{-1/2}]
			=-(1+L_{w})^{-1/2}[M_{\psi_{2}},(1+L_{w})^{1/2}](1+L_{w})^{-1/2}.
		\end{align*}
		This implies that
		\begin{align*}
			{\rm Term2}=M_{\psi_{1}}X_{k}^{w}(1+L_{w})^{-1/2}[M_{\psi_{2}},(1+L_{w})^{1/2}](1+L_{w})^{-1/2}M_{\psi_{3}}.
		\end{align*}
		Rewrite
		\begin{align*}
			{\rm Term2}=M_{\psi_{1}}X_{k}^{w}L_{w}^{-1/2}\cdot
			\frac{L_{w}^{1/2}}{(1+L_{w})^{1/2}}\cdot
			[M_{\psi_{2}},(1+L_{w})^{1/2}](1+L_{w})^{-1/2}M_{\psi_{3}}.
		\end{align*}
		Applying Lemma \ref{quasiriszbd}, the above first factor of ${\rm Term2}$ on the right hand side is bounded.
		By functional calculus, the second factor of ${\rm Term2}$ is bounded.
		Using Lemma \ref{tail2-compact}, the last factor of ${\rm Term2}$ belongs to $\K(\lp{2})$.
		Therefore, ${\rm Term2}$ is a product of bounded operators and compact operator and thus it is a compact operator.
		This completes the proof.
	\end{proof}

	\begin{proposition}\label{compactrieszquo}
		Let $\psi\in\cci$. Suppose that $w,v:\GG\rightarrow\glo$ are smooth functions with bounded inverse. If $w=v$ in some neighborhood of $\supp(\psi)$, then
		\begin{align*}
			M_{\psi}(\quasiriesz{w}{k}-\quasiriesz{v}{k})M_{\psi}\in\K(\lp{2}).
		\end{align*}
	\end{proposition}
	\begin{proof}
		Write
		\begin{align*}
			M_{\psi}(\quasiriesz{w}{k}-\quasiriesz{v}{k})M_{\psi}={\rm Term1}+{\rm Term2}+{\rm Term3},
		\end{align*}
		where
		\begin{align*}
			{\rm Term1}=M_{\psi}\Big(\quasiriesz{w}{k}-X_{k}^{w}(1+L_{w})^{-1/2}
			\Big)M_{\psi},\quad
			{\rm Term2}=M_{\psi}\Big(X_{k}^{v}(1+L_{v})^{-1/2}-\quasiriesz{v}{k}\Big)M_{\psi}
		\end{align*}
		and
		\begin{align*}
			{\rm Term3}
			=M_{\psi}X_{k}^{w}(1+L_{w})^{-1/2}M_{\psi}-M_{\psi}X_{k}^{v}(1+L_{v})^{-1/2}M_{\psi}.
		\end{align*}
		Applying Lemma \ref{rieszcompactt}, ${\rm Term1}$ and ${\rm Term2}$ are compact. It is left to show ${\rm Term3}\in\K(\lp{2})$.
		
		Since $w=v$ in some neighborhood of $\supp(\psi)$, it follows that $M_{\psi}X_{k}^{w}=M_{\psi}X_{k}^{v}$.
		Therefore,
		\begin{align*}
			{\rm Term3}=M_{\psi}X_{k}^{w}\Big((1+L_{w})^{-1/2}-(1+L_{v})^{-1/2}\Big)M_{\psi}.
		\end{align*}
		Using Lemma \ref{deltadecomo}, we have
		\begin{align*}
			{\rm Term3}
			&=\frac{2}{\pi}M_{\psi}X_{k}^{w}\int_{0}^{\infty}(\frac{1}{1+t^{2}+L_{w}}-\frac{1}{1+t^{2}+L_{v}})M_{\psi}\,dt\\
			&=\frac{2}{\pi}M_{\psi}X_{k}^{w}\int_{0}^{\infty}\frac{1}{1+t^{2}+L_{w}}(L_{v}-L_{w})\frac{1}{1+t^{2}+L_{v}}M_{\psi}\,dt.
		\end{align*}
		
		Write $\psi=\psi_{1}\psi_{2}$ with $\psi_{1},\psi_{2}\in\cci$ and $\supp\psi_{1}=\supp\psi$. Then,
		\begin{align*}
			(L_{v}-L_{w})\frac{1}{1+t^{2}+L_{v}}M_{\psi}&=(L_{v}-L_{w})M_{\psi_{1}}\frac{1}{1+t^{2}+L_{v}}M_{\psi_{2}}\\
			&+(L_{v}-L_{w})[\frac{1}{1+t^{2}+L_{v}},M_{\psi_{1}}]M_{\psi_{2}}.
		\end{align*}
		Clearly, $(L_{v}-L_{w})M_{\psi_{1}}=0$. This implies that
		\begin{align*}
			{\rm Term3}
			&=\frac{2}{\pi}M_{\psi}X_{k}^{w}\int_{0}^{\infty}\frac{1}{1+t^{2}+L_{w}}(L_{v}-L_{w})[\frac{1}{1+t^{2}+L_{v}},M_{\psi_{1}}]M_{\psi_{2}}\,dt.
		\end{align*}
		
		By \eqref{commutator-inverse},
		\begin{align*}
			[\frac{1}{1+t^{2}+L_{v}},M_{\psi_{1}}]=-\frac{1}{1+t^{2}+L_{v}}[L_{v},M_{\psi_{1}}]\frac{1}{1+t^{2}+L_{v}}.
		\end{align*}
		Note that $L_{v}-L_{w}=(1+t^{2}+L_{v})-(1+t^{2}+L_{w})$. Since $w=v$ in some neighborhood of $\supp(\psi)=\supp(\psi_{1})$, it follows that
		\begin{align*}
			{\rm Term3}
			&=\frac{2}{\pi}M_{\psi}\quasiriesz{v}{k}\int_{0}^{\infty}\frac{L_{v}}{1+t^{2}+L_{v}}\cdot L_{v}^{-1/2}[L_{v},M_{\psi_{1}}]\cdot\frac{1}{1+t^{2}+L_{v}}M_{\psi_{2}}\,dt\\
			&-\frac{2}{\pi}M_{\psi}\quasiriesz{w}{k}\int_{0}^{\infty}\frac{L_{w}}{1+t^{2}+L_{w}}\cdot L_{w}^{-1/2}[L_{v},M_{\psi_{1}}]\cdot\frac{1}{1+t^{2}+L_{v}}M_{\psi_{2}}\,dt\\
			&=:{\rm Term31}+{\rm Term32}.
		\end{align*}
		Clearly, 
		\begin{align*}
			[L_{v},M_{\psi_{1}}]=\sum_{k,j,l=1}^{n_{1}}M_{X_{l}(\overline{v_{jk}}v_{lk}(X_{j}\psi_{1}))}-
			\sum_{k,j,l=1}^{n_{1}}X_{l}M_{(\overline{v_{lk}}v_{jk}+\overline{v_{jk}}v_{lk})(X_{j}\psi_{1})},
		\end{align*}
		where $v_{jk}$ is the entries of matrix $v$.
		Therefore,
		\begin{align}\label{LM-commutator-decom}
			L_{v}^{-1/2}[L_{v},M_{\psi_{1}}]
			=\sum_{k,j,l=1}^{n_{1}}L_{v}^{-1/2}\Delta^{1/2}
			\cdot\Delta^{-1/2}M_{X_{l}(\overline{v_{jk}}v_{lk}(X_{j}\psi_{1}))}
			-\sum_{k,j,l=1}^{n_{1}}L_{v}^{-1/2}\Delta^{1/2}
			\cdot\Delta^{-1/2}X_{l}M_{(\overline{v_{lk}}v_{jk}+\overline{v_{jk}}v_{lk})(X_{j}\psi_{1})}.
		\end{align}
		It follows from Lemma \ref{laplacinbd} and Cwikel estimate in \cite{MSZ-Cwikel} that the above first factor on the right hand side is bounded.
		And Lemma \ref{laplacinbd} and Lemma \ref{quasiriszbd} implies the boundedness of the above second factor.
		Moreover, Write
		\begin{align*}
			\frac{1}{1+t^{2}+L_{v}}M_{\psi_{2}}
			=\frac{L_{v}^{1/4}}{1+t^{2}+L_{v}}\cdot L_{v}^{-1/4}M_{\psi_{2}}.
		\end{align*}
		Similar to \eqref{L14-compact}, $L_{v}^{-1/4}M_{\psi_{2}}$ is a compact operator.
		By \eqref{LL-spectral}, 
		\begin{align}\label{Lphi-bd-compact}
			\Norm{\frac{1}{1+t^{2}+L_{v}}M_{\psi_{2}}}_{\infty}
			\leq\frac{3}{4}3^{-1/4}(1+t^{2})^{-3/4}\Norm{L_{v}^{-1/4}M_{\psi_{2}}}_{\infty}.
		\end{align}
		Combining Lemma \ref{quasiriszbd}, functional calculus, \eqref{LM-commutator-decom} and \eqref{Lphi-bd-compact},  ${\rm Term31}$ is convergent in the uniform norm and even more it is a compact operator. 
		Replacing $L_{v}^{-1/2}$ in \eqref{LM-commutator-decom} by $L_{w}^{-1/2}$, it follows that $L_{w}^{-1/2}[L_{v},M_{\psi_{1}}]$ is a bounded operator.
		Applying a similar method of ${\rm Term31}$ to ${\rm Term32}$, ${\rm Term32}$ is also convergent in the uniform norm and even more it is compact. 
		Thus ${\rm Term3}$ is compact. This completes the proof.
	\end{proof}

	\subsection{Estimates on quasi-Riesz transforms}\label{Estimates on quasi-Riesz transforms}

	The topics of heat kernel associated to the heat semigroup on stratified Lie groups can be found in \cite[IV.4]{VSC1992} and \cite[Chapter 1. G]{FS1982}.
	\begin{lemma}\label{heat-matrix}
		Let $A\in\at.$ If $u^{A}$ is the $C^{\infty}$ solution of the hypo-elliptic equation $(\frac{\p}{\p t}+L_{A})u^A=0$ satisfying the initial condition $u^A(0,\cdot)=\delta,$ then
		\begin{align*}
			u^{A}(t,x)=|\det(A\rp)|u\left(t,A\rp x\right).
		\end{align*}
		Here, $u$ is the $C^{\infty}$ solution of the hypo-elliptic equation $(\frac{\p}{\p t}+\Delta)u=0$ satisfying the initial condition $u(0,\cdot)=\delta.$
	\end{lemma}
	\begin{proof} Obvious.
		
	\end{proof}
	
	The following assertion is the standard Calder\'{o}n-Zygmund theory in harmonic analysis on spaces of homogeneous type in \cite[Chapter 1. Section 1.4]{DH-Book2009} (see \cite[Chapter 5]{G2014c} for Euclidean analogue).
	
	\begin{lemma}\label{second-Riesz-bd}
		The operators $\{X_{l}X_{k}\Delta\rp\}_{l,k=1}^{n_{1}}$ are bounded on $\lp{2}.$ 
	\end{lemma}
	
	For $A\in\at$, we can identify it as a block-diagonal matrix associated to the stratification of $\g$ in the following form
	\begin{align}\label{diagonal}
		A=\diag(A_{(1)},A_{(2)},\cdots,A_{(\iota)}).
	\end{align}	
	\begin{lemma}\label{laplacianbdd}
		Suppose $A\in\at$. Then there is a constant $C$ independent of $A$ such that
		\begin{align*}
			\Norm{\Delta L_{A}^{-1}}_{\infty}\leq C \Norm{A_{(1)}\rp}_{\lm{n_{1}}}^{2},
		\end{align*}
		where $A_{(1)}$ is the first block of $A$ as in \eqref{diagonal}.
	\end{lemma}
	\begin{proof}
		Recall the formula
		\begin{align*}
			L_{A}^{-1}=\int_{0}^{\infty}e^{-tL_{A}}\,dt
		\end{align*}
		in the strong operator topology.
		By \cite[IV.4]{VSC1992}, one has $e^{-tL_{A}}f(x)=(h_{t}^{A}\ast f)(x)$ for $f\in\cci$.
		Then
		\begin{align*}
			\Delta L_{A}^{-1}f(x)=-\sum_{j=1}^{n_{1}}\int_{\GG}\int_{0}^{\infty}X_{j}^{2}h_{t}^{A}(y\rp x)\,dt\cdot f(y)\,dy.
		\end{align*}
		Denote
		\begin{align*}
			\Psi^{A}_{j}(x)=\int_{0}^{\infty}X_{j}^{2}h_{t}^{A}(x)\,dt,\quad x\in\GG.
		\end{align*}
		By a re-scaling in \cite[Proposition (1.68)]{FS1982}, it follows that $h_{t}^{A}(x)=t^{-Q/2}h_{1}^{A}(\delta_{t^{-1/2}}(x))$.
		Therefore, by Lemma \ref{heat-matrix}, one has
		\begin{align}\label{kernel-heat-A}
			\Psi^{A}_{j}(x )=|\det(A\rp)|
			\int_{0}^{\infty}t^{-Q/2}X_{j}^{2}\left(h_{1}\left(A\rp\delta_{t^{-1/2}}(x)\right)\right)\,dt.
		\end{align}
		Since $A\in\at$ has the form of \eqref{diagonal}, set
		\begin{align*}
			A\rp=\diag(A_{(1)}\rp,b_{rest})
			=\diag\left((b_{ij})_{i,j=1}^{n_{1}},b_{rest}\right).
		\end{align*}
		By \cite[Proposition 1.2.16]{BLU2007},
		\begin{align*}
			X_{j}\left(h_{1}\left(A\rp\delta_{t^{-1/2}}(x)\right)\right)
			&=\frac{d}{ds}\Big|_{s=0} h_{1}\left(A\rp\delta_{t^{-1/2}}(x\cdot se_{j})\right)\\
			&=\frac{d}{ds}\Big|_{s=0} h_{1}\left(A\rp\delta_{t^{-1/2}}(x)\cdot \delta_{t^{-1/2}}A\rp(se_{j})\right)\\
			&=t^{-1/2}\sum_{k=1}^{n_{1}}b_{kj}(X_{k}h_{1})\left(A\rp\delta_{t^{-1/2}}(x)\right).
		\end{align*}
Consequently,
		\begin{align}\label{derivative-heat-A}
			X_{j}^{2}\Big(h_{1}(A\rp\delta_{t^{-1/2}}(x))\Big)
			=t^{-1}\sum_{k,l=1}^{n_{1}}b_{kj}b_{lj}(X_{l}X_{k}h_{1})\left(A\rp\delta_{t^{-1/2}}(x)\right).
		\end{align}
		Let
		\begin{align*}
			\Psi_{l,k}(x)=\int_{0}^{\infty}(X_{l}X_{k}h_{1})(\delta_{t^{-1/2}}(x))\,t^{-Q/2-1}dt,\quad x\in\GG,
		\end{align*}
		and
		\begin{align*}
			T^{A}_{l,k}(f)(x)=\int_{\GG}\Psi_{l,k}(A\rp (y\rp x))f(y)\,dy.
		\end{align*}
		Let $\unit_{d\times d}$ be the $d\times d$ matrix unit. 
		Applying Lemma \ref{second-Riesz-bd}, there is constant $C$ independent of $A$ such that
		\begin{align*}
			\Norm{T^{\unit_{d\times d}}_{l,k}(f\circ A)}_{\lp{2}}
			\leq C \Norm{f\circ A}_{\lp{2}}.
		\end{align*}
		Observe that
		\begin{align*}
			\Norm{f\circ A}_{\lp{2}}=|\det(A\rp)|^{1/2}\Norm{f}_{\lp{2}}.
		\end{align*}
		Therefore, by $A\in\at$,
		\begin{align}\label{SIO-A}
			\Norm{T^{A}_{l,k}(f)}_{\lp{2}}
			=|\det(A)|^{3/2}\Norm{T^{\unit_{d\times d}}_{l,k}(f\circ A)}_{\lp{2}}
			\leq C |\det(A)|\Norm{f}_{\lp{2}}.
		\end{align}

		Since $A\in\at$, combining equalities \eqref{kernel-heat-A} and \eqref{derivative-heat-A}, we obtain 
		\begin{align*}
			\Psi^{A}_{j}(x)=|\det(A\rp)|\sum_{k,l=1}^{n_{1}}b_{kj}b_{lj}\Psi_{l,k}(A\rp x).
		\end{align*}
		Then
		\begin{align*}
			\Delta L_{A}^{-1}f(x)=-\sum_{j=1}^{n_{1}}\int_{\GG}\Psi^{A}_{j}(y\rp x)f(y)\,dy
			=-|\det(A\rp)|\sum_{j,k,l=1}^{n_{1}}b_{kj}b_{lj} T^{A}_{l,k}(f)(x).
		\end{align*}
		Thus, by triangle inequality and \eqref{SIO-A},
		\begin{align*}
			\Norm{\Delta L_{A}^{-1}f}_{\lp{2}}
			\leq|\det(A\rp)|\sum_{j,k,l=1}^{n_{1}}|b_{kj}||b_{lj}| \Norm{T^{A}_{l,k}(f)}_{\lp{2}}
			\leq C \sum_{j,k,l=1}^{n_{1}}|b_{kj}||b_{lj}| \Norm{f}_{\lp{2}}.
		\end{align*}
		This implies the desired result.
	\end{proof}

	\begin{lemma}\label{modulus}
		Let $f$ be a bounded smooth function on $\GG$. Then
		\begin{align*}
			\Norm{\Delta^{1/2}M_{f}\Delta^{-1/2}}_{\infty}
			\leq C(\Norm{f}_{\dot{W}^{1,Q}}+\Norm{f}_{\infty}).
		\end{align*}
	\end{lemma}
	\begin{proof}
		Notice that
		\begin{align*}
			\Norm{\Delta^{1/2}M_{f}\Delta^{-1/2}}_{\infty}
			\leq\Norm{\Delta^{1/2}\Bigg(\sum_{j=1}^{n_{1}}|X_{j}|\Bigg)\rp}_{\infty}
			\Norm{\sum_{j=1}^{n_{1}}|X_{j}|M_{f}\Delta^{-1/2}}_{\infty}.
		\end{align*}
		Since $\Delta^{1/2}(\sum_{j=1}^{n_{1}}|X_{j}|)\rp$ is a 0-homogeneous operator on $\GG$, it follows from \cite[Theorem 3.2.30]{FR-Book} that the first factor on the right hand side is finite.
		It is left to investigate the second factor.
		
		Clearly, $$X_{j}M_{f}\Delta^{-1/2}=M_{X_{j}f}\Delta^{-1/2}+M_{f}X_{j}\Delta^{-1/2}.$$
		By Cwikel estimate in \cite{MSZ-Cwikel}, the operator $M_{X_{j}f}\Delta^{-1/2}$ is compact and belongs to $\L_{Q,\infty}$. Moreover,
		\begin{align*}
			\Norm{M_{X_{j}f}\Delta^{-1/2}}_{\infty}
			\leq\Norm{M_{X_{j}f}\Delta^{-1/2}}_{\L_{Q,\infty}}
			\lesssim\Norm{X_{j}f}_{\lp{Q}}.
		\end{align*}
		Therefore, by triangle inequality,
		\begin{align*}
			\Norm{X_{j}M_{f}\Delta^{-1/2}}_{\infty}
			&\leq\Norm{M_{X_{j}f}\Delta^{-1/2}}_{\infty}+\Norm{M_{f}R_{j}}_{\infty}\\
			&\lesssim\Norm{X_{j}f}_{\lp{Q}}+\Norm{f}_{\infty}.
		\end{align*}
		Applying polar decomposition, there is a partial isometry $u_{j}$ such that $|X_{j}|=u_{j}X_{j}$.
		Combining $\Norm{u_{j}}_{\infty}\leq1$, one has
		\begin{align*}
			\Norm{\sum_{j=1}^{n_{1}}|X_{j}|M_{f}\Delta^{-1/2}}_{\infty}
			\leq\sum_{j=1}^{n_{1}}\Norm{X_{j}M_{f}\Delta^{-1/2}}_{\infty}
			\lesssim\Norm{f}_{\dot{W}^{1,Q}}+\Norm{f}_{\infty}.
		\end{align*}
		This yields the desired result.
	\end{proof}

	The following assertion is a generalization of the Riesz transforms as 0-homogeneous operators.
	Since $\Delta^{\alpha/2}X_{k}\Delta^{-\alpha+1/2}$ is 0-homogeneous operators, it follows from \cite[Theorem 3.2.30]{FR-Book} (see also \cite{CG1984}) that
	\begin{align}\label{0-order-bd}
		\Delta^{\alpha/2}X_{k}\Delta^{-\alpha+1/2}\in\B(\lp{2})
		\quad \mbox{if}\quad\alpha\in\RR.
	\end{align}

	\begin{lemma}\label{w1qa-b}
		Let $v\in\glo$ and $w:\GG\rightarrow\glo$ be a bounded smooth function. Then
		\begin{align*}
			\Norm{\Delta^{-1/4}(L_{w}-L_{v})\Delta^{-3/4}}_{\infty}
			\lesssim &\Norm{w-v}_{\lgm{n_{1}}} (\Norm{w}_{\lgm{n_{1}}}+\Norm{v}_{\lm{n_{1}}})\\
			&+\Norm{w}_{\lgm{n_{1}}}\Norm{w}_{\sgm}.
		\end{align*}
		Here, the constant only depends on $\GG$.
	\end{lemma}
	\begin{proof}
		Let $v=(v_{jk})$ and $w=(w_{jk})$.
		Denote
		\begin{align*}
			\Omega_{k,j_{1},j_{2}}=\Delta^{1/4}M_{v_{j_{1}k}v_{j_{2}k}-w_{j_{1}k}w_{j_{2}k}}\Delta^{-1/4}.
		\end{align*}
		Applying Lemma \ref{modulus} to $f_{k,j_{1},j_{2}}=v_{j_{1}k}v_{j_{2}k}-w_{j_{1}k}w_{j_{2}k}$, 
		the maximum modulus principal yields that
		\begin{align*}
			\Norm{\Omega_{k,j_{1},j_{2}}}_{\infty}
			&\leq\max\left\{\Norm{\Delta^{0}M_{f_{k,j_{1},j_{2}}}\Delta^{-0}}_{\infty},\Norm{\Delta^{1/2}M_{f_{k,j_{1},j_{2}}}\Delta^{-1/2}}_{\infty}\right\}\\
			&\lesssim\Norm{f_{k,j_{1},j_{2}}}_{\infty}+\Norm{f_{k,j_{1},j_{2}}}_{\dot{W}^{1,Q}}.
		\end{align*}
		It is easy to see
		\begin{align*}
			\Norm{f_{k,j_{1},j_{2}}}_{\infty}
			\leq\Norm{w-v}_{\lgm{n_{1}}}(\Norm{w}_{\lgm{n_{1}}}+\Norm{v}_{\lm{n_{1}}}).
		\end{align*}
		Since $X_{l}(v_{j_{1}k}v_{j_{2}k})=0$, it follows that
		\begin{align*}
			\Norm{f_{k,j_{1},j_{2}}}_{\dot{W}^{1,Q}}^{Q}
			=\sum_{l=1}^{n_{1}}\Norm{X_{l}(w_{j_{1}k}w_{j_{2}k})}_{\lp{Q}}^{Q}
			&\lesssim\Norm{w}_{\lgm{n_{1}}}^{Q}\sum_{l=1}^{n_{1}}
			(\Norm{X_{l}w_{j_{1}k}}_{\lp{Q}}^{Q}+\Norm{X_{l}w_{j_{2}k}}_{\lp{Q}}^{Q})\\
			&\lesssim\Norm{w}_{\lgm{n_{1}}}^{Q}\Norm{w}_{\sgm}^{Q}.
		\end{align*}
		Therefore
		\begin{align}\label{sobolevo}
			\Norm{\Omega_{k,j_{1},j_{2}}}_{\infty}
			\lesssim\Norm{w-v}_{\lgm{n_{1}}}(\Norm{w}_{\lgm{n_{1}}}+\Norm{v}_{\lm{n_{1}}})
			+\Norm{w}_{\lgm{n_{1}}}\Norm{w}_{\sgm}.
		\end{align}

		Denote
		\begin{align*}
			\Phi_{j}=\Delta^{-1/4}X_{j}\Delta^{-1/4}\quad{\rm and}\quad\Psi_{j}=\Delta^{1/4}X_{j}\Delta^{-3/4}.
		\end{align*}
		By \eqref{0-order-bd}, it follows that $\Phi_{j}$ and $\Psi_{j}$ are bounded on $\lp{2}$.
		A simple calculation shows that
		\begin{align*}
			\Delta^{-1/4}(L_{w}-L_{v})\Delta^{-3/4}
			&=\sum_{k=1}^{n_{1}}\sum_{j_{1},j_{2}=1}^{n_{1}}
			\Delta^{-1/4}X_{j_{1}}M_{v_{j_{1}k}v_{j_{2}k}-w_{j_{1}k}w_{j_{2}k}}X_{j_{2}}\Delta^{-3/4}\\
			&=\sum_{k=1}^{n_{1}}\sum_{j_{1},j_{2}=1}^{n_{1}}\Phi_{j_{1}}\Omega_{k,j_{1},j_{2}}\Psi_{j_{2}}.
		\end{align*}
		Thus, by triangle inequality and \eqref{sobolevo},
		\begin{align*}
			\Norm{\Delta^{-1/4}(L_{w}-L_{v})\Delta^{-3/4}}_{\infty}
			&\leq\sum_{k=1}^{n_{1}}\sum_{j_{1},j_{2}}^{n_{1}}\Norm{\Phi_{j_{1}}}_{\infty}
			\Norm{\Omega_{k,j_{1},j_{2}}}_{\infty}\Norm{\Psi_{j_{2}}}_{\infty}\\
			&\lesssim \Norm{w-v}_{\lgm{n_{1}}}(\Norm{w}_{\lgm{n_{1}}}+\Norm{v}_{\lm{n_{1}}})
			+\Norm{w}_{\lgm{n_{1}}}\Norm{w}_{\sgm}.
		\end{align*}
		This is the desired result.
	\end{proof}

	The following lemma is established in \cite[Thm 3.5.1]{DAO4}.
	
	\begin{lemma}\label{farc43}
		Let $\H$ be a Hilbert space and let $A,B\in\B(\H)$ be positive self-adjoint operators with the same domains and trivial kernels. Then
		\begin{align*}
			\Norm{(A^{1/2}-B^{1/2})A^{-1/2}}_{\infty}
			\leq c_{abs}\Norm{B^{-1/4}(B-A)A^{-3/4}}_{\infty}.
		\end{align*}
	\end{lemma}

	\begin{proposition}\label{rieszquotient}
		Assume that $a:\GG\rightarrow\at$ is a bounded smooth function with bounded (pointwise) inverse. For $y\in\GG$, there is a constant $C$ depending on $\GG$ such that
		\begin{align*}
			\Norm{\quasiriesz{a}{k}-\quasiriesz{a(y)}{k}}_{\infty}
			\leq C\Norm{a_{(1)}\rp}_{\lgm{n_{1}}}\Big(1+\Norm{a_{(1)}}_{\lgm{n_{1}}}^{2}\Norm{a_{(1)}\rp}_{\lgm{n_{1}}}^{2}\Big)\\
			\cdot\Big(\Norm{a_{(1)}-a_{(1)}(y)}_{\lgm{n_{1}}}+\Norm{a_{(1)}}_{\sgm}\Big),
		\end{align*}
		where $a_{(1)}$ is the first block of $a$ as in \eqref{diagonal}.
	\end{proposition}
	\begin{proof}
		Write
		\begin{align*}
			\quasiriesz{a}{k}-\quasiriesz{a(y)}{k}
			=(X_{k}^{a}-X_{k}^{a(y)})L_{a}^{-1/2}+X_{k}^{a(y)}(L_{a}^{-1/2}-L_{a(y)}^{-1/2})
			=:J_{1}+J_{2}.
		\end{align*}
		
		By \eqref{diagonal}, set
		\begin{align*}
			a=\diag(a_{(1)},a_{rest})=\diag\left((a_{jk})_{j,k=1}^{n_{1}},a_{rest}\right).
		\end{align*}
		It follows that $\displaystyle X_{k}^{a}-X_{k}^{a(y)}=\sum_{j=1}^{n_{1}}M_{a_{jk}-a_{jk}(y)}X_{j}$.
		Employing H\"{o}lder's inequality, Lemma \ref{laplacinbd} and triangle inequality, we obtain
		\begin{align*}
			\Norm{J_{1}}_{\infty}
			&\leq\Norm{(X_{k}^{a}-X_{k}^{a(y)})\Delta^{-1/2}}_{\infty}
			\Norm{\Delta^{1/2}L_{a}^{-1/2}}_{\infty}\\
			&\leq\Norm{a_{(1)}\rp}_{\lgm{n_{1}}}\sum_{j=1}^{n_{1}}
			\Norm{M_{a_{jk}-a_{jk}(y)}}_{\infty}\Norm{X_{j}\Delta^{-1/2}}\\
			&\leq C\Norm{a_{(1)}\rp}_{\lgm{n_{1}}}\Norm{a_{(1)}-a_{(1)}(y)}_{\lgm{n_{1}}}.
		\end{align*}
		This is the estimate for $J_{1}$.
		
		For $J_{2}$, using Lemma \ref{laplacian-bd-14}, one has
		\begin{align}\label{LD41}
			\Norm{L_{a}^{-1/4}\Delta^{1/4}}_{\infty}\leq\Norm{a_{(1)}\rp}_{\lgm{n_{1}}}^{1/2}.
		\end{align}
		According to Lemma \ref{laplacinbd} and Lemma \ref{laplacianbdd}, the Hadamard's three lines lemma \cite[Lemma 1.3.5]{G2014c} implies that
		\begin{align}\label{LD43}
			\Norm{\Delta^{3/4}L_{a(y)}^{-3/4}}_{\infty}
			\leq C\Norm{a_{(1)}\rp}_{\lgm{n_{1}}}^{3/2}.
		\end{align}
		Write
		\begin{align*}
			L_{a}^{-1/4}(L_{a}-L_{a(y)})L_{a(y)}^{-3/4}
			=L_{a}^{-1/4}\Delta^{1/4}\cdot\Delta^{-1/4}(L_{a}-L_{a(y)})\Delta^{-3/4}
			\cdot\Delta^{3/4}L_{a(y)}^{-3/4}.
		\end{align*}
		By \eqref{LD41}, \eqref{LD43} and Lemma \ref{w1qa-b},
		\begin{align}\label{LL-sobolev}
			\nonumber\Norm{L_{a}^{-1/4}(L_{a}-L_{a(y)})L_{a(y)}^{-3/4}}_{\infty}
			&\leq 
			C\Norm{a_{(1)}\rp}_{\lgm{n_{1}}}^{2}\Norm{\Delta^{-1/4}(L_{a}-L_{a(y)})\Delta^{-3/4}}_{\infty}\\
			&\leq C\Norm{a_{(1)}\rp}_{\lgm{n_{1}}}^{2}\Norm{a_{(1)}}_{\lgm{n_{1}}}
			\Big(\Norm{a_{(1)}-a_{(1)}(y)}_{\lgm{n_{1}}}+\Norm{a_{(1)}}_{\sgm}\Big).
		\end{align}
		Observing
		\begin{align*}
			J_{2}=X_{k}^{a(y)}\Delta^{-1/2}\cdot \Delta^{1/2}L_{a}^{-1/2}
			\cdot (L_{a(y)}^{1/2}-L_{a}^{1/2})L_{a(y)}^{-1/2},
		\end{align*}
		it follows from H\"{o}lder's inequality, Lemma \ref{laplacinbd}, Lemma \ref{farc43}  and inequality \eqref{LL-sobolev} that
		\begin{align*}
			\Norm{J_{2}}_{\infty}
			&\leq\Norm{X_{k}^{a(y)}\Delta^{-1/2}}_{\infty}\Norm{\Delta^{1/2}L_{a}^{-1/2}}_{\infty}
			\Norm{(L_{a(y)}^{1/2}-L_{a}^{1/2})L_{a(y)}^{-1/2}}_{\infty}\\
			&\leq C
			\Norm{a_{(1)}}_{\lgm{n_{1}}}\Norm{a_{(1)}\rp}_{\lgm{n_{1}}}
			\Norm{L_{a}^{-1/4}(L_{a}-L_{a(y)})L_{a(y)}^{-3/4}}_{\infty}\\
			&\leq C\Norm{a_{(1)}}_{\lgm{n_{1}}}^{2}\Norm{a_{(1)}\rp}_{\lgm{n_{1}}}^{3}
			\Big(\Norm{a_{(1)}-a_{(1)}(y)}_{\lgm{n_{1}}}+\Norm{a_{(1)}}_{\sgm}\Big).
		\end{align*}
		This estimate and the estimate of $J_{1}$ yield the desired result.
	\end{proof}
	\begin{corollary}\label{rieszquotientcorr}
		Let $\{a_i\}_{i\in I}$ be a family of smooth functions $\GG\to\at$ indexed by a set $I$, and suppose that $\{(a_i)_{(1)} : i\in I\}$ is a bounded subset of $\lgm{n_1}$. Then, there is $C>0$ such that, for all $i\in I$ and all $\gamma\in\GG$,
		$$\norm{R^{a_i}_k-R^{a_i(\gamma)}_k}_\infty\leq C\left (\norm{(a_i)_{(1)}-(a_i(\gamma))_{(1)}}_{\lgm{n_1}}+\norm{(a_i)_{(1)}}_{\sgm}\right ).$$
		\begin{proof}
			Since the $C$ in Proposition \ref{rieszquotient} only depends on $\GG$ and the $\lgm{n_1}$-norms of the $a(\cdot, \gamma, t)_{(1)}$ and $a(\cdot, \gamma, t)_{(1)}^{-1}$ are uniformly bounded, the Corollary follows directly.
		\end{proof}
	\end{corollary}

	\subsection{Approximation property}\label{Approximation property}

	A lattice $S$ of a Lie group $G$ is a discrete subgroup such that the quotient space $G/S$ is compact. Not every Lie group admits a lattice, even in the nilpotent case. By, \cite[Theorem 2.12]{R-book1972} a simply connected nilpotent Lie group $G$ admits a lattice if and only if its Lie algebra admits a basis for which the structure constants are rational numbers. Therefore, to prove Theorem \ref{symbol on manifold} on a general $\GG$-filtered manifold, we cannot rely on the assumption that $\GG$ admits a lattice. Therefore, we use the following lemma which follows easily from e.g., \cite[Lemma 5.7.5]{FR-Book} (see also \cite[Lemma 3.3]{LMSZ2024} for a generalisation on regular metric spaces).
	
	\begin{lemma}\label{bd multiplicity covering}
		There exists a countable set $\Gamma\subset\GG$, a family $\{\eta^\gamma\}_{\gamma\in\Gamma}$ of non-negative, smooth, compactly supported functions, and a constant $M>0$ such that
		\begin{enumerate}
			\item \itemcase{T}he balls $B(\gamma,1)$, for ${\gamma\in\Gamma}$, cover $\GG$;
			\item \itemcase{F}or all $\gamma\in\Gamma$, $\eta^\gamma$ is supported in $B(\gamma, 1)$;
			\item \itemcase{T}he functions $(\eta^\gamma)^2$, for $\gamma\in\Gamma$, form a partition of unity;
			\item \itemcase{F}or all $\gamma\in\Gamma$, one has $\eta^\gamma\equiv 1$ on a neighborhood of $\gamma$, and
			\item \itemcase{F}or each $C\geq 1$, no point of $\GG$ is contained in more than $(MC)^Q$ of the sets $B(\gamma, C)$, for ${\gamma\in\Gamma}$.
		\end{enumerate}
	\end{lemma}
	We will continue to refer to these $\Gamma$, $\eta^\gamma$ and $M$ throughout the remainder of this subsection. Moreover, if $\epsilon>0$, we will write $\Gamma_\epsilon=\delta_\epsilon(\Gamma)$ and if also $\gamma\in\Gamma_\epsilon$, we will write $\eta^\gamma_\epsilon=\eta^{\delta_\epsilon^{-1}(\gamma)}\circ\delta_{\epsilon}^{-1}$. Finally, we will denote the matrix unit in $\gln$ by $\unit_{n\times n}$.
	
	\begin{lemma}\label{ww-sobolev}
		If $n\geq1$ and $w:\GG\rightarrow\gln$ is a bounded smooth function, then there is $C>0$ such that
		\begin{align*}
			\Norm{w(x)-w(y)}_{\lm{n}}\leq C \Norm{w}_{\wkp{\infty}}\rho(y\rp x).
		\end{align*}
		for all $x,y\in\GG$.
	\end{lemma}
	\begin{proof}
		By \cite[Proposition 5.4]{F1975}, for each $1\leq k,j\leq n$, there is $C_{kj}>0$ such that
		\begin{equation}\label{meanvalue}
			\sup_{a,b\in\GG} |w_{kj}(ab)-w_{kj}(a)|\leq C_{kj}\sum_{j=1}^{n_1}\norm{X_j w_{kj}}_\infty\rho(b),
		\end{equation}
		for all $1\leq k,j\leq n$. Let $x,y\in\GG$. Since all norms involved are translation invariant, we may assume that neither $y$ nor $y^{-1}x$ are $o$. Therefore, applying \eqref{meanvalue} with $a=y$ and $b=y^{-1}x$, we get
		\begin{equation}\label{meanvalueap}
			|w_{kj}(x)-w_{kj}(y)|\leq C_{kj}\sum_{j=1}^{n_1}\norm{X_j w_{kj}}_\infty\rho(y^{-1}x),
		\end{equation}
		for all $1\leq k,j\leq n$. Now taking the sup over $k$ and $j$, and writing $\smash{C=\sup_{k,j} C_{kj}}$, we get
		\begin{align*}
			\Norm{w(x)-w(y)}_{\lm{n}}&=\sup_{1\leq k,j\leq n} |w_{kj}(x)-w_{kj}(y)|\leq \sup_{1\leq k,j\leq n} C_{kj}\sum_{j=1}^{n_1}\norm{X_j w_{kj}}_\infty\rho(y\rp x)\\
			&\leq C \sup_{1\leq k,j\leq n}\norm{X_j w_{kj}}\rho(y\rp x)\\
			&\leq C \norm{w}_{\wkp{\infty}}\rho(y\rp x),
		\end{align*}
		as desired.
	\end{proof}

	\begin{lemma}\label{log-exp-xk}
		Let $n\geq1$ and let $w:\Omega\subset\GG\rightarrow\gln$ be a bounded smooth function.
		\begin{enumerate}[\rm (1)]
			\item\label{log} \itemcase{I}f $\Norm{w(x)-\unit_{n\times n}}_{\lm{n}}\leq 1/2$ for $x\in\Omega$, then 
			\begin{align*}
				\Norm{\log(w(x))}_{\lm{n}}&\leq2\Norm{w(x)-\unit_{n\times n}}_{\lm{n}},\,\text{and}\\
				\Norm{X_{j}\log(w(x))}_{\lm{n}}&\leq2\Norm{X_{j}w(x)}_{\lm{n}}.
			\end{align*}
			
			\item\label{exp} \itemcase{I}f $1\leq p\leq\infty$, then
			$$\Norm{X_{j}\Exp(w)}_{\mlp{p}}\leq \exp(\Norm{w}_{\mlp{p}})\Norm{X_{j}w}_{\mlp{p}}.$$
		\end{enumerate}
	\end{lemma}
	\begin{proof}
		Let us prove \eqref{log} at first. For $\Norm{w(x)-\unit_{n\times n}}_{\lm{n}}<1$, one has the definition, in the sense of $\lm{n}$-norm,
		\begin{align*}
			\log(w(x))=\sum_{k=0}^{\infty}\frac{(-1)^{k}}{k+1}(w(x)-\unit_{n\times n})^{k+1}.
		\end{align*}
		It follows from triangle inequality and $\Norm{w(x)-\unit_{n\times n}}_{\lm{n}}\leq\frac{1}{2}$ that
		\begin{align*}
			\Norm{\log(w(x))}_{\lm{n}}
			&\leq\sum_{k=0}^{\infty}\Norm{w(x)-\unit_{n\times n}}_{\lm{n}}^{k+1}\\
			&\leq\sum_{k=0}^{\infty}\frac{1}{2^{k}}\Norm{w(x)-\unit_{n\times n}}_{\lm{n}}
			\leq2\Norm{w(x)-\unit_{n\times n}}_{\lm{n}}.
		\end{align*}
		This is the first inequality of \eqref{log}. 
		
		For the second inequality of \eqref{log}, employing the Leibniz rule, one has
		\begin{align*}
			X_{k}\Big((w(x)-\unit_{n\times n})^{k+1}\Big)
			=\sum_{l=0}^{k}(w(x)-\unit_{n\times n})^{l}\cdot X_{j}w(x)\cdot (w(x)-\unit_{n\times n})^{k-l}.
		\end{align*}
		Therefore, by triangle inequality,
		\begin{align*}
			\Norm{X_{k}\Big((w(x)-\unit_{n\times n})^{k+1}\Big)}_{\lm{n}}
			&\leq\sum_{l=0}^{k}\Norm{w(x)-\unit_{n\times n}}_{\lm{n}}^{k} \Norm{X_{j}w(x)}_{\lm{n}}\\
			&\leq\frac{k+1}{2^{k}}\Norm{X_{j}w(x)}_{\lm{n}}.
		\end{align*}
		Thus,
		\begin{align*}
			\Norm{X_{k}\log(w(x))}_{\lm{n}}
			\leq\sum_{k=0}^{\infty}\frac{1}{k+1}\Norm{X_{k}\Big((w(x)-\unit_{n\times n})^{k+1}\Big)}_{\lm{n}}
			\leq2\Norm{X_{j}w(x)}_{\lm{n}}.
		\end{align*}
		This is the proof of \eqref{log}.
		
		For \eqref{exp}, using the Leibniz rule again, one obtains
		\begin{align*}
			X_{j}\Exp(w)=\sum_{k=1}^{\infty}\frac{1}{k!}\sum_{l=0}^{k-1}w^{l}\cdot X_{j}w\cdot w^{k-1-l}.
		\end{align*}
		It follows from triangle inequality that
		\begin{align*}
			\Norm{X_{j}\Exp(w)}_{\mlp{p}}
			\leq\Norm{X_{j}w}_{\mlp{p}}\sum_{k=1}^{\infty}\frac{1}{(k-1)!}\Norm{w}_{\mlp{p}}^{k-1}
			=\exp(\Norm{w}_{\mlp{p}})\Norm{X_{j}w}_{\mlp{p}}.
		\end{align*}
		This completes the proof.
	\end{proof}

	\begin{lemma}\label{approximation-inner}
		Let $n\geq1$ and $\gamma\in\GG$. Suppose that $w:\GG\rightarrow\gln$ is a smooth function which is constant outside of some ball. 		There exists $C_{w}>0$ such that, if $0<\epsilon<\epsilon_w$, the formula
		\begin{align*}
			\xi^{\gamma}_{w}(x)=\log\Big(w(x)w(\gamma)\rp\Big)
		\end{align*}
		gives a well-defined, smooth function $\xi^{\gamma}_{w}$ on $B(\gamma,\epsilon)$ satisfying
		\begin{align*}
			\Norm{\xi^{\gamma}_{w}(x)}_{\lm{n}}\leq C_{w}\epsilon\leq 1
			\quad\text{and}\quad
			\Norm{(X_{j}\xi^{\gamma}_{w})(x)}_{\lm{n}}\leq 2C_{w}
		\end{align*}
		for all $x\in B(\gamma, \epsilon)$.
	\end{lemma}
	\begin{proof} Let $C$ be as in Lemma \ref{ww-sobolev}. Let $C_{w}=\max\{C,1\}\cdot \Norm{w}_{\wkp{\infty}}\Norm{w\rp}_{\lgm{n}}$. Set $\epsilon_w=(2C_w)^{-1}$.

		By Lemma \ref{ww-sobolev},
		\begin{align}\label{wwrp-wkp}
			\nonumber\Norm{w(x)w(\gamma)\rp-\unit_{n\times n}}_{\lm{n}}
			&\leq\Norm{w(x)-w(\gamma)}_{\lm{n}}\Norm{w(\gamma)\rp}_{\lm{n}}\\
			&\leq C\Norm{w}_{\wkp{\infty}}\Norm{w\rp}_{\lgm{n}}\rho(x\rp \gamma)\\
			\nonumber			&\leq C_w\epsilon.
		\end{align}
		
		If $0<\epsilon<\epsilon_w$ and $x\in B(\gamma,\epsilon)$, then one obtains
		\begin{align*}
			\Norm{w(x)w(\gamma)\rp-\unit_{n\times n}}_{\lm{n}}\leq\frac{1}{2}.
		\end{align*}
		Therefore, $w(x)w(\gamma)^{-1}$ is inside the radius of convergence for the power series defining $\log$ around $\unit_{n\times n}$, and hence $\log\big(w(x)w(\gamma)\rp\big)$ is well defined for all $x\in B(\gamma,\epsilon)$. Moreover, it is clearly smooth.
		
		By Lemma \ref{log-exp-xk} and \eqref{wwrp-wkp},
		\begin{align*}
			\Norm{\xi^{\gamma}_{w}(x)}_{\lm{n}}
			\leq2\Norm{w(x)w(\gamma)\rp-\unit_{n\times n}}_{\lm{n}}
			\leq 2C_{w}\epsilon\leq1.
		\end{align*}
		By Lemma \ref{log-exp-xk} and H\"{o}lder's inequality,
		\begin{align*}
			\Norm{(X_{j}\xi^{\gamma}_{w})(x)}_{\lm{n}}
			\leq2\Norm{(X_{j}w)(x)w(\gamma)\rp}_{\lm{n}}
			\leq2\Norm{w}_{\wkp{\infty}}\Norm{w\rp}_{\lgm{n}}\leq 2C_w.
		\end{align*}
		This yields the desired result.
	\end{proof}
	
	\begin{lemma}\label{approximation}
		Let $n\geq 1$, let $w:\GG\rightarrow\gln$ be a smooth function, constant outside of a ball. Let $\epsilon_w$ be as in Lemma \ref{approximation-inner}. Define $\Lambda=\{(\epsilon,\gamma) : \epsilon\in(0,\epsilon_w/2], \gamma\in\Gamma_\epsilon\}$. There is a family $\{w^\gamma_\epsilon\}_{(\epsilon,\gamma)\in\Lambda}$ of smooth functions $w^\gamma_\epsilon:\GG\to\gln$, such that
		\begin{enumerate}
			\item\label{equal} $w^\gamma_\epsilon$ agrees with $w$ on $B(\gamma, \epsilon)$, for all $(\epsilon,\gamma)\in\Lambda$;
			\item\label{constant} $w^\gamma_\epsilon$ equals $w(\gamma)$ on $\GG\setminus B(\gamma, 2\epsilon)$, for all $(\epsilon,\gamma)\in\Lambda$;
			\item\label{appro-gamma} $\smash{\lim_{\epsilon\downarrow 0}(\sup_{\gamma\in\Gamma_\epsilon} ||w^\gamma_\epsilon-w(\gamma)||_{\lgm{n}})=0}$;
			\item\label{appro-sobolev}  $\smash{\lim_{\epsilon\downarrow 0}(\sup_{\gamma\in\Gamma_\epsilon} ||w^\gamma_\epsilon||_{\sgn})=0}$, and
			\item\label{bound} $\{(w^\gamma_\epsilon)^{j} : (\epsilon,\gamma)\in\Lambda, j\in\{-1,1\}\}$ is a bounded subset of $\lgm{n}$.
		\end{enumerate}
	\end{lemma}
	\begin{proof} Let $\theta\in\cci$ be such that $0\leq\theta\leq 1$ on $\mathbb{G}.$ Suppose that $\theta=1$ on $B(o,1)$ and that $\theta$ is supported on $B(o,2).$ Write $\theta_{\gamma,\epsilon}(x)=\theta(\delta_{\epsilon\rp}(\gamma^{-1}x))$ for any $\gamma,x\in\GG$ and $\epsilon>0$. Let $(\epsilon, \gamma)\in\Lambda$. Define $w^\gamma_\epsilon:B(\gamma, 2\epsilon)\to\gln$ by
		\begin{align*}
			w^\gamma_\epsilon(x)=\Exp\Big(\theta_{\gamma,\epsilon}(x)\xi^{\gamma}_{w}(x)\Big)\cdot w(\gamma),
		\end{align*}
		where $\xi^{\gamma}_{w}$ is given in Lemma \ref{approximation-inner}. It follows from Lemma \ref{approximation-inner} that $w^\gamma_\epsilon$ is a well-defined and smooth. Note that $\theta_{\gamma,\epsilon}=1$ on $B(\gamma,\epsilon)$ and, therefore, 
		\begin{align*}
			w^\gamma_\epsilon(x)=\Exp\Big(\!\log\Big(w(x)w(\gamma)\rp\Big)\Big)\cdot w(\gamma)=w(x)
		\end{align*}
		for every $x\in B(\gamma,\epsilon).$ This is the assertion \eqref{equal}. Moreover, since $\theta_{\gamma,\epsilon}=0$ outside $B(\gamma,2\epsilon),$, we may extend $w^\gamma_\epsilon$ to all of $\GG$ by setting its value outside of $B(\gamma, 2\epsilon)$ to be $w(\gamma)$, giving \eqref{constant}. Given $v\in\gln$, it is elementary that
		\begin{align*}
			\Norm{\Exp(v)-\unit_{n\times n}}_{\lm{n}}
			\leq\Norm{v}_{\lm{n}}\sum_{k=1}^{\infty}\frac{1}{k!}\Norm{v}_{\lm{n}}^{k-1}
			\leq \exp\big(\Norm{v}_{\lm{n}}\big)\Norm{v}_{\lm{n}}.
		\end{align*}
		Letting $x\in B(\gamma,2\epsilon)$, it follows from H\"{o}lder's inequality, Lemma \ref{approximation-inner} and $0\leq\theta\leq1$ that
		\begin{align*}
			\Norm{\Exp\big(\theta_{\gamma,\epsilon}(x)\xi^{\gamma}_{w}(x)\big)-\unit_{n\times n}}_{\lm{n}}
			&\leq \exp(\Norm{\theta}_{\infty})\Norm{\theta}_{\infty}\Norm{\xi^{\gamma}_{w}(x)}_{\lm{n}}\leq6C_{w}\epsilon.
		\end{align*}
		Thus, for every $x\in B(\gamma,2\epsilon),$
		\begin{align*}
			\Norm{w^\gamma_\epsilon(x)-w(\gamma)}_{\lm{n}}
			&\leq\Norm{\Exp\big(\theta_{\gamma,\epsilon}(x)\xi^{\gamma}_{w}(x)\big)-\unit_{n\times n}}_{\lm{n}}
			\Norm{w(\gamma)}_{\lm{n}}\\
			&\leq6\Norm{w}_{\lgm{n}}C_{w}\epsilon.
		\end{align*}
		For every $x\notin B(\gamma,2\epsilon),$ we have $w^\gamma_\epsilon(x)=w(\gamma).$ Thus,
		\begin{align*}
			\Norm{w^\gamma_\epsilon(x)-w(\gamma)}_{\lm{n}}&\leq6\Norm{w}_{\lgm{n}}C_{w}\epsilon,\quad x\in\mathbb{G}.
		\end{align*}
		This yields the assertion \eqref{appro-gamma}.

		Let $1\leq j\leq n_1$. By the Leibniz rule, 
		\begin{align*}
			X_{j}(\theta_{\gamma,\epsilon}\xi^{\gamma}_{w})=
			X_{j}(\theta_{\gamma,\epsilon})\xi^{\gamma}_{w}+\theta_{\gamma,\epsilon}X_{j}(\xi^{\gamma}_{w}).
		\end{align*}
		Observe that $X_{j}$ is left-invariant and $\delta_{r}$-homogeneous of degree $1$ on $\GG$.
		Then
		\begin{align*}
			(X_{j}\theta_{\gamma,\epsilon})(x)=\epsilon\rp (X_{j}\theta)(\delta_{\epsilon^{-1}}(\gamma^{-1}x)).
		\end{align*}
		Since the support of $\theta_{\gamma,\epsilon}$ is contained in $B(\gamma,2\epsilon)$,
		it follows from triangle inequality, H\"{o}lder's inequality and Lemma \ref{approximation-inner} that
		\begin{align*}
			\Norm{X_{j}(\theta_{\gamma,\epsilon}\xi^{\gamma}_{w})}_{\mlp{Q}}
			&\leq\Norm{X_{j}\theta_{\gamma,\epsilon}}_{L_{Q}(B(\gamma,2\epsilon))}
			\Norm{\xi^{\gamma}_{w}}_{L_{\infty}(B(\gamma,2\epsilon);\MM_{n})}+\\
			&\quad+\Norm{\theta_{\gamma,\epsilon}}_{L_{Q}(B(\gamma,2\epsilon))}
			\Norm{X_{j}\xi^{\gamma}_{w}}_{L_{\infty}(B(\gamma,2\epsilon);\MM_{n})}\\
			&\leq\Norm{X_{j}\theta}_{\lp{Q}}\cdot 2C_{w}\epsilon+\epsilon\Norm{\theta}_{\lp{Q}}\cdot 2C_{w}\\
			&=:C_{\theta,w}\epsilon.
		\end{align*}
		Employing Lemma \ref{approximation-inner} again,
		\begin{align*}
			\Norm{\theta_{\gamma,\epsilon}\xi^{\gamma}_{w}}_{\mlp{Q}}
			\leq\Norm{\theta_{\gamma,\epsilon}}_{L^{Q}(B(\delta_{\epsilon}(\gamma),\epsilon))}
			\Norm{\xi^{\gamma}_{w}}_{L_{\infty}(B(\delta_{\epsilon}(\gamma),\epsilon);\MM_{n})}
			\leq \Norm{\theta}_{\lp{Q}}C_{w}\epsilon^{2}.
		\end{align*}
		Therefore, applying Lemma \ref{log-exp-xk} to $p=Q$, we have
		\begin{align*}
			\Norm{X_{j}\Exp\big(\theta_{\gamma,\epsilon}\xi^{\gamma}_{w}\big)}_{\mlp{Q}}
			&\leq \exp\Big(\Norm{\theta_{\gamma,\epsilon}\xi^{\gamma}_{w}}_{\mlp{Q}}\Big) \Norm{X_{j}(\theta_{\gamma,\epsilon}\xi^{\gamma}_{w})}_{\mlp{Q}}\\
			&\leq \exp\big(\Norm{\theta}_{\lp{Q}}C_{w}\epsilon^{2}\big)C_{\theta,w}\epsilon.
		\end{align*}
		The assertion \eqref{appro-sobolev} now easily follows.

		Finally, by Lemma \ref{approximation-inner},
		\begin{align*}
			\norm{w^\gamma_\epsilon}_{\lgm{n}}&\leq \norm{\Exp(\theta_{\gamma,\epsilon}\xi^\gamma_w)}_{\lgm{n}}\norm{w(\gamma)}_{\lm{n}}\\
			&\leq\exp\Big(\norm{\theta_{\gamma,\epsilon}\xi^\gamma_w}_{\lgm{n}}\Big)\norm{w(\gamma)}_{\lm{n}}\\
			&\leq \exp\Big(\norm{\theta}_\infty \norm{\xi^\gamma_w}_{\lgm{n}}\Big)\norm{w}_{\lgm{n}}\\
			&\leq \exp(\norm{\theta}_\infty)\norm{w}_{\lgm{n}},
		\end{align*}
		and similarly,
		\begin{align*}
			\norm{(w^\gamma_\epsilon)^{-1}}_{\lgm{n}}&\leq \norm{w(\gamma)^{-1}}_{\lm{n}}\norm{\Exp(\theta_{\gamma,\epsilon}\xi^\gamma_w)^{-1}}_{\lgm{n}}\\
			&=\norm{w(\gamma)^{-1}}_{\lm{n}}\norm{\Exp(-\theta_{\gamma,\epsilon}\xi^\gamma_w)}_{\lgm{n}}\\
			&\leq \exp(-\norm{\theta}_\infty)||w^{-1}||_{\lgm{n}}.
		\end{align*}
		Since these bounds are independent of $(\epsilon,\gamma)$, assertion \ref{bound} follows.
	\end{proof}

	\begin{proposition}\label{limitzero}
		Let $a: \GG\rightarrow\at$ be a smooth function, constant outside of a ball. Let $\{a^\gamma_\epsilon\}_{(\epsilon,\gamma)\in\Lambda}$ be as in Lemma \ref{approximation}. For all $k\in \{1,\,\dotsc, n_1\}$, $\epsilon>0$ and $\gamma\in\Gamma_\epsilon$, define
		$$A^a_{k,\epsilon}(\gamma)=M_{\eta^\gamma_\epsilon} R^{a^\gamma_\epsilon}_kM_{\eta^\gamma_\epsilon},\qquad R^a_{k,\epsilon}(y)=M_{\eta^\gamma_\epsilon} R^{a(\gamma)}_kM_{\eta^\gamma_\epsilon}.$$
		Then, for all $k\in \{1,\,\dotsc, n_1\}$ and $\epsilon>0$, the sums $A^a_{k,\epsilon}:=\sum_{\gamma\in\Gamma_\epsilon} A^a_{k,\epsilon}(\gamma)$ and $R^a_{k,\epsilon}:=\sum_{\gamma\in\Gamma_\epsilon} R^a_{k,\epsilon}(\gamma)$ converge strongly in $B(L_2(\GG))$, and one has
		$$\lim_{\epsilon\downarrow 0} \norm{A^a_{k,\epsilon}-R^a_{k,\epsilon}}_\infty=0$$
	\end{proposition}
	\begin{proof}
		First of all, the strong convergence of the sums follows from the local finiteness of the $\eta^\gamma_\epsilon$. Let $k\in \{1,\,\dotsc, n_1\}$ and $\epsilon>0$. There is $N\in\NN$ such that, for all $x\in\GG$, one has $\eta^\gamma_\epsilon(x)\neq 0$ for at most $N$ values of $\gamma\in\Gamma_\epsilon$. Since also $0\leq \eta^\gamma_\epsilon\leq 1$ for all $\gamma\in\Gamma_\epsilon$, if $b\in\ell_2(\Gamma_\epsilon)$, by Cauchy-Schwartz we have
		\begin{align*}
			\left|\sum_{\gamma\in\Gamma_\epsilon}\eta^\gamma_\epsilon(x)b_\gamma\right|^2\leq \sum_{\gamma\in\Gamma_\epsilon}|\eta^\gamma_\epsilon(x)|^2 \sum_{\gamma\in\Gamma_\epsilon}|b_\gamma|^2\leq N\sum_{\gamma\in\Gamma_\epsilon}|b_\gamma|^2.
		\end{align*}
		In particular, if $f\in\cci$, then
		\begin{align}
			\nonumber \norm{(A^a_{k,\epsilon}-R^a_{k,\epsilon})f}_2^2&=\int_G \left|\sum_{\gamma\in\Gamma_\epsilon}\eta^\gamma_\epsilon(x)((R^{a^\gamma_\epsilon}_k-R^{a(\gamma)}_{k}) M_{\eta^\gamma_\epsilon} f)(x)\right|^2d x\\
			\nonumber    &\leq N\!\int_G \sum_{\gamma\in\Gamma_\epsilon}\left|(R^{a^\gamma_\epsilon}_{k}-R^{a(\gamma)}_{k}) M_{\eta^\gamma_\epsilon} f\right|^2d \mu\\
			\label{lzfatou}&\leq N\!\sum_{\gamma\in\Gamma_\epsilon}\int_G \left|(R^{a^\gamma_\epsilon}_{k}-R^{a(\gamma)}_{k}) M_{\eta^\gamma_\epsilon} f\right|^2d\mu\\
			\nonumber&\leq N\sup_{\gamma\in\Gamma_\epsilon}\norm{R^{a^\gamma_\epsilon}_{k}-R^{a(\gamma)}_{k}}^2_\infty \biggl(\sum_{\gamma\in\Gamma_\epsilon}\norm{M_{\eta^\gamma_\epsilon} f}_2^2\biggr)\\
			\label{lzmonotone}&=N\sup_{\gamma\in\Gamma_\epsilon}\norm{R^{a^\gamma_\epsilon}_{k}-R^{a(\gamma)}_{k}}^2_\infty\int_G\biggl(\sum_{\gamma\in\Gamma_\epsilon}|\eta^\gamma_\epsilon|^2 |f|^2\biggr)\,d\mu\\
			\label{lzpartition}&=N\sup_{\gamma\in\Gamma_\epsilon}\norm{R^{a^\gamma_\epsilon}_{k}-R^{a(\gamma)}_{k}}^2_\infty\norm{f}_2^2
		\end{align}
		where \eqref{lzfatou} follows by Fatou's lemma, \eqref{lzmonotone} by monotone convergence, and \eqref{lzpartition} by the fact that $\{(\eta^\gamma_\epsilon)^2\}_{\gamma\in\Gamma_\epsilon}$ is a partition of unity on $\GG$. Therefore, by Lemma \ref{approximation}\eqref{bound}, Corollary \ref{rieszquotientcorr} of Proposition \ref{rieszquotient} applies, hence there is $C>0$ such that
		$$\norm{(A^a_{k,\epsilon}-R^a_{k,\epsilon})f}_2^2\leq N\sup_{\gamma\in\Gamma_\epsilon}C\left(\norm{a^\gamma_\epsilon-a(\gamma)}_{\lgm{n_1}}+\norm{a^\gamma_\epsilon}_{\sgm}\right)\norm{f}_2^2.$$
		Therefore, by \eqref{appro-gamma} and \eqref{appro-sobolev} of Lemma \ref{approximation},
		\begin{align*}
			\lim_{\epsilon\downarrow 0} \norm{A^a_{k,\epsilon}-R^a_{k,\epsilon}}_\infty\leq \sqrt{NC} \sup_{\norm{f}_2\leq 1}\left(\lim_{\epsilon\downarrow 0} \sup_{\gamma\in\Gamma_\epsilon}\norm{a^\gamma_\epsilon-a(\gamma)}_{\lgm{n_1}}+\lim_{\epsilon\downarrow 0} \sup_{\gamma\in\Gamma_\epsilon}\norm{a^\gamma_\epsilon}_{\sgm}\right)^{1/2}=0,
		\end{align*}
		where the supremum in $f$ is taken over $f\in C_c^\infty(\GG)$.
	\end{proof}

	\begin{proposition}\label{compactAk}
		Let $\epsilon>0$, let $k\in\{1, \hdots, n_1\}$, let $a$ and $A^{a}_{k,\epsilon}$ be as in Proposition \ref{limitzero} and let $\psi\in\cci$. Then
		\begin{align*}
			M_{\psi}(\quasiriesz{a}{k}-A^{a}_{k,\epsilon})M_{\psi}\in\K(\lp{2}).
		\end{align*}
	\end{proposition}
	
	\begin{proof}
		If $\gamma\in\Gamma_\epsilon$, then
		\begin{equation}\label{M eta psi}
			M_{(\eta^\gamma_\epsilon)^{2}}M_{\psi}R^a_k
			=M_{\eta^\gamma_\epsilon\psi}M_{\eta^\gamma_\epsilon}R^a_k=M_{\eta^\gamma_\epsilon\psi}\bigl(R^a_kM_{\eta^\gamma_\epsilon}+[M_{\eta^\gamma_\epsilon},R^a_k]\bigr).
		\end{equation}
		By construction, the set $I_\epsilon:=\{\gamma\in\Gamma_\epsilon : \eta^\gamma_\epsilon\psi\neq 0\}$ is finite, thus if we sum the right-hand side of \eqref{M eta psi} over $\gamma\in\Gamma_\epsilon$, only finitely many of the summands are non-zero, so the sum converges. Since also $\sum_{\gamma\in\Gamma_\epsilon}(\eta^\gamma_\epsilon)^2=1$, for any $L_2$-function $f$,
		\begin{align*}
			\biggl(M_{\psi}R^a_k\biggr)f&=\biggl(\sum_{\gamma\in\Gamma_\epsilon}(\eta^\gamma_\epsilon)^2\biggr)\biggl(M_{\psi}R^a_kf\biggr)=\biggl(\sum_{\gamma\in I_\epsilon}(\eta^\gamma_\epsilon)^2\psi \biggr)\biggl(R^a_kf\biggr)=\biggl(\sum_{\gamma\in I_\epsilon}M_{(\eta^\gamma_\epsilon)^{2}}M_{\psi}R^a_k\biggr)f,
		\end{align*}
		and therefore
		$$M_\psi R^a_k=\sum_{\gamma\in I_\epsilon}M_{\eta^\gamma_\epsilon\psi}(R^a_kM_{\eta^\gamma_\epsilon} + [M_{\eta^\gamma_\epsilon}, R^a_k]),$$
		hence
		\begin{align}
			\nonumber M_{\psi}(\quasiriesz{a}{k}-A^{a}_{k,\epsilon})M_{\psi}
   &=M_{\psi}\quasiriesz{a}{k}M_\psi-M_\psi A^{a}_{k,\epsilon}M_{\psi}\\
   \nonumber&=\biggl(\sum_{\gamma\in I_\epsilon}M_{\eta^\gamma_\epsilon\psi}\bigl(R^a_kM_{\eta^\gamma_\epsilon} + [M_{\eta^\gamma_\epsilon}, R^a_k]\bigr)\biggr)M_\psi - M_\psi A^{a}_{k,\epsilon}M_{\psi}\\
			\nonumber &=\sum_{\gamma\in I_\epsilon}M_{\eta^\gamma_\epsilon\psi}\bigl(R^a_kM_{\eta^\gamma_\epsilon} + [M_{\eta^\gamma_\epsilon}, R^a_k]\bigr)M_\psi-M_\psi \biggl(\sum_{\gamma\in\Gamma_\epsilon} M_{\eta^\gamma_\epsilon} R^{a^\gamma_\epsilon}_kM_{\eta^\gamma_\epsilon}\biggr)M_{\psi}\\
			\nonumber &=\sum_{\gamma\in I_\epsilon}M_{\eta^\gamma_\epsilon\psi}R^a_kM_{\eta^\gamma_\epsilon\psi} + \sum_{\gamma\in I_\epsilon}M_{\eta^\gamma_\epsilon\psi}[M_{\eta^\gamma_\epsilon}, R^a_k]M_\psi-\sum_{\gamma\in I_\epsilon} M_{\eta^\gamma_\epsilon\psi} R^{a^\gamma_\epsilon}_kM_{\eta^\gamma_\epsilon\psi}\\
			\label{cptsum}&=\sum_{\gamma\in I_\epsilon}M_{\eta^\gamma_\epsilon\psi}(R^a_k-R^{a^\gamma_\epsilon}_k)M_{\eta^\gamma_\epsilon\psi} + \sum_{\gamma\in I_\epsilon}M_{\eta^\gamma_\epsilon\psi}[M_{\eta^\gamma_\epsilon}, R^a_k]M_\psi.
		\end{align}
		By Proposition \ref{compactrieszquo}, the first sum appearing in in \eqref{cptsum} is compact, and by Proposition \ref{compactrieszcommu}, so too is the second.
	\end{proof}

	Let us identify the $C^{*}$-algebra $(C_{0}+\CC)(\GG)\otimes_{\rm min}\A_{2}$ with the algebra $(C_{0}+\CC)(\GG,\A_{2})$.
	\begin{proposition}\label{limitpsiriesz}
		Let $k\in\{1, \hdots, n_1\}$, let $a$ and $A^{a}_{k,\epsilon}$ be as in Proposition \ref{limitzero}, for all $\epsilon>0$, and let $\psi\in \cci$. Then $M_{\psi}\quasiriesz{a}{k,\epsilon}M_{\psi}\in\Pi$ for all $\epsilon>0$, and
		$$\lim_{\epsilon\downarrow0}\sym(M_{\psi}\quasiriesz{a}{k,\epsilon}M_{\psi})(x)=\psi(x)^{2} \quasiriesz{a(x)}{k}.$$
	\end{proposition}
	
	\begin{proof}
		For all $\epsilon>0$, let $I_{\epsilon}=\{\gamma\in\Gamma_\epsilon:\eta^\gamma_\epsilon\psi\neq 0\}$. Since the multipliers $M_{\psi\eta^\gamma_\epsilon}$ and all quasi-Riesz transforms are in $\Pi$, and since $I_\epsilon$ is finite, one has
		\begin{align}\label{finitesum}
			M_{\psi}\quasiriesz{a}{k,\epsilon}M_{\psi}=\sum_{\gamma\in I_{\epsilon}}
			M_{\psi\eta^\gamma_\epsilon}R_{k}^{a(\gamma)}M_{\psi\eta^\gamma_\epsilon}\in\Pi
		\end{align}
		for any $\epsilon>0$.
		
		Let $\epsilon>0$. We have
		\begin{align}
			\label{homo,def} \sym(M_{\psi}\quasiriesz{a}{k,\epsilon}M_{\psi})(x)
			&=\sum_{\gamma\in I_\epsilon}\psi(x)^{2}\eta^\gamma_\epsilon(x)^{2}R_{k}^{a(\gamma)}\\
   \nonumber&=\psi(x)^{2}\sum_{\gamma\in I_\epsilon}\eta^\gamma_\epsilon(x)^{2}(R_{k}^{a(x)}+R_{k}^{a(\gamma)}-R_{k}^{a(x)})\\
			\label{partunity}    &=\psi(x)^{2}\Bigl(R^{a(x)}_k+\sum_{\gamma\in I_\epsilon}\eta^\gamma_\epsilon(x)^{2}\bigl(R_{k}^{a(\gamma)}-R_{k}^{a(x)}\bigr)\Bigr)\\
			\label{equalonball}    &=\psi(x)^{2}\Bigl(R^{a(x)}_k+\sum_{\gamma\in I_\epsilon}\eta^\gamma_\epsilon(x)^{2}\bigl(R_{k}^{a^\gamma_\epsilon(\gamma)}-R_{k}^{a^\gamma_\epsilon(x)}\bigr)\Bigr)
		\end{align}
		Where \eqref{homo,def} follows from the definition of $\sym$ and the fact that it is a homomorphism, \eqref{partunity} follows from that $\sum_{\gamma\in\Gamma_\epsilon}(\eta^\gamma_\epsilon)^2=1$, and \eqref{equalonball} follows from that fact that, if $x\in B(\gamma, \epsilon)$, then $a(x)=a^\gamma_\epsilon(x)$, and if $x\notin B(\gamma,\epsilon)$, then $\eta^\gamma_\epsilon(x)=0$.
		
		By Lemma \ref{approximation}\eqref{bound}, we may apply Corollary \ref{rieszquotientcorr} Proposition \ref{rieszquotient} and obtain a $C>0$ such that, for all $\epsilon>0$ and all $\gamma\in\Gamma_\epsilon$, one has
		$$\norm{R^{a^\gamma_\epsilon}_k-R^{a^\gamma_\epsilon(\gamma)}_k}_\infty\leq C\left (\norm{(a^\gamma_\epsilon)_{(1)}-(a^\gamma_\epsilon(\gamma))_{(1)}}_{\lgm{n_1}}+\norm{(a^\gamma_\epsilon)_{(1)}}_{\sgm}\right ),$$
		and therefore, writing $\smash{N_1(\gamma, \epsilon):=||(a^\gamma_\epsilon)_{(1)}-(a^\gamma_\epsilon(\gamma))_{(1)}||_{\lgm{n_1}}}$ and $\smash{N_2(\gamma,\epsilon):=||(a^\gamma_\epsilon)_{(1)}||_{\sgm}}$, by \eqref{appro-gamma} and \eqref{appro-sobolev} of Lemma \ref{approximation} we have $\lim_{\epsilon\downarrow 0} \sup_{\gamma\in\Gamma_\epsilon}(N_1(\gamma,\epsilon)+N_2(\gamma,\epsilon))=0$ and hence
		\begin{align*}
			\norm{\lim_{\epsilon\downarrow 0} \sum_{\gamma\in I_\epsilon}\eta^\gamma_\epsilon(x)^{2}\Bigl(R_{k}^{a^\gamma_\epsilon(\gamma)}-R_{k}^{a^\gamma_\epsilon(x)}\Bigr)}_\infty&\leq C\lim_{\epsilon\downarrow 0}\sum_{\gamma\in I_\epsilon} \eta^\gamma_\epsilon(x)^{2}(N_1(\gamma,\epsilon)+N_2(\gamma, \epsilon))\\
			&\leq C\lim_{\epsilon\downarrow 0}\Bigl(\sum_{\gamma\in I_\epsilon} \eta^\gamma_\epsilon(x)^{2}\Bigr)\left(\sup_{\gamma\in I_\epsilon}(N_1(\gamma,\epsilon)+N_2(\gamma, \epsilon)\right)\\
			&\leq C \lim_{\epsilon\downarrow 0} \sup_{\gamma\in\Gamma_\epsilon} \bigl(N_1(\gamma,\epsilon)+N_2(\gamma,\epsilon)\bigr)\\
			&=0
		\end{align*}
		Therefore,
		$$\sym(M_{\psi}\quasiriesz{a}{k,\epsilon}M_{\psi})(x)=\psi(x)^{2}\Bigl(R^{a(x)}_k+\lim_{\epsilon\downarrow 0}\sum_{\gamma\in I_\epsilon}\eta^\gamma_\epsilon(x)^{2}\bigl(R_{k}^{a^\gamma_\epsilon(\gamma)}-R_{k}^{a^\gamma_\epsilon(x)}\bigr)\Bigr)=\psi(x)^{2} \quasiriesz{a(x)}{k},$$
		as desired.
	\end{proof}

	\subsection{Proof of Theorem \ref{symbolcalcu}}\label{Proof of Theorem3.1}

	\begin{proof}
		Let $A^{a}_{k,\epsilon}$ and $\quasiriesz{a}{k,\epsilon}$ be given as in Proposition \ref{limitzero}.
		Write
		\begin{align*}
			M_{\psi}\quasiriesz{a}{k}M_{\psi}=M_{\psi}(\quasiriesz{a}{k}-A^{a}_{k,\epsilon})M_{\psi}
			+M_{\psi}(A^{a}_{k,\epsilon}-\quasiriesz{a}{k,\epsilon})M_{\psi}
			+M_{\psi}\quasiriesz{a}{k,\epsilon}M_{\psi}.
		\end{align*}
		Applying Proposition \ref{compactAk}, the first factor on the right hand side is a compact operator and thus belongs to $\Pi$. 
		Applying Proposition \ref{limitzero}, the second factor tends to zero in the uniform norm. 
		It follows that $M_{\psi}\quasiriesz{a}{k}M_{\psi}$ lies in the closure of $\Pi$ with respect to the uniform norm. So $M_{\psi}\quasiriesz{a}{k}M_{\psi}$ belongs to $\Pi$.
		Following Proposition \ref{limitpsiriesz}, the third factor belongs to $\Pi$. 
		
		Using Proposition \ref{compactAk} and Theorem \ref{symbol}, one has
		\begin{align*}
			\sym(M_{\psi}(\quasiriesz{a}{k}-A^{a}_{k,\epsilon})M_{\psi})=0.
		\end{align*}
		Since $\sym:\Pi\rightarrow\A_{1}\otimes_{\rm min}\A_{2}$ is a $*$-homomorphism, it follows from Proposition \ref{limitzero} that, in the uniform norm,
		\begin{align*}
			\lim_{\epsilon\downarrow 0}\sym(M_{\psi}(A^{a}_{k,\epsilon}-\quasiriesz{a}{k,\epsilon})M_{\psi})=0.
		\end{align*}  
		Therefore, by Proposition \ref{limitpsiriesz},
		\begin{align*}
			\sym(M_{\psi}\quasiriesz{a}{k}M_{\psi})(x)=\lim_{\epsilon\downarrow 0}\sym(M_{\psi}\quasiriesz{a}{k,\epsilon}M_{\psi})(x)
			=\psi(x)^{2} \quasiriesz{a(x)}{k}.
		\end{align*}
		This completes the proof.
	\end{proof}

	\section{Equivariance of principal symbol}\label{Equivariance-of-PS}
	In this section we investigate the equivariance of principal symbol mapping $\sym$ under $\GG$-diffeomorphism.
	\begin{definition}
		Let $\Omega,\Omega'\subset\GG$ be open sets and $\Phi:\Omega\to\Omega'$ a diffeomorphism. We obtain an operator $U_{\Phi}:\D'(\Omega')\to\D'(\Omega)$ given by
		\begin{align*}
			U_{\Phi}f=\Jdet{\Phi}^{1/2}V_\Phi(f)
		\end{align*}
		for all $f\in C^\infty_c(\Omega)$, where $V_\Phi$ is as in Theorem \ref{Vdef}, and where $\Jdet{\Phi}:\Omega\to\RR$ is the Jacobian determinant
		$$\Jdet{\Phi}(x)=\lvert\det(J_\Phi(x))\rvert$$
		where $J_\Phi$ is the Jacobian matrix of $\Phi$.
	\end{definition}
	
	When $\Omega=\Omega'$, the restriction of $U_{\Phi}$ to $L_2(\Omega)$, also denoted by $U_\Phi$, is a unitary operator on $L_{2}(\Omega)$.

	\begin{definition}
		Let $\Omega,\Omega'\subset\GG$ and let $\Phi:\Omega\to\Omega'$ be a mapping. Then $\Phi$ is said to be \textit{affine}, if $\Phi$ is the restriction to $\Omega$ of a mapping $\GG\to\GG$ of the form $y\mapsto A(xy)$ for some $x\in\GG$ and some $A\in\at$.
	\end{definition}
	
	\begin{definition}
		An operator $T\in\Pi$ is said to be compactly supported in a set $\Omega\subset\GG$ if there is some $\psi\in\cci$ such that $\supp(\psi)\subset \Omega$ and
		$$M_\psi T=T=TM_\psi.$$
		If $T$ is compactly supported in $\Omega$ for some $\Omega\subset\GG$, we say that $T$ is \textit{compactly supported}.
	\end{definition}
	The following theorem is proved in subsection \ref{Proof of diffeomorphism-section}. 
	\begin{theorem}\label{diffeomorphism}
		Let $T\in\Pi$ be compactly supported in a compact set $K\subset\GG$. 
		Suppose that $\Phi:\GG\to\GG$ is an diffeomorphism such that
		\begin{enumerate}[\rm(1)]
			\item $\Phi$ is a $\GG$-diffeomorphism in some neighborhood of $K$;
			
			\item $\Phi$ is affine outside of some ball.
		\end{enumerate}
		Then
		\begin{align*}
			U_{\Phi}\rp T U_{\Phi}\in\Pi.
		\end{align*}
		Furthermore,
		\begin{align*}
			\sym(U_{\Phi}\rp T U_{\Phi})=\pi_{H^{\Phi}}(\sym(T))\circ\Phi\rp.
		\end{align*}
	\end{theorem}

	Next we investigate the behaviours of the algebra and the principal symbol mapping (locally) under change of coordinates. We need the following notations.
	\begin{definition}
		Let $\H$ be a Hilbert space and let $p\in\B(\H)$ be a projection.
		\begin{enumerate}[\rm(1)]
			\item If $T\in\B(\H)$ with $pTp=T$, then define $\Rest_{p}(T)\in\B(p\H)$ by setting $\Rest_{p}(T)=T|_{p\H}$. 
			\item If $T\in\B(p\H)$, then define $\Ext_{p}(T)\in\B(\H)$ by setting $\Ext_{p}(T)=T\circ p$.
		\end{enumerate}
	\end{definition}

	\begin{notation}\label{Ext-def on G}
		Let $\Omega\subset\GG$ be a subset. 
		\begin{enumerate}[\rm(1)]
			\item If $T\in\B(L_{2}(\GG))$ satisfies $T=M_{\chi_{\Omega}}T M_{\chi_{\Omega}}$, then $\Rest_{\Omega}(T)\in\B(L_{2}(\Omega))$ is a shorthand for $\Rest_{M_{\chi_{\Omega}}}(T)$.
			\item If $T\in\B(L_{2}(\Omega))$, then $\Ext_{\Omega}(T)\in\B(L_{2}(\GG))$ is a shorthand for $\Ext_{M_{\chi_{\Omega}}}(T)$.
		\end{enumerate}
	\end{notation}

	The following theorem is our second result in this section, which is proved in subsection \ref{Local equinvariance-section}. 
	\begin{theorem}\label{Ext-Rest invariant}
		Let $\Omega,\Omega'\subset\GG$ be two open sets and let $\Phi:\Omega\to\Omega'$ be an $\GG$-diffeomorphism. If $T\in\Pi$ is compactly supported in $\Omega$ and $\det(J_{\Phi})\neq0$ on the support of $T$, then $$\Ext_{\Omega'}\left(U_{\Phi}\rp \Rest_{\Omega}(T) U_{\Phi}\right)\in\Pi.$$
		Moreover, $$\sym\left(\Ext_{\Omega'}\left(U_{\Phi}\rp \Rest_{\Omega}(T) U_{\Phi}\right)\right)=\pi_{H^{\Phi}}(\sym(T))\circ\Phi\rp.$$
	\end{theorem}
	Actually, this theorem can be viewed as a corollary of Theorem \ref{diffeomorphism} but it plays an important role in the construction of the principal symbol mapping on $\GG$-filtered manifolds.

	\subsection{$\GG$-Diffeomorphisms}\label{vphivsec}
	In this subsection we prove Theorem \ref{vphiv}. Recall Definition \ref{stratifiedfiltration}, in particular the definition of the grading, $\gra n x$ for $n\in\N$, $x\in\GG$, of $T_x\GG$.
	For the proof of Theorem \ref{vphiv}, it is easy to see (1) $\Longrightarrow$ (2) $\Longrightarrow$ (3) $\Longrightarrow$ (4) $\Longrightarrow$ (1). Next, let us prove the equivalence of \eqref{David's G-diffeomorphism}$\Longleftrightarrow$\eqref{Dima's G-diffeomorphism} and \eqref{David's G-diffeomorphism}$\Longleftrightarrow$\eqref{Fulin's G-diffeomorphism} in Theorem \ref{vphiv}.
	
	The following simple linear algebra lemma is not known by the authors to have previously appeared in the literature.
	\begin{lemma}\label{blockidentity}
		Let $E$ be a vector space, $d=\dim E$, let $(V_i)_{i=1}^d$ and $(W_i)_{i=1}^d$ be bases for $E$, and let $(\omega_{ij})_{i,j=1}^d$ be the change of basis matrix; i.e., for all $1\leq i\leq d$,
		$$W_i=\sum_{j=1}^d \omega_{ij}V_j$$
		Let $1\leq n\leq d$, write $F:=\Span\{W_i:1\leq i\leq n\}$ and suppose that $\omega_{ij}=\delta_{ij}$ for all $1\leq i,j\leq n$. Let $\displaystyle V=\sum_i a_iV_i\in E$ and set $\displaystyle W=\sum_i a_iW_i$. Then $V\in F$ if and only if $V=W$.
	\end{lemma}
	\begin{proof}
		One direction is trivial; since $W\in F$, if $V=W$, then $V\in F$. Suppose that $V\in F$. Then there are scalars $b_1,\hdots, b_n$ such that $\displaystyle V=\sum_i b_iW_i$. Therefore,
		\begin{align*}
			0&=V-V=\sum_{j=1}^d a_jV_j-\sum_{i=1}^n b_iW_j=\sum_{j=1}^d a_jV_j-\sum_{i=1}^n b_i\biggl(\sum_{j=1}^d \omega_{ij}V_j\biggr)=\sum_{j=1}^d \biggl(a_j-\sum_{i=1}^n  b_i\omega_{ij}\biggr)V_j
		\end{align*}
		and hence $\displaystyle{a_j=\sum_{j=1}^n  b_i\omega_{ij}}$ for all $1\leq j\leq d$. Since $\omega_{ij}=\delta_{ij}$ for all $1\leq i,j\leq n$, it follows that $a_i=b_i$ for all $1\leq i\leq n$, and therefore that $V=W$.
	\end{proof}

	\begin{proof}[Proof of \eqref{David's G-diffeomorphism}$\Longleftrightarrow$\eqref{Dima's G-diffeomorphism} in Theorem \ref{vphiv}:] 
		
		Let $x\in\Omega$ and $1\leq k\leq n_1$. By Proposition \ref{pffacts} (\ref{pfvectorcoord}),
		$$(\Phi_*X_k)_{\Phi(x)}=\sum_{j=1}^d (X_k\Phi_j)(x)\frac{\partial}{\partial x_j}\bigg|_{\Phi(x)}$$
		By inspection of the coefficient functions in the change of basis in Remark \ref{basis-rep remark}, one sees that the bases $(\partial/\partial x_i|_{\Phi(x)})_{i=1}^d$ and $(Z_j|_{\Phi(x)})_{j=1}^d$ for $T_{\Phi(x)}\GG$ satisfy the conditions of Lemma \ref{blockidentity} with $n=n_1$, in which case $F=\gra 1 {\Phi(x)}$, and therefore $(\Phi_*)_x(X_k|_x)\in \gra 1 {\Phi(x)}$ if and only if $\displaystyle(\Phi_*X_k)|_{\Phi(x)}=\smash{\sum_{j=1}^{n_1}} (X_k\Phi_j)(x)X_j|_{\Phi(x)}$. Since $(\Phi_*)_x$, it follows that
		$$\Phi_*(\gra 1 x)\subset \gra 1 {\Phi(x)}\,\,\,\forall x\in\Omega\quad\text{if and only if}\quad\Phi_*X_k=\sum_{j=1}^{n_1} ((X_k\Phi_j)\circ\Phi^{-1})X_j\,\,\,\forall\,k\in\{1,\,\dotsc, n_1\},$$
		which is exactly the desired result.
	\end{proof}
	
	\begin{definition}
		Let $\xi_x:T_x\GG\to\RR^d$, for all $x\in\GG$, be the map taking a tangent vector at $x$ to its vector of coordinates relative to $\partial/\partial x_j|_x$. That is, if $V\in T_x\GG$ and
		$\displaystyle V=\sum_{i=1}^d a_i\partial/\partial x_i|_x,$
		then
		$$\xi_x(V)=(a_1,\,\dotsc,a_d).$$
	\end{definition}

	\begin{lemma}\label{itscoordinates}
		For all $x\in\GG$, we have $F(x)=\xi_x(\g_1^x)$, where $F(x)$ is as defined in Theorem \ref{vphiv}.
	\end{lemma}
	\begin{proof}
		Let $x\in\GG$, $2\leq h\leq \iota$, $1\leq s\leq n_h$. Let $y\in F(x)$. The orthogonal complement $E_{h,s}(x)^\perp$ is the set of those $y\in\RR^d$ satisfying
		\begin{equation}\label{ehsorth}
			y_{s}^{(h)}-\sum_{j=1}^{n_{1}}y_{j}^{(1)}a_{j,s}^{(1,h)}(x)=0.
		\end{equation}
		Therefore, it has a basis given by all those $\smash{e^{(i)}_j}$ for which $i\geq 2$ and $(i,j)\neq (h,s)$, and the vectors
		$$e^{(1)}_j+a^{(1,h)}_{j,s}(x)e^{(h)}_s,$$
		for $1\leq j\leq n_1$. Therefore $E_{h,s}(x)^\perp$ also contains $\smash{\alpha^{(1,h)}_{j,k}(x)e^{(i)}_j}$ whenever $i\geq 2$ and $(i,j)\neq (h,s)$ and therefore it also contains the sum
		\begin{equation*}
			h_j(x):=\bigl(e^{(1)}_j+a^{(1,h)}_{j,s}(x)e^{(h)}_s\bigr)+\biggl(\mathop{\sum\sum}\limits_{(l,j)\neq (h,s)}\alpha^{(1,l)}_{j,k}(x)e^{(l)}_{k}\biggr)=e^{(1)}_j+\sum_{l=2}^\iota\sum_{k=1}^{n_l}\alpha^{(1,l)}_{j,k}(x)e^{(l)}_{k}.
		\end{equation*}
		The vectors $h_j(x)$ are linearly independent, since, for every $1\leq k\leq n_1$, there is exactly one vector in the set $\{h_j(x):1\leq l\leq n_1\}$ which has a nonzero $\smash{e^{(1)}_j}$ coordinate. By (\ref{ehsorth}), the $\smash{e^{(h)}_s}$ coordinate of any element of $E_{h,s}(x)$ is determined by the $\smash{e^{(1)}_j}$ coordinates, $1\leq j\leq n_1$, and hence the coordinates of any element of $F(x)$ is determined by the first $n_1$ coordinates. More explicitly, the linear map $F(x)\to \RR^{n_1}$, given by projecting onto the first $n_1$ entries, is injective. Therefore, $F(x)$ has dimension at most $n_1$. Now, the $h_j(x)$ satisfy (\ref{ehsorth}) for every $1\leq j\leq n_1$ regardless of the value of $h$ and $s$. and hence $h_j(x)\in F(x)$ for every $1\leq j\leq n_1$. By linear independence and dimension counting, the $h_j(x)$ form a basis for $F(x)$. The proof is now completed by \eqref{pushforward} 
		and by using \eqref{basis-rep} with the form $h_j(x)=\xi_x(X_j|_x)$.
	\end{proof}

	\begin{proof}[Proof of \eqref{David's G-diffeomorphism}$\Longleftrightarrow$\eqref{Fulin's G-diffeomorphism} in Theorem \ref{vphiv}:]
		
		Let $x\in\Omega$ and $1\leq k\leq n_1$. Since $\xi_x$ is a linear isomorphism, by applying \eqref{pfptcoord} with the form $(\Phi_{*})_{x}(V)=\xi_{\Phi(x)}\rp\left(J_{\Phi}(x)\xi_{x}(V)\right)$ for $V\in T_{x}\GG$ and Lemma \ref{itscoordinates}, we have a chain of equivalences
		\begin{align*}
			(\Phi_{*})_{x} (X_k|_x)\in \g_1^{\Phi(x)}&\iff\xi_{\Phi(x)}((\Phi_{*})_{x} (X_k|_x))\in \xi_{\Phi(x)}\big(\g_1^{\Phi(x)}\big)\\
			&\iff J_\Phi(x)\xi_x(X_k|_x)\in  \xi_{\Phi(x)}\big(\g_1^{\Phi(x)}\big)\\
			&\iff J_\Phi(x)\xi_x(X_k|_x)\in  F(\Phi(x))
		\end{align*}
		Therefore, by Remark \ref{spanner}, Lemma \ref{itscoordinates}, and linearity, we have
		$$(\Phi_{*})_{x} (\g_1^x)\subset \g_1^{\Phi(x)}\iff J_\Phi(x)\xi_x(F(x))\subset F(\Phi(x)),$$
		and since $x$ was arbitrary, the proof is completed.
	\end{proof}

	\subsection{Compactness \expandafter{\romannumeral2}}
	\
	\newline

	The main task is to show that principal symbol mapping is equivariant under $\GG$-diffeomorphism.
	
	\begin{lemma}\label{multiplication}
		Let $\Omega, \Omega'\subset\GG$ be open, let $\Phi:\Omega\to\Omega'$ be a diffeomorphism and let $f\in \lp{\infty}$. Then, in $\B(L_2(\Omega))$ we have
		\begin{equation}\label{UconjMult}
			U_{\Phi}\rp M_{f}U_{\Phi}=V_{\Phi}\rp M_{f}V_{\Phi}=\Phi_*(M_f)=M_{f\circ\Phi\rp}
		\end{equation}
	\end{lemma}
	\begin{proof}
		By Proposition \ref{pffacts}, one has
		$$U_{\Phi}\rp M_{f}U_{\Phi}=(M_{\Jdet{\Phi}}^{1/2}V_{\Phi})^{-1} M_{f}(M_{\Jdet{\Phi}}^{1/2}V_{\Phi})=V_{\Phi}\rp M_{\Jdet{\Phi}^{-1/2}f\Jdet{\Phi}^{1/2}}V_{\Phi}=V_{\Phi}\rp M_{f}V_{\Phi},$$
		which establishes the first equality in (\ref{UconjMult}). The second equality is trivial, and the third follows from Proposition \ref{pffacts} (\ref{pfmultiplier}).
	\end{proof}
	
	\begin{lemma}\label{derivation on distributions}
		Let $\Omega\subset\GG$ be open and let $X\in\Gamma(T\Omega)$ be a vector field on $\Omega$. Then, for all $f,g\in \D'(\Omega)$, if $f$ is smooth or $g$ is smooth, the product distributions $fg$, $(Xf)g$ and $f(Xg)$ are defined and one has
		$$X(fg)=(Xf)g+f(Xg).$$
	\end{lemma}
	\begin{proof}
		Suppose that $g\in C^\infty(\Omega)$. The operator $f\mapsto X(fg)-(Xf)g-f(Xg)$ on $\D'(\Omega)$ is continuous and vanishes on the dense set $C_c^\infty(\Omega)$, and therefore it vanishes on $\D'(\Omega)$, hence $X(fg)=(Xf)g+f(Xg)$. The proof in the case that $f$ is smooth is analogous.
	\end{proof}
	
	\begin{lemma}\label{vector field multiplier commutator}
		Let $\Omega\subset\GG$ be open, let $X\in\Gamma(T\Omega)$ be a vector field on $\Omega$ and let $f\in L_\infty(\Omega)$. Then we have the equality of operators on $W^{1,2}(\Omega)$,
		$$[X, M_f]=M_{Xf}$$
	\end{lemma}
	\begin{proof}
		If $g\in W^{1,2}(\Omega)\cap C^\infty_c(\Omega)$ is smooth, then by Lemma \ref{derivation on distributions},
		\begin{align*}
			[X, M_f]g=(XM_f)g-(M_fX)g=X(fg)-f(Xg)=(Xf)g+f(Xg)-f(Xg)=M_{Xf}g.
		\end{align*}
		Therefore the operator $f\mapsto [X, M_f]-M_{Xf}g$ vanishes on a dense subset of $W^{1,2}(\Omega)$, and hence vanishes everywhere by continuity, therefore $[X, M_f]=M_{Xf}$.
	\end{proof}
	
	\begin{lemma}\label{UXUplus}
		Let $\Omega, \Omega'\subset\GG$ be open and let $\Phi:\Omega\to\Omega'$ be a diffeomorphism. For all $1\leq k\leq n_1$, define
		\begin{align*}
			a^{\Phi}_{k}:=(\Jdet{\Phi}^{-1/2} X_{k}\Jdet{\Phi}^{1/2})\circ\Phi\rp,
			\quad1\leq k\leq n_{1}.
		\end{align*}
		Then we have the following equality of operators on $W^{1,2}(\Omega)$.
		\begin{align}\label{UXU-HPhi}
			U_{\Phi}^{-1}X_{k}U_{\Phi}=\Phi_*(X_{k})+a_{k}^{\Phi}
		\end{align}
		
	\end{lemma}
	\begin{proof}
		By Lemma \ref{vector field multiplier commutator}, Lemma \ref{multiplication} and  Lemma \ref{pffacts}, 
		\begin{align*}
			U_{\Phi}^{-1}X_{k}U_{\Phi}&=V_\Phi^{-1}M_{\Jdet\Phi}^{-1/2}X_k M_{\Jdet\Phi}^{1/2}V_\Phi=\Phi_*\Bigl(M_{\Jdet\Phi}^{-1/2}X_k M_{\Jdet\Phi}^{1/2}\Bigr)\\
			&=\Phi_*\Bigl(X_k+M_{\Jdet\Phi}^{-1/2}\bigl[X_k, M_{\Jdet\Phi}^{1/2}\bigr]\Bigr)\\
			&=\Phi_*(X_k)+\Phi_*\Bigl(M_{\Jdet\Phi}^{-1/2}M_{X_k\Jdet{\Phi}^{1/2}}\Bigr)\\
			&=\Phi_*(X_k)+a_{k}^{\Phi}
		\end{align*}
		as desired.
	\end{proof}

	\begin{lemma}\label{deltadecomposition}
		Let $\psi\in \cci$ and $\Phi:\GG\rightarrow \GG$ be a diffeomorphism. 
		Let $a:\GG\rightarrow\at$ be a smooth mapping. Suppose that
		\begin{enumerate}
			\item $\Phi$ is a $\GG$-diffeomorphism in the neighborhood of $\supp(\psi\circ\Phi)$;
			
			\item $\Phi$ is affine outside of some ball;
			
			\item $a$ and $a\rp$ are bounded;
			
			\item $a=H^{\Phi}\circ\Phi\rp$ on $\supp(\psi)$.
		\end{enumerate}
		If $A\in\at$, then there are $T_{1},T_{2}\in\B(\lp{2})$ and $\psi_{2}\in\cci$ such that
		\begin{align*}
			\Big((1+U_{\Phi}\rp{L_{A}} U_{\Phi})^{-1/2}-(1+ L_{aA})^{-1/2}\Big)M_{\psi}
			&=\int_{0}^{\infty}(1+t^{2}+ L_{aA})^{-1}T_{1}\frac{ L_{aA}^{1/2}}{1+t^{2}+ L_{aA}}M_{\psi_{2}}\,dt\\
			&+\int_{0}^{\infty}U_{\Phi}\rp(1+t^{2}+{L_{A}})^{-1}U_{\Phi}T_{2}\frac{ L_{aA}^{1/2}}{1+t^{2}+ L_{aA}}M_{\psi_{2}}\,dt.
		\end{align*}
	\end{lemma}
	\begin{proof}
		Denote ${\rm LHS}=\Big((1+U_{\Phi}\rp{L_{A}} U_{\Phi})^{-1/2}-(1+ L_{aA})^{-1/2}\Big)M_{\psi}$. 
		Applying Lemma \ref{deltadecomo} to $1+U_{\Phi}\rp{L_{A}} U_{\Phi}$ and $1+ L_{aA}$, it follows that
		\begin{align*}
			{\rm LHS}
			=\frac{2}{\pi}\int_{0}^{\infty}\Big((1+t^{2}+U_{\Phi}\rp{L_{A}} U_{\Phi})^{-1}-(1+t^{2}+ L_{aA})^{-1}\Big)M_{\psi}\,dt
			=:\frac{2}{\pi}\int_{0}^{\infty}I_{t}\,dt.
		\end{align*}
		
		Write $\psi=\psi_{1}\psi_{2}$ with $\psi_{1},\psi_{2}\in \cci$ and $\supp(\psi_{1})=\supp(\psi)$. 
		Then
		\begin{align}\label{equalityo}
			\nonumber I_{t}&=(1+t^{2}+U_{\Phi}\rp{L_{A}} U_{\Phi})^{-1}( L_{aA}-U_{\Phi}\rp{L_{A}} U_{\Phi})M_{\psi_{1}}(1+t^{2}+ L_{aA})^{-1}M_{\psi_{2}}+\\
			\nonumber&+(1+t^{2}+U_{\Phi}\rp{L_{A}} U_{\Phi})^{-1}( L_{aA}-U_{\Phi}\rp{L_{A}} U_{\Phi})[(1+t^{2}+ L_{aA})^{-1},M_{\psi_{1}}]M_{\psi_{2}}\\
			&=:I_{1,t}-I_{2,t}.
		\end{align}
		
		By \eqref{commutator-inverse},
		\begin{align*}
			&I_{2,t}=\\
			&(1+t^{2}+U_{\Phi}\rp{L_{A}} U_{\Phi})^{-1}( L_{aA}-U_{\Phi}\rp{L_{A}} U_{\Phi})(1+t^{2}+ L_{aA})^{-1}[ L_{aA},M_{\psi_{1}}](1+t^{2}+ L_{aA})^{-1}M_{\psi_{2}}
		\end{align*}
		where we use equality $[1+t^{2}+ L_{aA},M_{\psi_{1}}]=[ L_{aA},M_{\psi_{1}}]$.
		Note that
		$$ L_{aA}-U_{\Phi}\rp{L_{A}} U_{\Phi}=(1+t^{2}+ L_{aA})-(1+t^{2}+U_{\Phi}\rp{L_{A}} U_{\Phi}).$$
		It yields that
		\begin{align}\label{equalityt}
			\nonumber I_{2,t}&=(1+t^{2}+U_{\Phi}\rp{L_{A}} U_{\Phi})^{-1}[ L_{aA},M_{\psi_{1}}](1+t^{2}+ L_{aA})^{-1}M_{\psi_{2}}+\\
			&-(1+t^{2}+ L_{aA})^{-1}[ L_{aA},M_{\psi_{1}}](1+t^{2}+ L_{aA})^{-1}M_{\psi_{2}}.
		\end{align}
		
		Let
		\begin{align}\label{T1-bd}
			T_{1}=\frac{2}{\pi}[ L_{aA},M_{\psi_{1}}]\Delta^{-1/2}\cdot 
			\Delta^{1/2} L_{aA}^{-1/2}.
		\end{align}
		Similar to \eqref{ML-commutator}, it follows that $[ L_{aA},M_{\psi_{1}}]$ is a differential operator with bounded compactly supported coefficients. By Cwikel estimate in \cite{MSZ-Cwikel} and Lemma \ref{laplacinbd}, one concludes that $T_{1}$ is a bounded operator.
		
		Let
		\begin{align}\label{T2-bd}
			T_{2}=\frac{2}{\pi}( L_{aA}-U_{\Phi}\rp{L_{A}} U_{\Phi})M_{\psi_{1}} L_{aA}^{-1/2}-\frac{2}{\pi}[ L_{aA},M_{\psi_{1}}] L_{aA}^{-1/2}.
		\end{align}
		By $a=H^{\Phi}\circ\Phi\rp$ on $\supp(\psi_{1})$ and Lemma \ref{UXUplus}, one has
		\begin{align*}
			U_{\Phi}\rp{L_{A}} U_{\Phi}M_{\psi_{1}}
			&=-\sum_{k,j,i=1}^{n_{1}}A_{jk}A_{ik}
			\left(X_{j}^{a}X_{i}^{a}+X_{j}^{a}M_{a_{i}^{\Phi}}+M_{a_{j}^{\Phi}}X_{i}^{a}+M_{a_{j}^{\Phi}a_{i}^{\Phi}}\right)M_{\psi_{1}}\\
			&=L_{aA}M_{\psi_{1}}-\sum_{k,j,i=1}^{n_{1}}A_{jk}A_{ik}
			\left(X_{j}^{a}M_{a_{i}^{\Phi}}+M_{a_{j}^{\Phi}}X_{i}^{a}+M_{a_{j}^{\Phi}a_{i}^{\Phi}}\right)M_{\psi_{1}}.
		\end{align*}
		One concludes that $(L_{aA}-U_{\Phi}\rp{L_{A}} U_{\Phi})M_{\psi_{1}}$ is a differential operator with bounded compactly supported coefficients.
		Combining with assertion on $T_{1}$, it yields that $T_{2}$ is also a bounded operator.
		
		Combining \eqref{equalityo}, \eqref{equalityt}, \eqref{T1-bd} and \eqref{T2-bd}, we have
		\begin{align*}
			{\rm LHS}=\frac{2}{\pi}\int_{0}^{\infty}(I_{1,t}-I_{2,t})\,dt
			&=\int_{0}^{\infty}(1+t^{2}+ L_{aA})^{-1}T_{1}\frac{ L_{aA}^{1/2}}{1+t^{2}+ L_{aA}}M_{\psi_{2}}\,dt+\\
			&+\int_{0}^{\infty}(1+t^{2}+U_{\Phi}\rp{L_{A}} U_{\Phi})^{-1}T_{2}\frac{ L_{aA}^{1/2}}{1+t^{2}+ L_{aA}}M_{\psi_{2}}\,dt.
		\end{align*}
		This yields the desired result.
	\end{proof}

	\begin{lemma}\label{term4}
		Let $\psi$, $\Phi$ and $a$ be given as in Lemma \ref{deltadecomposition}. If $A\in\at$, then
		$$M_{\psi}U_{\Phi}\rp X_{k}^{A}U_{\Phi}\Big((1+U_{\Phi}\rp L_{A} U_{\Phi})^{-1/2}-(1+{L_{aA}})^{-1/2}\Big)M_{\psi}\in\K(\lp{2}).$$
	\end{lemma}
	\begin{proof}
		Denote ${\rm LHS}$ the above term on the left hand side.
		By  Lemma \ref{deltadecomposition}, one has
		\begin{align*}
			{\rm LHS}
			&=\int_{0}^{\infty}M_{\psi}U_{\Phi}\rp X_{k}^{A}U_{\Phi}(1+t^{2}+{L_{aA}})^{-1}T_{1}\frac{L_{aA}^{1/2}}{1+t^{2}+{L_{aA}}}M_{\psi_{2}}\,dt+\\
			&+\int_{0}^{\infty}M_{\psi}U_{\Phi}\rp X_{k}^{A}(1+t^{2}+{L_{A}})^{-1}U_{\Phi}T_{2}\frac{L_{aA}^{1/2}}{1+t^{2}+{L_{aA}}}M_{\psi_{2}}\,dt\\
			&=:{\rm Tail1}+{\rm Tail2},
		\end{align*}
		where $T_{1},T_{2}$ and $\psi_{2}$ are given in  Lemma \ref{deltadecomposition}.
		
		Since $a=H^{\Phi}\circ\Phi\rp$ on $\supp(\psi)$, it follows from Lemma \ref{UXUplus} that
		\begin{align*}
			M_{\psi}U_{\Phi}\rp X_{k}^{A}U_{\Phi}\Delta^{-1/2}
			=\sum_{j=1}^{n_{1}}\sum_{i=1}^{n_{1}}A_{jk}M_{\psi a_{ij}}X_{i}\Delta^{-1/2}
			+\sum_{j=1}^{n_{1}}A_{jk}M_{\psi a_{j}^{\Phi}}\Delta^{-1/2}.
		\end{align*}
		Therefore $M_{\psi}U_{\Phi}\rp X_{k}^{A}U_{\Phi} L_{A}^{-1/2}$ is a bounded operator by Cwikel estimate in \cite{MSZ-Cwikel}.
		By the spectral theorem,
		\begin{align}\label{lapineq}
			\Norm{\frac{ L_{aA}^{1/2}}{1+t^{2}+ L_{aA}}}_{\infty}\leq\sup_{s>0}\frac{s}{1+t^{2}+s^{2}}
			\lesssim\frac{1}{(1+t^{2})^{1/2}}.
		\end{align}
		Since $T_{1}\in\B(\lp{2})$,  it follows from Lemma \ref{laplacinbd} that
		\begin{align*}
			&\Norm{M_{\psi}U_{\Phi}\rp X_{k}^{A}U_{\Phi}(1+t^{2}+ L_{aA})^{-1}T_{1}\frac{ L_{aA}^{1/2}}{1+t^{2}+ L_{aA}}M_{\psi_{2}}}_{\infty}\\
			&\leq\Norm{M_{\psi}U_{\Phi}\rp X_{k}^{A}U_{\Phi} \Delta^{-1/2}}_{\infty}\Norm{ \Delta^{1/2} L_{aA}^{-1/2}}_{\infty}
			\Norm{\frac{ L_{aA}^{1/2}}{1+t^{2}+ L_{aA}}}_{\infty}^{2}\Norm{T_{1}}_{\infty}\Norm{\psi_{2}}_{\infty}\\
			&\lesssim\frac{1}{1+t^{2}}.
		\end{align*}
		This implies that the integrand ${\rm Tail1}$ is convergent in the uniform norm.

		Applying Lemma \ref{pffacts} and the spectral theorem, one obtains
		\begin{align*}
			\Norm{X_{k}^{A}(1+t^{2}+ L_{A})^{-1}}_{\infty}
			\leq\Norm{X_{k}^{A} L_{A}^{-1/2}}_{\infty}\Norm{\frac{ L_{A}^{1/2}}{1+t^{2}+ L_{A}}}_{\infty}
			\lesssim\frac{1}{(1+t^{2})^{1/2}}.
		\end{align*}
		Moreover, by $T_{2}\in\B(\lp{2})$ and \eqref{lapineq},
		\begin{align*}
			&\Norm{M_{\psi}U_{\Phi}\rp X_{k}^{A}(1+t^{2}+ L_{A})^{-1}U_{\Phi}T_{2}\frac{ L_{aA}^{1/2}}{1+t^{2}+ L_{aA}}M_{\psi_{2}}}_{\infty}\\
			&\leq\Norm{\psi}_{\infty}\Norm{\psi_{2}}_{\infty}\Norm{X_{k}^{A}(1+t^{2}+ L_{A})^{-1}}_{\infty}
			\Norm{T_{2}}_{\infty}\Norm{\frac{ L_{aA}^{1/2}}{1+t^{2}+ L_{aA}}}_{\infty}\\
			&\lesssim\frac{1}{1+t^{2}}.
		\end{align*}
		This implies that integrand ${\rm Tail2}$ is convergent in the uniform norm.
		
		Write
		$$\frac{ L_{aA}^{1/2}}{1+t^{2}+ L_{aA}}M_{\psi_{2}}=\frac{ L_{aA}}{1+t^{2}+ L_{aA}}\cdot L_{aA}^{-1/2} \Delta^{1/2}\cdot \Delta^{-1/2}M_{\psi_{2}}.$$
		By Cwikel estimate in \cite{MSZ-Cwikel}, it is a compact operator. 
		Since integrands ${\rm Tail1}$ and ${\rm Tail2}$ are convergent in the uniform norm, it follows that they are compact.
		The proof is completed.
	\end{proof}

	\begin{proposition}\label{rieszcompact}
		Let $\psi$, $\Phi$ and $a$ be given as in Lemma \ref{deltadecomposition}. If $A\in\at$, then
		\begin{align*}
			M_{\psi}(U_{\Phi}\rp \quasiriesz{A}{k}U_{\Phi}-\quasiriesz{aA}{k})M_{\psi}\in\K(\lp{2}).
		\end{align*}
	\end{proposition}
	
	\begin{proof}
		Write
		\begin{align*}
			M_{\psi}(U_{\Phi}\rp \quasiriesz{A}{k}U_{\Phi}-\quasiriesz{aA}{k})M_{\psi}={\rm Term1}+{\rm Term2}+{\rm Term3}+{\rm Term4},
		\end{align*}
		where
		$${\rm Term1}=M_{\psi}U_{\Phi}\rp X_{k}^{A}\Big(L_{A}^{-1/2}-(1+L_{A})^{-1/2}\Big)U_{\Phi}M_{\psi},$$
		$${\rm Term2}=M_{\psi}(U_{\Phi}\rp X_{k}^{A}U_{\Phi}- X_{k}^{aA})(1+L_{aA})^{-1/2}M_{\psi},$$
		$${\rm Term3}=M_{\psi} X_{k}^{aA}\Big((1+L_{aA})^{-1/2}-L_{aA}^{-1/2}\Big)M_{\psi}$$
		and
		$${\rm Term4}=M_{\psi}U_{\Phi}\rp X_{k}^{A}U_{\Phi}\Big((1+U_{\Phi}\rp L_{A} U_{\Phi})^{-1/2}-(1+ L_{aA})^{-1/2}\Big)M_{\psi}.$$
		It is left to show that these four terms are separately compact.
		
		Write
		\begin{align*}
			{\rm Term1}
			&=M_{\psi}U_{\Phi}\rp X_{k}^{A}L_{A}^{-1/2}\cdot\Big(1-\frac{L_{A}^{1/2}}{(1+L_{A})^{1/2}}\Big)U_{\Phi}M_{\psi}\\
			&=M_{\psi}U_{\Phi}\rp X_{k}^{A}L_{A}^{-1/2}\cdot V_{A}\rp\Big(1-\frac{\Delta^{1/2}}{(1+\Delta)^{1/2}}\Big)M_{\psi\circ\Phi}V_{A}U_{\Phi}.
		\end{align*}
		Let $\phi(t)=1-\frac{t^{1/2}}{(1+t)^{1/2}}$. 
		One has $\phi\in L_{Q,\infty}(\RR_{+},t^{Q/2-1}dt)$. 
		It follows from Cwikel estimate \cite[Corollary 4.5]{MSZ-Cwikel} that $\phi(\Delta)M_{\psi\circ\Phi}=\Big(1-\frac{\Delta^{1/2}}{(1+\Delta)^{1/2}}\Big)M_{\psi\circ\Phi}$ is a compact operator. 
		So ${\rm Term1}$ is a compact operator.
		
		By the assumption on $a$ and Lemma \ref{UXUplus}, one has
		\begin{align*}
			M_{\psi}U_{\Phi}\rp X_{k}^{A}U_{\Phi}
			&=\sum_{j=1}^{n_{1}}A_{jk}U_{\Phi}\rp X_{j}U_{\Phi}\\
			&=\sum_{j=1}^{n_{1}}\sum_{i=1}^{n_{1}}A_{jk}M_{\psi a_{ij}}X_{i}
			+\sum_{j=1}^{n_{1}}A_{jk}M_{\psi a_{j}^{\Phi}}\\
			&=M_{\psi}X_{k}^{aA}+\sum_{j=1}^{n_{1}}A_{jk}M_{\psi a_{j}^{\Phi}}.
		\end{align*}
		That is
		\begin{align*}
			{\rm Term2}=
			\sum_{j=1}^{n_{1}}A_{jk}M_{\psi a_{j}^{\Phi}}\Delta^{-1/2}
			\cdot\Delta^{1/2} L_{aA}^{-1/2}
			\cdot L_{aA}^{1/2}(1+ L_{aA})^{-1/2}M_{\psi}.
		\end{align*}
		Applying Lemma \ref{laplacinbd} and functional calculus, the second and third factors are bounded.
		Using Cwikel estimate in \cite{MSZ-Cwikel}, the first factor is a compact operator and therefore so is ${\rm Term2}$.
		
		Write
		\begin{align*}
			{\rm Term3}=-M_{\psi}X_{k}^{aA}L_{aA}^{-1/2}\cdot
			\frac{ L_{aA}^{1/2}}{(1+ L_{aA})^{1/2}((1+ L_{aA})^{1/2}+ L_{aA}^{1/2})}
			\cdot L_{aA}^{-1/2}\Delta^{1/2}\cdot\Delta^{-1/2}M_{\psi},
		\end{align*}
		where we use the equality $1-\frac{L_{aA}^{1/2}}{(1+ L_{aA})^{1/2}}=\frac{1}{(1+ L_{aA})^{1/2}((1+ L_{aA})^{1/2}+ L_{aA}^{1/2})}$.
		Applying Lemma \ref{quasiriszbd} and  Lemma \ref{laplacinbd}, the first and third factors are bounded. 
		The boundedness of second factor follows from functional calculus.
		By Cwikel estimate in \cite{MSZ-Cwikel}, the last factor is compact and so is ${\rm Term3}$.
		
		The compactness of ${\rm Term4}$ follows from Lemma \ref{term4}.
		This completes the proof. 
	\end{proof}

	\subsection{Equivariance of quasi-Riesz transforms}
	
	\begin{lemma}\label{productcontain}
		Let $\{f_{k}\}_{k=1}^{m}\subset\A_{1}$ and $\{g_{k}\}_{k=1}^{m}\subset\A_{2}$. Then
		\begin{align*}
			\prod_{k=1}^{m}\po(f_{k})\pt(g_{k})\in\po(\prod_{k=1}^{m}f_{k})\pt(\prod_{k=1}^{m}g_{k})+\K(\lp{2}).
		\end{align*}
	\end{lemma}
	\begin{proof}
		We prove the assertion by induction on $m$. For $m=1$, there is nothing to show. 
		
		Let us prove the assertion for $m=2$. 
		We have
		\begin{align*}
			\po(f_{1})\pt(g_{1})\po(f_{2})\pt(g_{2})
			&=[\pt(g_{1}),\po(f_{1}f_{2})]\pt(g_{2})+\\
			&+[\po(f_{1}),\po(g_{1})]\po(f_{2})\pt(g_{2})+\po(f_{1}f_{2})\pt(g_{1}g_{2}).
		\end{align*}
		By Lemma \ref{compactness},
		\begin{align*}
			[\po(f_{1}),\po(g_{1})],[\pt(g_{1}),\po(f_{1}f_{2})]\in\K(\lp{2}).
		\end{align*}
		Therefore,
		\begin{align*}
			\po(f_{1})\pt(g_{1})\po(f_{2})\pt(g_{2})\in \po(f_{1}f_{2})\pt(g_{1}g_{2})+\K(\lp{2}).
		\end{align*}
		This shows the assertion for $m=2$.
		
		It remains to show the step of induction. Suppose the assertion holds for $m\geq2$ and let us show it holds for $m+1$.
		Clearly,
		\begin{align*}
			\prod_{k=1}^{m+1}\po(f_{k})\pt(g_{k})=\po(f_{1})\pt(g_{1})\cdot \prod_{k=2}^{m+1}\po(f_{k})\pt(g_{k}).
		\end{align*}
		Using the inductive assumption, one has
		\begin{align*}
			\prod_{k=1}^{m+1}\po(f_{k})\pt(g_{k})
			\in \po(f_{1})\pt(g_{1})\po(\prod_{k=2}^{m+1}f_{k})\pt(\prod_{k=2}^{m+1}g_{k})+\K(\lp{2}).
		\end{align*}
		Applying the assertion for $m=2$, one obtains
		\begin{align*}
			\po(f_{1})\pt(g_{1})\po(\prod_{k=2}^{m+1}f_{k})\pt(\prod_{k=2}^{m+1}g_{k})
			\in \po(\prod_{k=1}^{m+1}f_{k})\pt(\prod_{k=1}^{m+1}g_{k})+\K(\lp{2}).
		\end{align*}
		Thus,
		\begin{align*}
			\prod_{k=1}^{m+1}\po(f_{k})\pt(g_{k})\in \po(\prod_{k=1}^{m+1}f_{k})\pt(\prod_{k=1}^{m+1}g_{k})+\K(\lp{2}).
		\end{align*}
		This establishes the step of induction. It completes the proof.
	\end{proof}
	

	\begin{lemma}\label{rieszpi}
		Let $\phi\in\cci$ and $\Phi:\GG\rightarrow \GG$ be a diffeomorphism. 
		Suppose that
		\begin{enumerate}
			\item $\Phi$ is a $\GG$-diffeomorphism in some neighborhood of $\supp(\phi)$;
			
			\item $\Phi$ is affine outside of some ball.
			
			
		\end{enumerate}
		If $A\in\at$, then
		\begin{align*}
			U_{\Phi}\rp M_{\phi}\quasiriesz{A}{k}M_{\phi}U_{\Phi}\in\Pi.
		\end{align*}
		Moreover,
		\begin{align*}
			\sym(U_{\Phi}\rp M_{\phi}\quasiriesz{A}{k}M_{\phi}U_{\Phi})=\pi_{H^{\Phi}}(\sym(M_{\phi}\quasiriesz{A}{k}M_{\phi}))\circ\Phi\rp.
		\end{align*}
		Similarly,
		\begin{align*}
			U_{\Phi}\rp M_{\phi}(\quasiriesz{A}{k})^{*}M_{\phi}U_{\Phi}\in\Pi.
		\end{align*}
		and
		\begin{align*}
			\sym(U_{\Phi}\rp M_{\phi}(\quasiriesz{A}{k})^{*}M_{\phi}U_{\Phi})
			=\pi_{H^{\Phi}}(\sym(M_{\phi}(\quasiriesz{A}{k})^{*}M_{\phi}))\circ\Phi\rp.
		\end{align*}
	\end{lemma}
	\begin{proof}
		Noting that $\Phi$ is a $\GG$-diffeomorphism in some neighborhood of $\supp(\phi)$,  Theorem \ref{admissibility} gives the existence of $H^{\Phi}$. Let $a\in\cca$ be a smooth function constant outside of a ball and satisfying that $a=H^{\Phi}$ near the support of $\phi$.
		By  Lemma \ref{multiplication}, 
		\begin{align}\label{UMRMU-decompose}
			\nonumber U_{\Phi}\rp M_{\phi}\quasiriesz{A}{k}M_{\phi}U_{\Phi}
			&=M_{\phi\circ\Phi\rp} U_{\Phi}\rp \quasiriesz{A}{k}U_{\Phi} M_{\phi\circ\Phi\rp}\\
			&=M_{\phi\circ\Phi\rp} (U_{\Phi}\rp \quasiriesz{A}{k}U_{\Phi}-\quasiriesz{(a\circ\Phi\rp)A}{k}) M_{\phi\circ\Phi\rp}
			+M_{\phi\circ\Phi\rp} \quasiriesz{(a\circ\Phi\rp)A}{k} M_{\phi\circ\Phi\rp}.
		\end{align}
		Applying  Proposition \ref{rieszcompact} to $\psi=\phi\circ\Phi\rp$, it follows that the first factor in \eqref{UMRMU-decompose} is a compact operator.

		Applying Theorem \ref{symbolcalcu}, one has
		\begin{align*}
			M_{\phi\circ\Phi\rp} \quasiriesz{(a\circ\Phi\rp)A}{k} M_{\phi\circ\Phi\rp}\in\Pi
		\end{align*}
		and
		\begin{align*}
			\sym(M_{\phi\circ\Phi\rp} \quasiriesz{(a\circ\Phi\rp)A}{k} M_{\phi\circ\Phi\rp})(x)
			=(\phi\circ\Phi\rp(x))^{2}\quasiriesz{(a\circ\Phi\rp(x))A}{k},\quad x\in \GG.
		\end{align*}
		It follows from Theorem \ref{symbol} that the compact operators are contained in $\Pi$ and their symbols are zero. 
		Therefore, $$U_{\Phi}\rp M_{\phi}\quasiriesz{A}{k}M_{\phi}U_{\Phi}\in\Pi.$$ 
		Moreover,
		\begin{align*}
			\sym(U_{\Phi}\rp M_{\phi}\quasiriesz{A}{k}M_{\phi}U_{\Phi})(x)
			=\sym(M_{\phi\circ\Phi\rp} \quasiriesz{(a\circ\Phi\rp)A}{k} M_{\phi\circ\Phi\rp})(x)
			=(\phi\circ\Phi\rp(x))^{2}\quasiriesz{(a\circ\Phi\rp(x))A}{k}.
		\end{align*}
		Observing that $V_{A}V_{B}=V_{BA}$, and therefore, for $y=\Phi\rp(x)$,
		\begin{align*}
			\quasiriesz{a(y)A}{k}
			=V_{a(y)A}\rp R_{k} V_{a(y)A}
			=V_{a(y)}\rp V_{A}\rp R_{k} V_{A}V_{a(y)},
		\end{align*}
		one obtains
		\begin{align*}
			\phi(y)^{2}\quasiriesz{a(y)A}{k}
			=\phi(y)^{2}V_{H^{\Phi}(y)}\rp \quasiriesz{A}{k}V_{H^{\Phi}(y)}
			=V_{H^{\Phi}(y)}\rp\sym(M_{\phi}\quasiriesz{A}{k}M_{\phi})(y)V_{H^{\Phi}(y)}.
		\end{align*}
		This yields that
		\begin{align*}
			\sym(U_{\Phi}\rp M_{\phi}\quasiriesz{A}{k}M_{\phi}U_{\Phi})(x)
			=\pi_{H^{\Phi}}(\sym(M_{\phi}\quasiriesz{A}{k}M_{\phi}))(\Phi\rp(x)).
		\end{align*}
		The proof is completed.
		
	\end{proof}

	\begin{lemma}\label{ufginclusion}
		Let $f\in C_{c}(\GG)$ and $\Phi:\GG\rightarrow \GG$ be a diffeomorphism.
		Suppose that 
		\begin{enumerate}
			\item $\Phi$ is a $\GG$-diffeomorphism in some neighborhood of $\supp(f)$;
			
			\item $\Phi$ is affine outside of some ball. 
		\end{enumerate}
		If $g\in\A_{2}$, then
		\begin{align*}
			U_{\Phi}\rp \po(f)\pt(g)U_{\Phi}\in\Pi.
		\end{align*}
		Moreover,
		\begin{align*}
			\sym(U_{\Phi}\rp \po(f)\pt(g)U_{\Phi})=\pi_{H^{\Phi}}(\sym(\po(f)\pt(g)))\circ\Phi\rp.
		\end{align*}
	\end{lemma}
	\begin{proof}
		Let ${\rm Alg}(\set{R_{k}}_{k=1}^{n_{1}})$ be the $*$-algebra generated by $\{\quasiriesz{A}{k}:A\in\at\}_{k=1}^{n_{1}}$.
		There is a sequence $\set{g_{k}}_{k\geq1}\subset {\rm Alg}(\set{R_{k}}_{k=1}^{n_{1}})$ such that $g_{k}\rightarrow g$ in the uniform norm.
		It is valuable to consider the case when $g$ is a monomial.
		
		Suppose first that $g$ is a monomial.
		Let $\displaystyle g=\prod_{j=1}^{m}r_{j}$ for $r_{j}\in\{\quasiriesz{A}{k},(\quasiriesz{A}{k})^{*}:A\in\at\}_{k=1}^{n_{1}}$.
		As $f\in C_{c}(\GG)$, select a positive $\phi\in\cci$ such that $f\phi=f$. Clearly,
		\begin{align}\label{fgrepo}
			\po(f)\pt(g)=\po(f)\po(\phi^{2m})\pt(g).
		\end{align}
		Let $\psi_{j}=\phi^{2}$ for $1\leq j\leq m$. By Lemma \ref{productcontain},
		\begin{align*}
			\po(\phi^{2m})\pt(g)=\po(\prod_{j=1}^{m}\psi_{j})\pt(\prod_{j=1}^{m}r_{j})\in\prod_{j=1}^{m}\po(\psi_{j})\pt(r_{j})+\K(\lp{2}).
		\end{align*}
		Note that
		\begin{align*}
			\po(\psi_{j})\pt(r_{j})=\po(\phi)\pt(r_{j})\po(\phi)+\po(\phi)[\po(\phi),r_{j}].
		\end{align*}
		Since $\phi\in\cci$ and $r_{j}\in\A_{2}$, it follows from Lemma \ref{compactness} that
		\begin{align*}
			\po(\psi_{j})\pt(r_{j})\in \po(\phi)\pt(r_{j})\po(\phi)+\K(\lp{2}).
		\end{align*}
		Thus,
		\begin{align*}
			\po(\phi^{2m})\pt(g)\in\prod_{j=1}^{m}\po(\phi)\pt(r_{j})\po(\phi)+\K(\lp{2}).
		\end{align*}
		This implies that there is a $\kappa\in\K(\lp{2})$ satisfying that
		\begin{align}\label{fgrept}
			\po(\phi^{2m})\pt(g)=\prod_{j=1}^{m}\po(\phi)\pt(r_{j})\po(\phi)+\kappa.
		\end{align}
		Combining equalities (\ref{fgrepo}) and (\ref{fgrept}), one obtains
		\begin{align}\label{monomial-sum}
			U_{\Phi}\rp\po(f)\pt(g)U_{\Phi}=
			U_{\Phi}\rp\po(f)U_{\Phi}\prod_{j=1}^{m}U_{\Phi}\rp\po(\phi)\pt(r_{j})\po(\phi)U_{\Phi}
			+U_{\Phi}\rp\po(f)\kappa U_{\Phi}.
		\end{align}
		Clearly,
		\begin{align*}
			U_{\Phi}\rp\po(f)\kappa U_{\Phi}\in\K(\lp{2}).
		\end{align*}
		By Theorem \ref{symbol},
		\begin{align}\label{monomial-zero}
			\sym(U_{\Phi}\rp\po(f)\kappa U_{\Phi})=0.
		\end{align}
		Observe that $\Pi$ is a Banach algebra.
		By  Lemma \ref{multiplication} and  Lemma \ref{rieszpi},
		\begin{align*}
			U_{\Phi}\rp\po(f)U_{\Phi}\prod_{j=1}^{m}U_{\Phi}\rp\po(\phi)\pt(r_{j})\po(\phi)U_{\Phi}\in\Pi.
		\end{align*}
		This concludes that $$U_{\Phi}\rp \po(f)\pt(g)U_{\Phi}\in\Pi.$$
		Combining Lemma \ref{multiplication} and  Theorem \ref{symbol}, equality in \eqref{transfer} implies that
		\begin{align}\label{monomial-multiplication}
			\sym(U_{\Phi}\rp\po(f)U_{\Phi})=(f\circ \Phi\rp) \otimes\one=\pi_{H^{\Phi}}(f\otimes\one)\circ\Phi\rp.
		\end{align}
		By Lemma \ref{rieszpi}, 
		\begin{align}\label{monomilal-riesz}
			\sym(U_{\Phi}\rp\po(\phi)\pt(r_{j})\po(\phi)U_{\Phi})
			&=\pi_{H^{\Phi}}(\sym(\po(\phi)\pt(r_{j})\po(\phi)))\circ\Phi\rp\\
			\nonumber&=\pi_{H^{\Phi}}(\phi^{2}\otimes r_{j})\circ\Phi\rp.
		\end{align}
		Since $\sym$ is a $*$-homomorphism, it follows from \eqref{monomial-sum}, \eqref{monomial-zero}, \eqref{monomial-multiplication} and \eqref{monomilal-riesz} that
		\begin{align*}
			\sym(U_{\Phi}\rp \po(f)\pt(g)U_{\Phi})
			=\pi_{H^{\Phi}}(f\phi^{2m}\otimes \prod_{j=1}^{m}r_{j})\circ\Phi\rp
			=\pi_{H^{\Phi}}(\sym(\po(f)\pt(g)))\circ\Phi\rp.
		\end{align*}
		This proves the assertion for the case when $g$ is a monomial. 
		By linearity, the same assertion holds if $g\in{\rm Alg}(\set{R_{k}}_{k=1}^{n_{1}})$.
		
		For the general case, let $g\in\A_{2}$ and take a sequence $\set{g_{k}}_{k\geq1}\subset{\rm Alg}(\set{R_{k}}_{k=1}^{n_{1}})$ such that $g_{k}\rightarrow g$ in the uniform norm.
		By the assertion in the monomial case, we have
		$$U_{\Phi}\rp \po(f)\pt(g_{k})U_{\Phi}\in\Pi.$$
		Note that
		\begin{align*}
			\Norm{U_{\Phi}\rp \po(f)\pt(g)U_{\Phi}-U_{\Phi}\rp \po(f)\pt(g_{k})U_{\Phi}}_{\infty}
			\leq\Norm{f}_{\infty}\Norm{g-g_{k}}_{\infty}\rightarrow0.
		\end{align*}
		Since $\Pi$ is closed in the uniform norm, it follows that
		\begin{align*}
			U_{\Phi}\rp \po(f)\pt(g)U_{\Phi}\in\Pi.
		\end{align*}
		Moreover,
		\begin{align*}
			\sym(U_{\Phi}\rp \po(f)\pt(g)U_{\Phi})&=\lim_{k\rightarrow\infty}\sym(U_{\Phi}\rp \po(f)\pt(g_{k})U_{\Phi})\\
			&=\lim_{k\rightarrow\infty}\pi_{H^{\Phi}}(\sym(\po(f)\pt(g_{k})))\circ\Phi\rp\\
			&=\pi_{H^{\Phi}}(\sym(\po(f)\pt(g)))\circ\Phi\rp.
		\end{align*}
		This completes the proof.
	\end{proof}
	

	\subsection{Proof of Theorem \ref{diffeomorphism}}\label{Proof of diffeomorphism-section}

	\begin{proof}
		Since $\Pi$ is a $C^{*}$-algebra, there exist a sequence $\set{T_{k}}_{k\in\NN}$ in the $*$-algebra generated by $\po(\A_{1})$ and $\pt(\A_{2})$ such that $T_{k}\rightarrow T$ in the uniform norm.
		Write
		\begin{align*}
			T_{k}=\sum_{\gamma=1}^{\gamma_{k}}\prod_{ l=1}^{ l_{k}}\po(f_{k,\gamma, l})\pt(g_{k,\gamma, l}),
			\quad f_{k,\gamma, l}\in\A_{1},\quad g_{k,\gamma, l}\in\A_{2}.
		\end{align*}
		By  Lemma \ref{productcontain}, 
		\begin{align}\label{elementexist}
			T_{k}\in\sum_{\gamma=1}^{\gamma_{k}}\po(\prod_{ l=1}^{ l_{k}}f_{k,\gamma, l})\pt(\prod_{ l=1}^{ l_{k}}g_{k,\gamma, l})+\K(\lp{2}).
		\end{align}
		Denote for brevity,
		\begin{align*}
			f_{k,\gamma}=\prod_{ l=1}^{ l_{k}}f_{k,\gamma, l}\quad{\rm and}\quad g_{k,\gamma}=\prod_{ l=1}^{ l_{k}}g_{k,\gamma, l}.
		\end{align*}
		By formula \eqref{elementexist}, there exists an $\omega_{k}\in\K(\lp{2})$ such that
		\begin{align*}
			T_{k}=\sum_{\gamma=1}^{\gamma_{k}}\po(f_{k,\gamma})\pt(g_{k,\gamma})+\omega_{k}.
		\end{align*}
		Since $x$ is compactly supported in $K$, it follows that there is $\psi\in\cci$ with $\supp(\psi)=K$ satisfying $T=\po(\psi)T$.
		Therefore,
		\begin{align*}
			\Norm{T-\po(\psi)T_{k}}_{\infty}\leq\Norm{\psi}_{\infty}\Norm{T-T_{k}}_{\infty}\rightarrow0.
		\end{align*}
		Hence, it is valuable to take $\po(\psi)T_{k}$ into consideration. 
		Note that
		\begin{align}\label{psixk}
			\po(\psi)T_{k}=\sum_{\gamma=1}^{\gamma_{k}}\po(\psi f_{k,\gamma})\pt(g_{k,\gamma})+\po(\psi)\omega_{k}.
		\end{align}
		Combining Theorem \ref{symbol}, one has
		\begin{align*}
			U_{\Phi}\rp\po(\psi)\omega_{k}U_{\Phi}\in\K(\lp{2})
			\quad\mbox{and}\quad
			\sym(U_{\Phi}\rp\po(\psi)\omega_{k}U_{\Phi})=0.
		\end{align*}
		Applying  Lemma \ref{ufginclusion} to $\psi f_{k,\gamma}\in C_{c}(G)$, it follows that
		$$U_{\Phi}\rp\po(\psi f_{k,\gamma})\pt(g_{k,\gamma})U_{\Phi}\in\Pi$$
		and that
		\begin{align*}
			\sym(U_{\Phi}\rp\po(\psi f_{k,\gamma})\pt(g_{k,\gamma})U_{\Phi})=\pi_{H^{\Phi}}(\po(\psi f_{k,\gamma})\pt(g_{k,\gamma}))\circ\Phi\rp.
		\end{align*}
		Thus, by \eqref{psixk}, we have $U_{\Phi}\rp\po(\psi)T_{k}U_{\Phi}\in\Pi$ and
		\begin{align*}
			\sym(U_{\Phi}\rp\po(\psi)T_{k}U_{\Phi})
			&=\sum_{\gamma=1}^{\gamma_{k}}\sym(U_{\Phi}\rp\po(\psi f_{k,\gamma})\pt(g_{k,\gamma})U_{\Phi})\\
			&=\sum_{\gamma=1}^{\gamma_{k}}\pi_{H^{\Phi}}(\sym(\po(\psi f_{k,\gamma})\pt(g_{k,\gamma})))\circ\Phi\rp\\
			&=\pi_{H^{\Phi}}(\sym(\po(\psi)T_{k}))\circ\Phi\rp.
		\end{align*}
		
		Since $\Pi$ is a $C^{*}$-algebra and $U_{\Phi}\rp \po(\psi)T_{k}U_{\Phi}\rightarrow U_{\Phi}\rp T U_{\Phi}$ in the uniform norm, $U_{\Phi}\rp T U_{\Phi}\in\Pi$ and 
		\begin{align*}
			\sym(U_{\Phi}\rp\po(\psi)T_{k}U_{\Phi})\rightarrow\sym(U_{\Phi}\rp T U_{\Phi})
		\end{align*}
		in the uniform norm. In other words, 
		\begin{align*}
			\pi_{H^{\Phi}}(\sym(\po(\psi)T_{k}))\circ\Phi\rp\rightarrow\pi_{H^{\Phi}}(\sym(T))\circ\Phi\rp
		\end{align*}
		in the uniform norm. Observing that $\sym(\po(\psi)T_{k})\rightarrow\sym(T)$ in the uniform norm, it gives that
		\begin{align*}
			\sym(U_{\Phi}\rp T U_{\Phi})=\pi_{H^{\Phi}}(\sym(T))\circ\Phi\rp.
		\end{align*}
		This completes the proof.
	\end{proof}

	\subsection{Local equinvariance}\label{Local equinvariance-section}

	\begin{lemma}\label{local equal}
		Let $\Omega\subset\GG$ be an open set and let $\Phi:\Omega\to\GG$ be a smooth mapping. If $\xi\in\Omega$ satisfies $\det(J_{\Phi}(\xi))\neq0$, then there is a diffeomorphism $\Phi_{\xi}:\GG\to\GG$ such that
		\begin{enumerate}[\rm(1)]
			\item $\Phi_{\xi}=\Phi$ on some ball $B(\xi,r_{1,\xi})$ for some $r_{1,\xi}>0$;
			\item $\Phi_{\xi}$ is affine outside of a ball $B(\xi,r_{2,\xi})$ for some $r_{2,\xi}<\infty$.
		\end{enumerate}
	\end{lemma}
	\begin{proof}
		For $\xi\in\GG$, let $\psi\in C_{c}^{\infty}(\GG)$ be such that $\psi=1$ on the ball $B(\xi,1)$. Set $$\Psi_{r}(\zeta)=\zeta\rp\cdot\left(\psi(\delta_{r\rp}(\zeta\rp\xi))(\zeta\cdot\Phi(\zeta))\right).$$ 
		Obviously, $\det(J_{\Psi_{r}})$ never vanishes for sufficiently small $r>0$ and $\det(J_{\Psi_{r}})\to1$ as $r\to0$. It follows from \cite[Corollary 4.3]{Palais1959} that $\Psi_{r}:\RR^{d}\to\RR^{d}$ is a diffeomorphism for sufficiently small $r>0$. Choose any such $r$ associated to $\xi$ and denote it by $r_{\xi}$. Set $\Phi_{\xi}=\Psi_{r_{\xi}}$. This diffeomorphism satisfies the required properties.
	\end{proof}

	\begin{lemma}\label{local conjugate equal}
		Let $\Omega,\Omega'\subset\GG$ and let $\Phi:\Omega\to\Omega'$ be a diffeomorphism. Let $B\subset\Omega$ be a ball and let $\Phi_{0}:\GG\to\GG$ be a diffeomorphism such that $\Phi_{0}=\Phi$ on $B$. If $T\in\B(L_{2}(\GG))$ is supported in $B$, then
		$$\Ext_{\Omega'}\left(U_{\Phi}\rp\Rest_{\Omega}(T) U_{\Phi}\right)=U_{\Phi_{0}}\rp T U_{\Phi_{0}}.$$
	\end{lemma}
	\begin{proof}
		Using the same arguments in \cite[Lemma 6.3]{DAO3}, we obtain the desired result.
		
	\end{proof}

	Let us start the proof of Theorem \ref{Ext-Rest invariant}.
	\begin{proof}[Proof of Theorem \ref{Ext-Rest invariant}]
		Let $K$ be the compact set of $\Omega$ such that $T\in\Pi$ is supported in $K$. 
		For $x\in K$, since $\det(J_{\Phi}(x))\neq0$, there is a diffeomorphism $\Phi_{x}:\GG\to\GG$ satisfying arguments in Lemma \ref{local equal}. Without loss of generality, we assume that $r_{1,x}$ is so small such that $B(x,r_{1,x})\subset\Omega$. Therefore, $\Phi_{x}$ is an $\GG$-diffeomorphism on $B(x,r_{1,x})$.
		And the open balls $\{B(x,r_{1,x})\}_{K}$ cover $K$. By the compactness of $K$, there is a finite sub-cover covering $K$, saying $\{B(x_{k},r_{1,x_{k}})\}_{k=1}^{N}$.
		Let $\{\phi_{k}\}_{k=1}^{N}$ be such that $\phi_{k}\in C_{c}^{\infty}(B(x_{k},r_{1,x_{k}}))$ and $$\sum_{k=1}^{N}\phi_{k}^{2}=1\quad\mbox{on}\quad K.$$
		Set $$T_{0}=\sum_{k=1}^{N}M_{\phi_{k}}[M_{\phi_{k}},T],\quad T_{k}=M_{\phi_{k}}T M_{\phi_{k}},\quad 1\leq k\leq N.$$
		Therefore, we have
		\begin{align}\label{local decomposition}
			T=\sum_{k=1}^{N}M_{\phi_{k}^{2}}T=\sum_{k=0}^{N}T_{k}.
		\end{align}
		
		When $1\leq k\leq N$, $T_{k}$ is compactly supported in $B(x_{k},r_{1,x_{k}})$. By Lemma \ref{local conjugate equal}, we have $$\Ext_{\Omega'}\left(U_{\Phi}\rp \Rest_{\Omega}(T_{k}) U_{\Phi}\right)=U_{\Phi_{x_{k}}}\rp T_{k} U_{\Phi_{x_{k}}}.$$
		Since $\Phi_{x_{k}}$ is an $\GG$-diffeomorphism on $B(x_{k},r_{1,x_{k}})$ and is affine outside of some ball, it follows from Theorem \ref{admissibility} and Theorem \ref{diffeomorphism} that $$U_{\Phi_{x_{k}}}\rp T_{k} U_{\Phi_{x_{k}}}\in\Pi$$ and $$\sym\left(U_{\Phi_{x_{k}}}\rp T_{k} U_{\Phi_{x_{k}}}\right)=\pi_{H^{\Phi_{x_{k}}}}(\sym(T_{k}))\circ\Phi_{x_{k}}\rp.$$
		Observing that $\sym(T_{k})$ is supported in $B(x_{k},r_{1,x_{k}})$ and $\Phi_{x_{k}}=\Phi$ on $B(x_{k},r_{1,x_{k}})$, we obtain $$\pi_{H^{\Phi_{x_{k}}}}(\sym(T_{k}))\circ\Phi_{x_{k}}\rp=\pi_{H^{\Phi}}(\sym(T_{k}))\circ\Phi\rp.$$
		Thus, 
		\begin{align}\label{local ext Tk}
			\Ext_{\Omega'}\left(U_{\Phi}\rp \Rest_{\Omega}(T_{k}) U_{\Phi}\right)\in\Pi
		\end{align} 
		and
		\begin{align}\label{local sym Tk}
			\sym\left(\Ext_{\Omega'}\left(U_{\Phi}\rp \Rest_{\Omega}(T_{k}) U_{\Phi}\right)\right)=\pi_{H^{\Phi}}(\sym(T_{k}))\circ\Phi\rp.
		\end{align}

		When $k=0$, it follows from Lemma \ref{compactness} that $T_{0}$ is a compact operator and therefore $\sym(T_{0})=0$. Clearly, $T_{0}$ is compactly supported in $\Omega$. If $A\in\B(L_{2}(\Omega'))$ is compact, then $\Ext_{\Omega'}(A)\in\B(L_{2}(\GG))$ is also compact. Thus, 
		\begin{align}\label{local ext T0}
			\Ext_{\Omega'}\left(U_{\Phi}\rp \Rest_{\Omega}(T_{0}) U_{\Phi}\right)\in\K(L_{2}(\GG))\subset\Pi
		\end{align} 
		and
		\begin{align}\label{local sym T0}
			\sym\left(\Ext_{\Omega'}\left(U_{\Phi}\rp \Rest_{\Omega}(T_{0}) U_{\Phi}\right)\right)=0.
		\end{align}
		
		Combining \eqref{local decomposition}, \eqref{local ext Tk} and \eqref{local ext T0}, we have $$\Ext_{\Omega'}\left(U_{\Phi}\rp \Rest_{\Omega}(T) U_{\Phi}\right)\in\Pi.$$
		Moreover, by \eqref{local decomposition}, \eqref{local sym Tk} and \eqref{local sym T0},
		\begin{align*}
			\sym\left(\Ext_{\Omega'}\left(U_{\Phi}\rp \Rest_{\Omega}(T) U_{\Phi}\right)\right)
			&=\sum_{k=0}^{N}\sym\left(\Ext_{\Omega'}\left(U_{\Phi}\rp \Rest_{\Omega}(T_{k}) U_{\Phi}\right)\right)\\
			&=\sum_{k=1}^{N}\pi_{H^{\Phi}}(\sym(T_{k}))\circ\Phi\rp\\
			&=\pi_{H^{\Phi}}(\sym(T))\circ\Phi\rp-\pi_{H^{\Phi}}(\sym(T_{0}))\circ\Phi\rp\\
			&=\pi_{H^{\Phi}}(\sym(T))\circ\Phi\rp.
		\end{align*}
		This completes the proof.
	\end{proof}

	\section{Globalisation theorem}\label{Globalisation theorem}
	In this section, we investigate the globalisation theorem for the existence of a $*$-homomorphism which is originated from \cite{DAO3} in the setting of compact manifold.
	\subsection{Local algebras}
	\begin{definition}
		Let $M$ be a topological space. An open cover $\{U_i\}_{i\in\II}$ (atlas $\{(\U_i, h_i)\}_{i\in\II}$, partition of unity $\{\phi_i\}_{i\in\II}$) of $M$ is said to be \textit{locally finite} if, for each $x\in M$, there are only finitely many $i\in\II$ such that $x\in U_i$, ($x\in\U_i$, $\phi_i\neq 0$).
	\end{definition}
	
	It is well known that every manifold admits a locally finite atlas, and it is clear that any partition of unity subordinate to a locally finite open cover is itself locally finite.
	
	\begin{definition}\label{localalg}
		Let $M$ be a manifold and $\{(\U_i, h_i)\}_{i\in I}$ be a countable locally finite atlas where each $\U_i$ has compact closure. Let $\Borel$ be the Borel $\sigma$-algebra on $M$ and let $\nu$ be a countably additive positive measure on $\Borel$. We say that $\{\A_i\}_{i\in\II}$ are local algebras if, for every $i\in\II$,
		\begin{enumerate}[\rm(1)]
			\item $\A_i$ is a $*$-subalgebra of $\B(L_2(M, \nu))$;
			\item \itemcase{e}very element of $\A_i$ is compactly supported in $\U_i$;
			\item \itemcase{f}or every $j\in\II$, if $T\in\A_i$ is compactly supported in $\U_i\cap\U_j$, then $T\in \A_j$;
			\item \itemcase{i}f $T\in\K(L_2(M, \nu))$ and $T$ is compactly supported in $\U_i$ then $T\in\A_i$;
			\item \itemcase{i}f $\phi\in C_c(\U_i)$, then $M_\phi\in\A_i$;
			\item \itemcase{i}f $\phi\in C_c(\U_i)$, then the closure of $M_\phi\A_i M_\phi$ in the uniform norm of $\B(L_{2}(M,\nu))$ is contained in $\A_i$;
			\item \itemcase{i}f $T\in\A_i$ and $\phi\in C_c(\U_i)$, then $[T, M_\phi]\in \K(L_2(M, \nu))$.
		\end{enumerate}
	\end{definition}
	
	\begin{definition}\label{A}
		In the setting of Definition \ref{localalg}, we say that $T\in\A$ if
		\begin{enumerate}[\rm(1)]
			\item \itemcase{f}or every $i\in\II$ and for every $\phi\in C_c(\U_i)$, we have $M_\phi T M_\phi\in \A_i$;
			\item \itemcase{f}or every $\phi\in C_{c}(M)$, the commutator $[T, M_\phi]$ is compact on $L_{2}(M,\nu)$, 
			\item \itemcase{t}here exists a sequence $(\psi_n)$ in $C_{c}(M)$ such that $M_{\psi_n}T\to T$ uniformly in $\B(L_{2}(M,\nu))$.
		\end{enumerate}
	\end{definition}
	
	Recall that, for $M$ a topological space, we say $f\in C_0(M)$ if $f\in C(M)$ and for each $\epsilon>0$, there is a compact set $K_{\epsilon}\subset M$ such that $|f(x)|<\epsilon$ for all $x\in M\setminus K_{\epsilon}$.
	
	\begin{definition}\label{localhomo}
		Let $\B$ be a $C^*$-algebra. In the setting of Definition \ref{localalg}, $\{\homo_i\}_{i\in\II}$ are called local homomorphisms if
		\begin{enumerate}[\rm(1)]
			\item \itemcase{F}or every $i\in\II$, $\homo_i:\A_i\to \B$ is a $*$-homomorphism;
			\item \itemcase{f}or every $i,j\in I$, we have $\homo_i=\homo_j$ on $\A_i\cap\A_j$;
			\item $T\in\A_i$ is compact if and only if $\homo_i(T)=0$;
			\item \itemcase{t}here exists a $*$-homomorphism $\Homo:C_0(M)\to \B$ such that
			$$\Homo(\phi)=\homo_i(M_\phi),\qquad\phi\in C_c(\U_i),\,i\in\II.$$
		\end{enumerate}
	\end{definition}

	The following theorem is our globalisation theorem for the existence of a $*$-homomorphism in which we do not allow the compactness on manifold.
	\begin{theorem}\label{globalisation}
		In the setting of Definitions \ref{localalg}, \ref{localhomo} and \ref{A}, we have
		\begin{enumerate}[\rm(1)]
			\item\label{Cstar algebra} $\A$ is a $C^*$-algebra which contains $\A_i$ for every $i\in\II$ and $\K(L_2(M,\nu))$;
			\item\label{star homomorphism} \itemcase{t}here exists a $*$-homomorphism $\homo:\A\to\B$ such that
			\begin{enumerate}
				\item $\homo=\homo_i$ on $\A_i$ for every $i\in\II$;
				\item $\ker(\homo)=\K(L_2(M,\nu))$,
			\end{enumerate}
			\item\label{star homomorphism unique} \itemcase{t}he $*$-homomorphism as in (2) is unique.
		\end{enumerate}
	\end{theorem}
	
	We prove this theorem in parts. Let us introduce three lemmas at first.
	\begin{lemma}\label{convergy1}
		Suppose that $T\in \B(L_2(M,\nu))$, $(\psi_n)$ is a sequence of nonzero functions in $C_{c}(M)$, and that $M_{\psi_n}T\to T$ uniformly. Then, if $(\phi_n)$ is any sequence in $C_{c}(M)$ converging uniformly on compact sets to a bounded function $\phi$, we have $M_{\phi_n}T\to M_\phi T$ uniformly.
	\end{lemma}
	\begin{proof}
		The case $T=0$ is trivial, so assume $\norm{T}>0$. Let $\varepsilon>0$ and let $m\in\NN$ be such that $\norm{M_{\psi_m}T-T}<\frac{\varepsilon}{2\norm{\phi}_\infty}$. Since $K:=\supp(\phi_m)$ is compact, there is $N\in\NN$ such that $n\geq N$ implies $|\phi_n-\phi|<\frac{\varepsilon}{2\norm T\norm{\psi_m}_\infty}$ on $K$. Thus, if $n\geq N$, we have
		\begin{align*}
			\norm{M_{\phi_n}T-M_{\phi}T}&=\norm{M_{\phi_n-\phi}(T-M_{\psi_m}T+M_{\psi_m}T)}\\
			&\leq \norm{M_{\phi_n-\phi}(T-M_{\psi_m}T)}+\norm{M_{(\phi_n-\phi)\psi_m}T}\\
			&\leq \norm{\phi}_\infty\norm{T-M_{\psi_m}T}+\norm{(\phi_n-\phi)\psi_m}_\infty\norm{T}\\
			&<\varepsilon
		\end{align*}
		hence $M_{\phi_n}T\to M_\phi T$ uniformly.
	\end{proof}
	
	\begin{lemma}\label{convergy2}
		Let $M$ be a topological space and let $(\phi_i)_{i\in\NN}$ a locally finite partition of unity on $M$. Then the series $\displaystyle\sum_{i\in\NN} \phi_i$ converges uniformly on compact sets.
	\end{lemma}
	\begin{proof}
		Let $K\subset M$ be compact and let $\varepsilon>0$. For each $x\in K$, let
		$$N_x:=\max\{N\in\NN : \phi_i(x)\neq 0\}$$
		which is finite by the local finiteness condition. Since $\displaystyle\sum_{i=1}^N\phi_i$ is continuous, there is a neighbourhood $U_x$ of $M$ such that $\displaystyle\sum_{i=1}^N\phi_i(y)>1-\varepsilon$ for all $y\in U_x$. Since $K$ is compact, there is a finite sequence $x_1,\dotsc x_N\in K$ such that $\{U_{x_i}\}_{1\leq i\leq n}$ covers $K$. Set $N=\max\{N_{x_i}:1\leq i\leq n\}$. Then clearly $\displaystyle\sum_{i=1}^N\phi_i(x)>1-\varepsilon$ on $K$, hence $\displaystyle\sum_{i=1}^\infty\phi_i$ converges uniformly on $K$.
	\end{proof}
	
	\begin{lemma}\label{representation}
		Let $M$ be a topological space, let $T\in \B(H)$ and suppose that there exists a sequence $(\psi_n)$ of nonzero functions in $C_{c}(M)$, such that $M_{\psi_n}T\to T$ uniformly. Then, if $\II$ is countable and $(\phi_i)_{i\in\II}$ is a locally finite partition of unity with $\phi_{i}\in C_{c}(M)$, we have
		$$T=\sum_{i\in\II} M_{\phi_i}T=\sum_{i\in\II} (M_{\phi_i^{1/2}} TM_{\phi_i^{1/2}}+M_{\phi_i^{1/2}}[T, M_{\phi_i^{1/2}}])$$
		with the series converging uniformly.
	\end{lemma}
	\begin{proof}
		Without loss of generality, $\II=\NN$. By Lemma \ref{convergy2}, the sequence $\psi_n:=\sum_{i=1}^n\phi_i$ converges uniformly on compact sets to the constant function 1, hence
		\begin{align*}
			T=\lim_{n\to\infty} M_{\psi_n}T=\lim_{n\to\infty} \sum_{i=1}^nM_{\phi_i}T=\sum_{i=1}^\infty M_{\phi_i}T
		\end{align*}
		uniformly. The second identity follows from the observation that
		$$M_{\phi_i}T=M_{\phi_i^{1/2}} TM_{\phi_i^{1/2}}+M_{\phi_i^{1/2}}[T, M_{\phi_i^{1/2}}]$$
		for any $i\in\II$.
	\end{proof}

	\subsection{Proof of the globalisation theorem}
	
	Let us now fix notation and hypotheses for the rest of the proof. Let $M$ be a manifold and $\{(\U_i, h_i)\}_{i\in I}$ be a countable locally finite atlas where each $\U_i$ has compact closure. Let $(\phi_i)_{i\in\II}$ be a partition of unity subordinate to $\{\U_i\}_{i\in\II}$. Let $\Borel$ be the Borel $\sigma$-algebra on $M$, let $\nu$ be a countably additive positive measure on $\Borel$ and let $\{\A_i\}_{i\in\II}$ be local algebras. Finally, let $\B$ be a $C^*$-algebra and let $\{\homo_i\}_{i\in\II}$ be local homomorphisms.
	
	\begin{proposition}\label{globalisation-algebra}
		$\A$ is a $C^*$-algebra.
	\end{proposition}
	\begin{proof}
		The only nontrivial parts are closure under products and completeness. Suppose $T,S\in \A$. First, note that, if $\psi\in C_0(M)$,
		\begin{align*}
			[TS, M_\psi]&=TSM_\psi-M_\psi TS\\
			&=TSM_\psi-(T M_\psi S+[T, M_\psi]S)\\
			&=T[S, M_\psi]-[T, M_\psi]S
		\end{align*}
		which is compact by hypothesis, verifying condition (2) of Definition \ref{A} for $TS$. For condition (1) of Definition \ref{A} for $TS$, let $i\in\II$ and suppose $\phi\in C_c(\U_i)$. We may assume that $0\leq\phi$. We have
		\begin{align*}
			M_\phi TS M_\phi &= M_{\phi^{1/2}}M_{\phi^{1/2}}TSM_{\phi^{1/2}}M_{\phi^{1/2}}\\
			&=M_{\phi^{1/2}}(TM_{\phi^{1/2}}+[T,M_{\phi^{1/2}}])(SM_{\phi^{1/2}}+[S,M_{\phi^{1/2}}])M_{\phi^{1/2}}\\
			&=M_{\phi^{1/2}}TM_{\phi^{1/2}}M_{\phi^{1/2}}SM_{\phi^{1/2}}+K
		\end{align*}
		For some $K$. Now $\A_i$ is an algebra and $M_{\phi^{1/2}}TM_{\phi^{1/2}},M_{\phi^{1/2}}SM_{\phi^{1/2}}\in \A_i$ hence $M_{\phi^{1/2}}TM_{\phi^{1/2}}M_{\phi^{1/2}}SM_{\phi^{1/2}}\in\A_i$. Additionally, $K$ can easily be seen to be both compact and compactly supported in $\U_i$, hence lies in $A_i$. It follows that $M_\phi TSM_\phi\in\A_i$, as desired. We now look to condition (3) of Definition \ref{A} for $TS$. Let $(\psi_n)_{n\in\NN}$ be a sequence in $C_{c}(M)$ converging to the constant function 1 uniformly on compact sets. Then, by the continuity of products in $B(H)$,
		\begin{align*}
			TS=(\lim_{n\to\infty} M_{\psi_n} T)S=(\lim_{n\to\infty} M_{\psi_n} )TS
		\end{align*}
		hence $TS$ validates condition (3).
		
		To see that $\A$ is complete, suppose that $T_n\to T$ uniformly. Condition (2) of Definition \ref{A} for $T$ follows by the continuity of products and the closedness of $\K(L_2(M,\nu))$; if $\psi\in C_0(M)$, then
		$$\lim_{n\to\infty}[T_n, M_\psi]=\lim_{n\to\infty}(T_nM_\psi-M_\psi T_n)=TM_\psi-M_\psi T$$
		hence $[T, M_\psi]$ is compact. For (1) of Definition \ref{A} for $T$, let $i\in\II$ and $\phi\in C_c(\U_i)$. Take $\phi_0\in C_c(\U_i)$ such that $\phi\phi_0=\phi$. Then
		$$M_\phi T_n M_\phi=M_{\phi_0}M_\phi T_n M_\phi M_{\phi_0}$$
		now $M_\phi T_n M_\phi\to M_\phi TM_\phi$ in the norm, hence $M_\phi TM_\phi$ lies in the norm closure of $M_{\phi_0} \A_iM_{\phi_0}$. But by hypothesis (Definition \ref{localalg} (6)), this closure is contained in $\A_i$. Thus $M_\phi TM_\phi\in A_i$. Finally, suppose that $(\psi_n)_{n\in\NN}$ is a sequence in $C_{c}(M)$ converging to the constant function 1 uniformly on compact sets. Let $\varepsilon>0$. There is $m\in\NN$ such that $\norm{T-T_m}<\frac{\varepsilon}{3}$. Having chosen $m$, observe that, since $T_m\in A$, there is $N\in\NN$ such that $n\geq N$ implies $\norm{M_{\psi_n}T_m-T_m}<\frac{\varepsilon}{3}$. Thus if $n\geq N$, we has
		\begin{align*}
			\norm{M_{\psi_n} T-T}&=\norm{M_{\psi_n} (T-T_m+T_m)-T}\\
			&\leq \norm{M_{\psi_n}}\norm{T-T_m}+\norm{M_{\psi_n}T_m-T}\\
			&\leq \norm{T-T_m}+\norm{M_{\psi_n}T_m-T_m}+\norm{T_m-T}\\
			&<\varepsilon
		\end{align*}
		thus $M_{\psi_n}T\to T$ uniformly so $T$ satisfies condition (3) of Definition \ref{A}.
	\end{proof}
	
	\begin{proposition}\label{homconst}
		Suppose that $\homo$ satisfies the conditions in Theorem \ref{globalisation} (2). Then, for all $T\in\A$,
		$$\homo(T)=\sum_{i\in\II} \homo_i(M_{\phi_i^{1/2}}TM_{\phi_i^{1/2}}).$$
		In particular, $\homo$ is unique.
	\end{proposition}
	\begin{proof}
		Let $T\in A$. Since $\homo$ is a $*$-homomorphism between $C^*$-algebras - and therefore continuous with respect to the norm, hence applying $\homo$ to the formula given in Lemma \ref{representation}, one obtains
		\begin{align*}
			\homo(T)=\sum_{i\in\II}\left(\homo(M_{\phi_i^{1/2}} TM_{\phi_i^{1/2}})+\homo(M_{\phi_i^{1/2}}[T, M_{\phi_i^{1/2}}])\right).
		\end{align*}
		The result now follows from the observations that $M_{\phi_i^{1/2}} TM_{\phi_i^{1/2}}\in\A_i$ for $i\in\II$ so $\homo(M_{\phi_i^{1/2}} TM_{\phi_i^{1/2}})=\homo_i(M_{\phi_i^{1/2}} TM_{\phi_i^{1/2}})$and that $[T, M_{\phi_i^{1/2}}]$ is compact hence $\homo(M_{\phi_i^{1/2}}[T, M_{\phi_i^{1/2}}])=0$.
	\end{proof}
	
	\begin{lemma}\label{homconst}
		The formula
		\begin{equation}\label{homdef}
			\homo(T):=\sum_{i\in\II} \homo_i(M_{\phi_i^{1/2}}TM_{\phi_i^{1/2}})
		\end{equation}
		gives a well-defined mapping $\A\to\B$.
	\end{lemma}
	\begin{proof}
		This follows from Definition \ref{A} (3), Lemma \ref{convergy2} and the completeness of $\mathcal B$.
	\end{proof}
	
	From now, fix $\homo$ to be defined by (\ref{homdef}). It is clear that $\homo$ vanishes on compact operators, since each $\homo_i$ does.
	
	\begin{lemma}\label{multipliersout}
		For each $T\in\A$ and $\phi\in C_0(M)$, 
		$$\homo(TM_\phi)=\homo(T)\homo(M_\phi)=\homo(M_\phi)\homo(T)=\homo(M_\phi T)$$
	\end{lemma}
	\begin{proof}
		For each $i\in\II$, let $\psi_i\in C_c(\U_i)$ be such that $\psi_i\phi_i=\phi_i$. Then, in addition, $\psi_i\phi_i^{1/2}=\phi_i^{1/2}$ and so
		\begin{align*}
			M_{\phi_i^{1/2}}(TM_\phi) M_{\phi_i^{1/2}}&=M_{\phi_i^{1/2}}TM_\phi (M_{\psi_i}M_{\phi_i^{1/2}})\\
			&=(M_{\phi_i^{1/2}}TM_{\phi_i^{1/2}})M_{\phi\psi_i}
		\end{align*}
		Since $\homo_i$ is a $*$-homomorphism for all $i\in\II$, one gets
		$$\homo_i((M_{\phi_i^{1/2}}TM_{\phi_i^{1/2}}) M_{\phi\psi_i})=\homo_i(M_{\phi_i^{1/2}}TM_{\phi_i^{1/2}})\homo_i(M_{\phi\psi_i})$$
		Now,
		\begin{align*}
			\homo_i(M_{\phi\psi_i})&=\Homo(\phi\psi_i)=\Homo(\psi_i)\Homo(\phi)\\
			&=\homo_i(M_{\psi_i})\homo(M_\phi)
		\end{align*}
		Thus
		\begin{align*}
			\homo_i(M_{\phi_i^{1/2}}(TM_\phi) M_{\phi_i^{1/2}})&=\homo_i(M_{\phi_i^{1/2}}TM_{\phi_i^{1/2}})\homo_i(M_{\psi_i})\homo(M_\phi)\\
			&=\homo_i(M_{\phi_i^{1/2}}TM_{\phi_i^{1/2}}M_{\psi_i})\homo(M_\phi)\\
			&=\homo_i(M_{\phi_i^{1/2}}TM_{\phi_i^{1/2}})\homo(M_\phi)
		\end{align*}
		hence
		\begin{align*}
			\homo(TM_\phi)&=\sum_{i\in\II} \homo_i(M_{\phi_i^{1/2}}TM_\phi M_{\phi_i^{1/2}})\\
			&=\sum_{i\in\II} \homo_i(M_{\phi_i^{1/2}}T M_{\phi_i^{1/2}})\homo(M_\phi)\\
			&=\homo(T)\homo(M_\phi)
		\end{align*}
		An analogous argument shows that
		$$\homo(M_\phi T)=\homo(M_\phi)\cdot\homo(T),$$
		and finally, since $[T,M_\phi]$ is compact and $\homo$ vanishes on compact operators, one has
		$$\homo(M_\phi T)=\homo(TM_\phi+[T, M_\phi])=\homo(TM_\phi),$$
		which completes the proof.
	\end{proof}
	
	\begin{proposition}\label{globalisation-equal}
		For each $j\in\II$ and each $T\in\A_j$, one has $\homo(T)=\homo_j(T)$.
	\end{proposition}
	\begin{proof}
		For all $i\in\II$, the operator $M_{\phi_i^{1/2}}TM_{\phi_i^{1/2}}$ is compactly supported in both $\U_i\cap \U_j$ and hence lies in $\A_i\cap\A_j$ by Definition \ref{localalg} (3). By Definition \ref{localhomo} (2), one has
		\begin{equation}\label{ijinterchange}
			\homo_i(M_{\phi_i^{1/2}}TM_{\phi_i^{1/2}})=\homo_j(M_{\phi_i^{1/2}}TM_{\phi_i^{1/2}})
		\end{equation}
		for all $i\in\II$. Let $\phi\in C_c(\U_j)$ be such that $T=M_\phi T=M_\phi T$. Then, for $i\in\II$, since $\homo_j$ is a homomorphism, vanishes on compact operators and $[T, M_{\phi_i^{1/2}}]$ is compact,
		\begin{align*}
			\homo_j(M_{\phi_i^{1/2}}TM_{\phi_i^{1/2}})&=\homo_j(M_{\phi_i^{1/2}\phi}TM_{\phi\phi_i^{1/2}})\\
			&=\homo_j(TM_{\phi_i^{1/2}\phi}+[T, M_{\phi_i^{1/2}\phi}])\homo_j(M_{\phi\phi_i^{1/2}})\\
			&=\homo_j(T)\homo_j(M_{\phi_i^{1/2}\phi})\homo_j(M_{\phi\phi_i^{1/2}})\\
			&=\homo_j(T)\homo_j(M_{\phi_i\phi^2})
		\end{align*}
		Thus by (\ref{ijinterchange}), if $i\in\II$,
		$$\homo_i(M_{\phi_i^{1/2}}TM_{\phi_i^{1/2}})=\homo_j(T)\homo_j(M_{\phi_i\phi^2})$$
		hence
		\begin{align*}
			\homo(T)=\sum_{i\in\II} \homo_i(M_{\phi_i^{1/2}}TM_{\phi_i^{1/2}})=\sum_{i\in\II} \homo_j(T)\homo_j(M_{\phi_i\phi^2}).
		\end{align*}
		Since $\phi$ is compactly supported, the series $\sum_{i\in\II} \phi_i\phi^2$ converges uniformly to $\phi^2$, hence $\sum_{i\in\II} M_{\phi_i\phi^2}$ converges in norm. Applying continuity of products and of $\homo_j$, we get
		\begin{align*}
			\homo(T)&=\sum_{i\in\II} \homo_j(T)\homo_j(M_{\phi_i\phi^2})\\
			&=\homo_j(T)\sum_{i\in\II} \homo_j(M_{\phi_i\phi^2})\\
			&=\homo_j(T) \homo_j(\sum_{i\in\II}M_{\phi_i\phi^2})\\
			&=\homo_j(T) \homo_j(M_{\phi^2})
		\end{align*}
		Thus since $\homo_j$ is a homomorphism and $TM_{\phi^2}=T$,
		$$\homo(T)=\homo_j(TM_{\phi^2})=\homo_J(T)$$
		as desired.
	\end{proof}
	
	\begin{lemma}
		Let $T,S\in\A$ and suppose that $i\in\II$ and that $T$ is compactly supported in $\U_i$. Then
		$$\homo(TS)=\homo(T)\homo(S).$$
	\end{lemma}
	\begin{proof}
		Let $\phi\in C_c(\U_j)$ be such that $T=TM_\phi=M_\phi T$. Then
		\begin{align*}
			TS=(TM_\phi)S=T(M_{\phi^{1/2}} M_{\varphi^{1/2}}S)=T(M_{\phi^{1/2}} SM_{\varphi^{1/2}}+M_{\phi^{1/2}}[S, M_{\phi^{1/2}}]).
		\end{align*}
		Write $S_1=M_{\phi^{1/2}} SM_{\varphi^{1/2}}$ and $S_2=M_{\phi^{1/2}}[S, M_{\phi^{1/2}}]$ so that $TS=TS_1+TS_2$. Then $S_1\in\A_i$, and $T=M_\phi TM_\phi\in\A_i$, hence in particular $\homo(TS_1)=\homo_j(TS_1)=\homo(T)\homo(S_1)$. Additionally, $S_2$ is compact, hence so is $TS_2$ thus $\homo(TS_2)=0$. Therefore, by Lemma \ref{multipliersout}, since $T,M_\phi\in\A_j$,
		\begin{align*}
			\homo(TS)&=\homo(TS_1)+\homo(TS_2)=\homo(T)\homo(S_1)\\
			&=\homo(T)\homo(M_\phi)\homo(S)\\
			&=\homo(TM_\phi)\homo(S)\\
			&=\homo(TS)
		\end{align*}
		as desired.
	\end{proof}
	\begin{lemma}
		Let $T,S\in\A$. Then
		$$\homo(TS)=\homo(T)\homo(S).$$
	\end{lemma}
	\begin{proof}
		Let $i\in\II$. Then 
		\begin{align*}
			M_{\phi_i^{1/2}}TSM_{\phi_i^{1/2}}&=M_{\phi_i^{1/2}}TSM_{\phi_i^{1/4}}M_{\phi_i^{1/4}}\\
			&=M_{\phi_i^{1/2}}TM_{\phi_i^{1/4}}SM_{\phi_i^{1/4}}-M_{\phi_i^{1/2}}T[S,M_{\phi_i^{1/4}}]M_{\phi_i^{1/4}}
		\end{align*}
		Thus
		\begin{align*}
			\homo_i(M_{\phi_i^{1/2}}TSM_{\phi_i^{1/2}})&=\homo_i(M_{\phi_i^{1/2}}TM_{\phi_i^{1/4}}SM_{\phi_i^{1/4}})\\
			&=\homo_i(M_{\phi_i^{1/2}}TM_{\phi_i^{1/8}})\homo_i(M_{\phi_i^{1/8}}SM_{\phi_i^{1/4}})\\
			&=\homo_i(M_{\phi_i^{1/2}}TM_{\phi_i^{1/8}})\homo(M_{\phi_i^{1/8}}SM_{\phi_i^{1/4}})\\
			&=\homo_i(M_{\phi_i^{1/2}}TM_{\phi_i^{1/8}})\homo(M_{\phi_i^{1/8}})\homo(M_{\phi_i^{1/4}})\homo(S)\\
			&=\homo_i(M_{\phi_i^{1/2}}TM_{\phi_i^{1/8}})\homo_i(M_{\phi_i^{3/8}})\homo(S)\\
			&=\homo_i(M_{\phi_i^{1/2}}TM_{\phi_i^{1/2}})\homo(S)
		\end{align*}
		Therefore
		\begin{align*}
			\homo(TS)&=\sum_{i\in\II} \homo_i(M_{\phi_i^{1/2}}TSM_{\phi_i^{1/2}})\\
			&=\sum_{i\in\II} \homo_i(M_{\phi_i^{1/2}}TM_{\phi_i^{1/2}})\hom(S)\\
			&=\homo(T)\homo(S)
		\end{align*}
		as desired.
	\end{proof}
	
	\begin{proposition}\label{globalisation-compact}
		$\ker(\homo)=\K(L_2(M,\nu))$.
	\end{proposition}
	\begin{proof}
		Suppose that $T\in\ker(\homo)$. Recall that $T$ is compact if and only if $T^*T$ is compact. Since $\homo$ is a homomorphism, we have $\homo(T^*T)=0$. But
		$$\homo(T^*T)=\sum_{i\in\II}\homo_i(M_{\phi_i^{1/2}}T^*TM_{\phi_i^{1/2}})$$
		and each term in the sum on the right is positive, hence $\homo_i(M_{\phi_i^{1/2}}T^*TM_{\phi_i^{1/2}})=0$ for all $i\in\II$. This implies that $M_{\phi_i^{1/2}}T^*TM_{\phi_i^{1/2}}\in\ker(\homo_i)\subset\K(L_2(M,\nu))$ for all $i\in\II$. Since $M_{\phi_i^{1/2}}[T, M_{\phi_i^{1/2}}]$ is also compact, it follows that $M_{\phi_i}T^*T$ is compact and hence $\displaystyle\sum_{i\in\II} M_{\phi_i}T^*T=T^*T$ is compact.
	\end{proof}
	
	In sum, Proposition \ref{globalisation-algebra}, Proposition \ref{homconst}, Proposition \ref{globalisation-equal} and Proposition \ref{globalisation-compact} give the complete proof of Theorem \ref{globalisation}.

	\section{Principal symbol on $\GG$-filtered manifolds}\label{symbol on Carnot manifolds}
 
	In this section, we prove Theorem \ref{symbol on manifold} the existence of principal symbol mapping on $\GG$-filtered manifolds.
	
	%
	%

	\subsection{Construction of the principal symbol mapping}


	Before our proof of Theorem \ref{symbol on manifold}, we give some necessary preparations. For later use convenience, denote $H_{\Phi}(x)=(X_j\Phi_i(x))_{i,j=1}^{n_{1}}$ for a $C^{1}$-diffeomorphism $\Phi$ on some open set of $\GG$. Let $\Omega,\Omega',\Omega''\subset\GG$ be open sets. If $\Psi:\Omega\to\Omega'$ and $\Phi:\Omega'\to\Omega''$ are $C^{1}$-diffeomorphisms, then we have
	\begin{align}\label{horizontal Jacobian Leibniz-law}
		H_{\Phi\circ\Psi}(x)=H_{\Phi}(\Psi(x))H_{\Psi}(x)\quad\mbox{for}\quad x\in\Omega.
	\end{align}
	\begin{definition}
		Let $M$ be a $\GG$-filtered manifold, with $\GG$-atlas $\{(\U_{i},h_{i})\}_{i\in\II}$. Let $\Borel$ be the Borel $\sigma$-algebra on $M$ and let $\nu:\Borel\rightarrow\RR$ be a smooth positive density. Define an isometry $W_{i}:L_{2}(\U_{i},\nu)\rightarrow L_{2}(\Omega_{i},\nu\circ h_{i}\rp)$ by setting
		\begin{align*}
			W_{i}f=f\circ h_{i}\rp,\quad f\in L_{2}(\U_{i},\nu).
		\end{align*}
	\end{definition}
	This definition is well-defined because the mapping $h_{i}:(\U_{i},\nu)\rightarrow(\Omega_{i},\nu\circ h_{i}\rp)$ preserves the measure. If $T$ is an operator compactly supported in $\U_{i}$, then we have the following relation 
	$$L_{2}(\Omega_{i},\nu\circ h_{i}\rp) \stackrel{W_{i}\rp}{\longrightarrow}
	L_{2}(\U_{i},\nu) \stackrel{T}{\longrightarrow}
	L_{2}(\U_{i},\nu) \stackrel{W_{i}}{\longrightarrow}
	L_{2}(\Omega_{i},\nu\circ h_{i}\rp),$$
	i.e. $W_{i}TW_{i}\rp$ is understood as an element of the algebra $\B(L_{2}(\Omega_{i},\nu\circ h_{i}\rp))$. More precisely, this operator is compactly supported in $\Omega_{i}$.

	For notion $\Ext_{\Omega_{i}}$ below, we refer to Notation \ref{Ext-def on G}.
	\begin{definition}\label{local principal domain}
		Let $M$ be a $\GG$-filtered manifold, with $\GG$-atlas $\{(\U_{i},h_{i})\}_{i\in\II}$ and let $\nu$ be a smooth positive density on $M$. For each $i\in\II$, let 
		$$\Pi_{i}=\set{T\in\B(L_{2}(M,\nu)):T~is~ compactly~ supported~ in~ \U_{i}~and~ \Ext_{\Omega_{i}}(W_{i}TW_{i}\rp)\in\Pi}.$$
	\end{definition}
	For example, if $\phi\in C_{c}(\U_{i})$, the operator $M_{\phi}$  belongs to $\Pi_{i}$. In fact $\Pi_{i}$ is the domain $\Pi_{M}$ in the local coordinate $(\U_{i},h_{i})$.

	The following lemma verifies that the algebra $E_{hom}$ is a bundle of $C^{*}$-algebras.
	\begin{lemma}\label{form-mapping-triple}
		Let $M$ be a $\GG$-filtered manifold, with $\GG$-atlas $\{(\U_{i},h_{i})\}_{i\in\II}$ and let $\nu$ be a smooth positive density on $M$. For $i,j\in\II$, let the mappings $\pi_{i,j}:\U_{i}\cap\U_{j}\rightarrow\aut\left(C^{*}\Big(\{\quasiriesz{A}{k}:A\in\at\}_{k=1}^{n_{1}}\Big)\right)$ be defined by
		\begin{align*}
			\pi_{i,j}(x)=\pi_{H^{\Phi_{i,j}}(h_{i}(x))},\quad i,j\in\II \quad such\quad that\quad x\in\U_{i}\cap\U_{j}\neq\emptyset,
		\end{align*}
		where $\Phi_{i,j}$ is the transition mapping $h_j\circ h_i\rp$. Then the family $\omega=\{\pi_{i,j}\}_{i,j\in\II}$ satisfies the cocycle conditions \eqref{cocycle-two} and \eqref{cocycle-three}.
	\end{lemma}
	\begin{proof}
		Let us verify \eqref{cocycle-two} at first. Let $i,j\in\II$ be such that $\U_{i}\cap\U_{j}\neq\emptyset$. For $i,j\in\II$, by $\Phi_{i,j}=h_j\circ h_i\rp$, it is clear that
		$$\Phi_{j,i}\circ\Phi_{i,j}(\xi)=\xi \quad\mbox{for}\quad\xi\in\Omega_{j,i}.$$
		Since $M$ is a smooth manifold, by \eqref{horizontal Jacobian Leibniz-law}, one has
		\begin{align*}
			H_{\Phi_{j,i}\circ\Phi_{i,j}}(\xi)=H_{\Phi_{j,i}}(\Phi_{i,j}(\xi))H_{\Phi_{i,j}}(\xi).
		\end{align*}
		Recalling Remark \ref{basis-rep remark}, we have, $H_{\id}(\xi)=\II_{n_{1}}$ for any $\xi\in\GG$. In other words, substituting $\xi$ as $h_{i}(x)$ with $x\in\U_{i}\cap\U_{j}$, it follows that
		\begin{align*}
			H_{\Phi_{j,i}}(h_{j}(x))H_{\Phi_{i,j}}(h_{i}(x))=H_{\Phi_{j,i}}(\Phi_{i,j}(h_{i}(x))) H_{\Phi_{i,j}}(h_{i}(x))
			=H_{\Phi_{j,i}\circ\Phi_{i,j}}(h_{i}(x))=H_{\id}(h_{i}(x))=\II_{n_{1}}.
		\end{align*}
		Moreover, by Theorem \ref{admissibility}, there is an inverse matrix $b_{rest}(x)\in{\rm GL}(d-n_{1},\RR)$ such that
		\begin{align*}
			H^{\Phi_{i,j}}(h_{i}(x))H^{\Phi_{j,i}}(h_{j}(x))
			=\diag\left(\left(H_{\Phi_{j,i}}(h_{j}(x))H_{\Phi_{i,j}}(h_{i}(x))\right)^{t},b_{rest}(x)\right)
			=\diag(\II_{n_{1}},b_{rest}(x)).
		\end{align*}
		For $A\in\at$, decompose $A$ into the form of \eqref{diagonal}, i.e. there is an inverse matrix $A_{rest}\in{\rm GL}(d-n_{1},\RR)$ such that
		$A=\diag(A_{n_{1}},A_{rest})$.
		Therefore,
		\begin{align}\label{HH equal I}
			H^{\Phi_{i,j}}(h_{i}(x))H^{\Phi_{j,i}}(h_{j}(x))A
			=\diag\left(A_{n_{1}},b_{rest}(x)A_{rest} \right).
		\end{align}
		For $k\in\{1,\,\dotsc,n_{1}\}$ and $A\in\at$, using Lemma \ref{pffacts}, $V_{A}V_{B}=V_{BA}$ and \eqref{HH equal I}, one obtains
		\begin{align*}
			\pi_{i,j}(x)\pi_{j,i}(x)(\quasiriesz{A}{k})
			&=\pi_{H^{\Phi_{i,j}}(h_{i}(x))}\pi_{H^{\Phi_{j,i}}(h_{j}(x))}(\quasiriesz{A}{k})\\
			&=V_{H^{\Phi_{i,j}}(h_{i}(x))}\rp V_{H^{\Phi_{j,i}}(h_{j}(x))}\rp V_{A}\rp R_{k} V_{A} V_{H^{\Phi_{j,i}}(h_{j}(x))} V_{H^{\Phi_{i,j}}(h_{i}(x))}\\
			&=V_{H^{\Phi_{i,j}}(h_{i}(x))H^{\Phi_{j,i}}(h_{j}(x))A}\rp R_{k} V_{H^{\Phi_{i,j}}(h_{i}(x))H^{\Phi_{j,i}}(h_{j}(x))A}\\
			&=\quasiriesz{A}{k}.
		\end{align*}
		Since $C^{*}\Big(\{\quasiriesz{A}{k}:A\in\at\}_{k=1}^{n_{1}}\Big)$ is generated by $\quasiriesz{A}{k}$ with $k\in\{1,\,\dotsc,n_{1}\}$ and $A\in\at$, it follows that, for any $T\in C^{*}\Big(\{\quasiriesz{A}{k}:A\in\at\}_{k=1}^{n_{1}}\Big)$, $$\pi_{i,j}(x)\pi_{j,i}(x)(T)=T.$$ 
		This is the desired equality in \eqref{cocycle-two}.

		\bigskip

		Next let us verify \eqref{cocycle-three}. Let $i,j,k\in\II$ be such that $\U_{i}\cap\U_{j}\cap\U_{k}\neq\emptyset$. Obviously,
		$$\Phi_{k,i}\circ\Phi_{j,k}\circ\Phi_{i,j}(\xi)=\xi \quad\mbox{for}\quad \xi\in\Omega_{i,j}.$$
		Applying \eqref{horizontal Jacobian Leibniz-law} again, it yields that, for $\xi\in\Omega_{i,j}\cap\Omega_{i,k}$,
		\begin{align*}
			H_{\Phi_{k,i}}(\Phi_{i.k}(\xi))H_{\Phi_{j,k}}(\Phi_{i,j}(\xi))H_{\Phi_{i,j}}(\xi)
			=H_{\Phi_{k,i}\circ\Phi_{j,k}\circ\Phi_{i,j}}(\xi)
			=\II_{n_{1}}.
		\end{align*}
		Therefore, replacing the $\xi$ by $h_{i}(x)$ with $x\in\U_{i}\cap\U_{j}\cap\U_{k}$, there is a matrix $b'_{rest}(x)\in{\rm GL}(d-n_{1},\RR)$ such that
		\begin{align*}
			H^{\Phi_{i,j}}(h_{i}(x))H^{\Phi_{j,k}}(h_{j}(x))H^{\Phi_{k,i}}(h_{k}(x))A
			&=\diag\left(\left(H_{\Phi_{k,i}}(h_{k}(x))H_{\Phi_{j,k}}(h_{j}(x))H_{\Phi_{i,j}}(h_{i}(x))\right)^{t}A_{n_{1}},b'_{rest}(x)A_{rest} \right)\\
			&=\diag(A_{n_{1}},b'_{rest}(x)A_{rest}).
		\end{align*}
		Moreover, by Lemma \ref{pffacts} and $V_{A}V_{B}=V_{BA}$, one has
		\begin{align*}
			\pi_{i,j}(x)\pi_{j,k}(x)\pi_{k,i}(x)(\quasiriesz{A}{k})
			&=\pi_{H^{\Phi_{i,j}}(h_{i}(x))}\pi_{H^{\Phi_{j,k}}(h_{j}(x))}\pi_{H^{\Phi_{k,i}}(h_{k}(x))}(\quasiriesz{A}{k})\\
			&=V_{H^{\Phi_{i,j}}(h_{i}(x))}\rp V_{H^{\Phi_{j,k}}(h_{j}(x))}\rp V_{H^{\Phi_{k,i}}(h_{k}(x))}\rp V_{A}\rp R_{k} V_{A} V_{H^{\Phi_{k,i}}(h_{k}(x))} V_{H^{\Phi_{j,k}}(h_{j}(x))} V_{H^{\Phi_{i,j}}(h_{i}(x))}\\
			&=V_{H^{\Phi_{i,j}}(h_{i}(x))H^{\Phi_{j,k}}(h_{j}(x))H^{\Phi_{k,i}}(h_{k}(x))A}\rp R_{k} V_{H^{\Phi_{i,j}}(h_{i}(x))H^{\Phi_{j,k}}(h_{j}(x))H^{\Phi_{k,i}}(h_{k}(x))A}\\
			&=\quasiriesz{A}{k}.
		\end{align*}
		This claims \eqref{cocycle-three} on the generators of $C^{*}\Big(\{\quasiriesz{A}{k}:A\in\at\}_{k=1}^{n_{1}}\Big)$ and then the density gives the desired result.

	\end{proof}

	
	For the notion $\sym$ we refer to Theorem \ref{symbol}. The mapping $\Theta_{i}$ used in Definition \ref{local principal symbol} is defined in \eqref{natural embedding}.
	\begin{definition}\label{local principal symbol}
		Let $M$ be a $\GG$-filtered manifold, with $\GG$-atlas $\{(\U_{i},h_{i})\}_{i\in\II}$ and let $\nu$ be a smooth positive density on $M$. For each $i\in\II$, the mapping $\sym_{i}:\Pi_{i}\rightarrow C_{b}(E_{hom})$ is defined by the formula
		\begin{align*}
			\sym_{i}(T)=\Theta_{i}(\sym(\Ext_{\Omega_{i}}(W_{i}TW_{i}\rp))\circ h_{i}),\quad T\in\Pi_{i}.
		\end{align*}
	\end{definition}

	\begin{lemma}\label{subalgebra-homomorphism}
		Let $M$ be a $\GG$-filtered manifold, with $\GG$-atlas $\{(\U_{i},h_{i})\}_{i\in\II}$ and let $\nu$ be a smooth positive density on $M$. For each $i\in\II$, we have
		\begin{enumerate}[\rm(1)]
			\item $\Pi_{i}$ is a $*$-subalgebra in $\B(L_{2}(M,\nu))$;
			\item $\sym_{i}:\Pi_{i}\rightarrow C_{b}(E_{hom})$ is a $*$-homomorphism.
		\end{enumerate}
	\end{lemma}
	\begin{proof}
		The proof of (1): we need to show that $\Pi_{i}$ is a subalgebra of $\B(L_{2}(M,\nu))$ and it is closed in taking adjoint which is associated to the measure $\nu$ defined on $M$. 
		For $T_{1},T_{2}\in\Pi_{i}$, it is clear that $T_{1}T_{2}$ is compactly supported in $\U_{i}$. Moreover,
		\begin{align*}
			\Ext_{\Omega_{i}}(W_{i}T_{1}T_{2} W_{i}\rp)=\Ext_{\Omega_{i}}(W_{i}T_{1} W_{i}\rp)\Ext_{\Omega_{i}}(W_{i}T_{2} W_{i}\rp)\in\Pi.
		\end{align*}
		This implies that $\Pi_{i}$ is a subalgebra.
		
		Let $T\in\Pi_{i}$. Since $T$ is compactly supported in $\U_{i}$, it follows that $(W_{i}T W_{i}\rp)^{*}$ is compactly supported in $\Omega_{i}$, i.e. there is a $\phi\in C_{c}(\Omega_{i})$ such that
		\begin{align}\label{WiTWi compact support}
			(W_{i}T W_{i}\rp)^{*}=M_{\phi}(W_{i}T W_{i}\rp)^{*}=(W_{i}T W_{i}\rp)^{*}M_{\phi}.
		\end{align}
		Here the adjoint of $W_{i}T W_{i}\rp$ is taken with respect to the Lebesgue measure on $\GG$.
		Letting $\psi_{i}$ be the Radon-Nikodym derivative of $\nu\circ h_{i}\rp$ on $\Omega_{i}$, it follows that $\psi_{i}$ and $\psi_{i}\rp$ are smooth on $\Omega_{i}$. Thus, for two testing functions $f,g$ on $\Omega_{i}$ with respect to measure $d\nu\circ h_{i}\rp$, 
		\begin{align*}
			\int_{\Omega_{i}}(W_{i}T^{*}W_{i}\rp f)(w)g(w) \,d\nu\circ h_{i}\rp(w)
			&=\int_{\Omega_{i}}T^{*}(f\circ h_{i})(h_{i}\rp(w))g(w)\,d\nu\circ h_{i}\rp(w)\\
			&=\int_{\Omega_{i}}(f\circ h_{i})(h_{i}\rp(w)) T(g\circ h_{i})(h_{i}\rp(w))\,d\nu\circ h_{i}\rp(w)\\
			&=\int_{\Omega_{i}}f(w) (W_{i}TW_{i}\rp g)(w)\,d\nu\circ h_{i}\rp(w)\\
			&=\int_{\Omega_{i}}M_{\psi_{i}\rp} (W_{i}TW_{i}\rp)^{*} M_{\psi_{i}}f(w) g(w)\,d\nu\circ h_{i}\rp(w),
		\end{align*}
		i.e.
		\begin{align*}
			W_{i}T^{*} W_{i}\rp=M_{\psi_{i}\rp}(W_{i}T W_{i}\rp)^{*}M_{\psi_{i}}.
		\end{align*}
		Combining this equality with \eqref{WiTWi compact support}, one obtains
		\begin{align*}
			W_{i}T^{*} W_{i}\rp=M_{\psi_{i}\rp\phi}(W_{i}T W_{i}\rp)^{*}M_{\psi_{i}\phi}.
		\end{align*}
		Note that $\psi_{i}\phi,\psi_{i}\rp\phi\in C_{c}(\Omega_{i})$ and the adjoint of $W_{i}T W_{i}\rp$ is taken with respect to the Lebesgue measure on $\GG$. Therefore,
		\begin{align}\label{*-Ext-commutative}
			\Ext_{\Omega_{i}}(W_{i}T^{*} W_{i}\rp)
			=M_{\psi_{i}\rp\phi}\left(\Ext_{\Omega_{i}}(W_{i}T W_{i}\rp)\right)^{*}M_{\psi_{i}\phi}\in\Pi.
		\end{align}
		This implies that $T^{*}\in\Pi_{i}$.

		\bigskip
		
		The proof of (2): it is clear that $\sym_{i}$ is a homomorphism and it is left to show it is invariant in taking adjoint. By \eqref{*-Ext-commutative} and property of the $*$-homomorphism $\sym$, one has
		\begin{align*}
			\sym(\Ext_{\Omega_{i}}(W_{i}T^{*} W_{i}\rp))
			&=\sym(M_{\psi_{i}\rp\phi})\cdot \sym(\Ext_{\Omega_{i}}(W_{i}T W_{i}\rp)^{*})\cdot \sym(M_{\psi_{i}\phi})\\
			&=\sym(M_{\phi^{2}})\cdot \sym(\Ext_{\Omega_{i}}(W_{i}T W_{i}\rp))^{*}\\
			&=\sym(\Ext_{\Omega_{i}}(M_{\phi^{2}}W_{i}T W_{i}\rp))^{*}\\
			&=\sym(\Ext_{\Omega_{i}}(W_{i}T W_{i}\rp))^{*},
		\end{align*}
		where, in the last equality, we use
		$$M_{\phi^{2}}W_{i}T W_{i}\rp=W_{i}T W_{i}\rp.$$
		Since $\Theta_{i}$ is a $*$-homomorphism, it follows that
		\begin{align*}
			\sym_{i}(T^{*})&=\Theta_{i}(\sym(\Ext_{\Omega_{i}}(W_{i}T^{*} W_{i}\rp))\circ h_{i})\\
			&=\Theta_{i}(\sym(\Ext_{\Omega_{i}}(W_{i}T W_{i}\rp))^{*}\circ h_{i})\\
			&=\left(\Theta_{i}(\sym(\Ext_{\Omega_{i}}(W_{i}T W_{i}\rp))\circ h_{i})\right)^{*}\\
			&=\sym_{i}(T)^{*}.
		\end{align*}
		This is the desired result.
	\end{proof}

	\begin{lemma}\label{symi=symj}
		Let $M$ be a $\GG$-filtered manifold, with $\GG$-atlas $\{(\U_{i},h_{i})\}_{i\in\II}$ and let $\nu$ be a smooth positive density on $M$. Let $(\U_{i},h_{i})$ and $(\U_{j},h_{j})$ be two charts and let $T\in\B(L_{2}(M,\nu))$ be compactly supported in $\U_{i}\cap\U_{j}$. If $T\in\Pi_{i}$, then $T\in\Pi_{j}$ and $\sym_{i}(T)=\sym_{j}(T)$.
	\end{lemma}
	\begin{proof}
		Denote $$T_{i}=\Ext_{\Omega_{i}}(W_{i}T W_{i}\rp),\quad T_{j}=\Ext_{\Omega_{j}}(W_{j}T W_{j}\rp).$$
		Noting that $W_{j}=V_{\Phi_{i,j}}\rp W_{i}$, one has $$T_{j}=\Ext_{\Omega_{j}}(V_{\Phi_{i,j}}\rp W_{i}T W_{i}\rp V_{\Phi_{i,j}}),$$ 
		and therefore
		\begin{align}\label{Rest Tj-Ti}
			\Rest_{\Omega_{j}}(T_{j})=
			U_{\Phi_{i,j}}\rp M_{|J_{\Phi_{i,j}}|^{1/2}}\Rest_{\Omega_{i}}(T_{i}) M_{|J_{\Phi_{i,j}}|^{-1/2}} U_{\Phi_{i,j}}.
		\end{align}
		Since $T$ is compactly supported in $\U_{i}\cap\U_{j}$, it follows that $T_{i}$ is compactly supported in $\Omega_{i,j}$. Let $\Omega'\subset\Omega_{i,j}$ be compact such that
		$$M_{\chi_{\Omega'}}T_{i}M_{\chi_{\Omega'}}=W_{i}T W_{i}\rp.$$
		Using Tietze extension theorem, choose $\phi\in C(\GG)$ such that $\phi\rp\in C(\GG)$ and $\phi=|J_{\Phi_{i,j}}|^{1/2}$ on $\Omega'$. One obtains
		\begin{align*}
			\Ext_{\Omega_{i}}(M_{|J_{\Phi_{i,j}}|^{1/2}}\Rest_{\Omega_{i}}(T_{i}) M_{|J_{\Phi_{i,j}}|^{-1/2}})=M_{\phi}T_{i} M_{\phi\rp}=T_{i}+[M_{\phi},T_{i}]M_{\phi\rp}.
		\end{align*}
		Since $T_{i}\in\Pi$, it follows from Lemma \ref{compactness} that $[M_{\phi},T_{i}]M_{\phi\rp}\in\K(L_{2}(\GG))$. Thus,
		\begin{align*}
			\Ext_{\Omega_{i}}(M_{|J_{\Phi_{i,j}}|^{ 1/2}}\Rest_{\Omega_{i}}(T_{i}) M_{|J_{\Phi_{i,j}}|^{- 1/2}})\in T_{i}+\K(L_{2}(\GG)).
		\end{align*}
		Moreover,
		\begin{align}\label{include Tj-Ti }
			M_{|J_{\Phi_{i,j}}|^{1/2}}\Rest_{\Omega_{i}}(T_{i}) M_{|J_{\Phi_{i,j}}|^{-1/2}}\in \Rest_{\Omega_{i}}(T_{i})+\K(L_{2}(\Omega_{i})).
		\end{align}
		Observing that $\Phi_{i,j}:\Omega_{i,j}\to\Omega_{j,i}$ is a $C^{\infty}$ diffeomorphism, it follows that $U_{\Phi_{i,j}}:L_{2}(\Omega_{j,i})\to L_{2}(\Omega_{i,j})$ is unitary and therefore $$U_{\Phi_{i,j}}\rp \K(L_{2}(\Omega_{i}))U_{\Phi_{i,j}}\subset \K(L_{2}(\Omega_{j})).$$ Combining \eqref{Rest Tj-Ti} and \eqref{include Tj-Ti }, one has
		\begin{align*}
			\Rest_{\Omega_{j}}(T_{j})\in U_{\Phi_{i,j}}\rp \Rest_{\Omega_{i}}(T_{i})U_{\Phi_{i,j}}+\K(L_{2}(\Omega_{j})).
		\end{align*}
		Furthermore,
		\begin{align}\label{Tj-Ext-Rest-Ti}
			T_{j}\in \Ext_{\Omega_{j}}(U_{\Phi_{i,j}}\rp \Rest_{\Omega_{i}}(T_{i})U_{\Phi_{i,j}})+\K(L_{2}(\GG)).
		\end{align}
		Since $\{(\U_{i},h_{i})\}_{i\in\II}$ is a $\GG$-atlas, by Theorem \ref{Ext-Rest invariant}, the above first factor on the right hand side belongs to $\Pi$. By Theorem \ref{symbol}, $\K(L_{2}(\GG))$ belongs to $\Pi$. Therefore, $T_{j}\in\Pi$.
		
		Additionally, by Theorem \ref{Ext-Rest invariant},
		\begin{align*}
			\sym\left(\Ext_{\Omega_{j}}(U_{\Phi_{i,j}}\rp \Rest_{\Omega_{i}}(T_{i})U_{\Phi_{i,j}})\right)
			=\pi_{H^{\Phi_{i,j}}}(\sym(T_{i}))\circ\Phi_{i,j}\rp.
		\end{align*}
		Combining this equality with \eqref{Tj-Ext-Rest-Ti}, we obtain
		\begin{align*}
			\sym(T_{j})\circ h_{j}
			&=\pi_{H^{\Phi_{i,j}}}(\sym(T_{i}))\circ\Phi_{i,j}\rp\circ h_{j}\\
			&=\pi_{H^{\Phi_{i,j}}}(\sym(T_{i}))\circ h_{i}\\
			&=\pi_{H^{\Phi_{i,j}}\circ h_{i}}(\sym(T_{i})\circ h_{i})\\
			&=\pi_{i,j}(\sym(T_{i})\circ h_{i}).
		\end{align*}
		Hence,
		\begin{align*}
			\sym_{j}(T)=\Theta_{j}(\sym(T_{j})\circ h_{j})
			=\Theta_{j}\circ\pi_{i,j}(\sym(T_{i})\circ h_{i})
			=\Theta_{i}(\sym(T_{i})\circ h_{i})=\sym_{i}(T).
		\end{align*}
		This is the desired result.
	\end{proof}

	The following theorem states that the collection $\{\Pi_{i}\}_{i\in\II}$ of $*$-algebras and the collection $\{\sym_{i}\}_{i\in\II}$ of $*$-homomorphisms satisfy the conditions in Theorem \ref{globalisation}.
	\begin{theorem}\label{Pii-symi satisfy local}
		Let $M$ be a $\GG$-filtered manifold, with $\GG$-atlas $\{(\U_{i},h_{i})\}_{i\in\II}$ and let $\nu$ be a smooth positive density on $M$. Then
		\begin{enumerate}[\rm(1)]
			\item\label{Pii satisfy local} \itemcase{t}he collection $\{\Pi_{i}\}_{i\in\II}$ introduced in Definition \ref{local principal domain} satisfies the conditions in Definition \ref{localalg};
			\item\label{symi satisfy local} \itemcase{t}he collection $\{\sym_{i}\}_{i\in\II}$ introduced in Definition \ref{local principal symbol} satisfies the conditions in Definition \ref{localhomo}.
		\end{enumerate}
	\end{theorem}
	\begin{proof}
		Proof of (1): the condition (1) in Definition \ref{localalg} is claimed in Lemma \ref{subalgebra-homomorphism}. 
		The condition (2) in Definition \ref{localalg} is contained in Definition \ref{local principal domain}. 
		The condition (3) in Definition \ref{localalg} is claimed in Lemma \ref{symi=symj}.
		
		For condition (4) in Definition \ref{localalg}, let $T\in\K(L_{2}(M,\nu))$ be compactly supported in $\U_{i}$ and let $\varphi\in C_{c}^{\infty}(\Omega_{i})$ be such that $W_{i}T W_{i}\rp=M_{\varphi}W_{i}T W_{i}\rp M_{\varphi}$. Then $W_{i}T W_{i}\rp$ is compact on $L_{2}(\Omega_{i},\nu\circ h_{i}\rp)$. So that there are operators $\{T_{n}\}_{n}$ with finite rank such that $T_{n}$ tends to $W_{i}TW_{i}\rp$ in the norm of $\B(L_{2}(\Omega_{i},\nu\circ h_{i}\rp))$. Let us consider the operators $\{M_{\varphi}T_{n}M_{\varphi}\}_{n}$.
		Letting $\psi_{i}$ be the Radon-Nikodym derivative of $\nu\circ h_{i}\rp$ on $\Omega_{i}$, for $f\in L_{2}(\Omega_{i})$, one has
		\begin{align*}
			\Norm{M_{\varphi}T_{n}M_{\varphi}f}_{L_{2}(\Omega_{i})}^{2}
			&=\int_{\Omega_{i}}\left|\varphi(x)T_{n}M_{\varphi}f(x)\right|^{2}\psi_{i}(x)\psi_{i}(x)\rp \,dx\\
			&\leq\Norm{\varphi^{2}\psi_{i}\rp}_{\infty}\int_{\Omega_{i}}\left|T_{n}M_{\varphi}f(x)\right|^{2}\psi_{i}(x)\,dx\\
			&\leq\Norm{\varphi^{2}\psi_{i}\rp}_{\infty}\Norm{\varphi^{2}\psi_{i}}_{\infty}\Norm{T_{n}}_{\B(L_{2}(\Omega_{i},\nu\circ h_{i}\rp))}\Norm{f}_{L_{2}(\Omega_{i})},
		\end{align*}
		i.e. $M_{\varphi}T_{n}M_{\varphi}$ belongs to $\B(L_{2}(\Omega_{i}))$.
		Moreover,
		\begin{align*}
			\Norm{(M_{\varphi}T_{n}M_{\varphi}-W_{i}TW_{i}\rp)f}_{L_{2}(\Omega_{i})}^{2}
			&=\int_{\Omega_{i}}\left|\varphi(x)(T_{n}-W_{i}TW_{i}\rp)M_{\varphi}f(x)\right|^{2}\psi_{i}(x)\psi_{i}(x)\rp \,dx\\
			&\leq\Norm{\varphi^{2}\psi_{i}\rp}_{\infty}\Norm{\varphi^{2}\psi_{i}}_{\infty}\Norm{(T_{n}-W_{i}TW_{i}\rp)}_{\B(L_{2}(\Omega_{i},\nu\circ h_{i}\rp))}\Norm{f}_{L_{2}(\Omega_{i})}.
		\end{align*}
		This implies that $W_{i}TW_{i}\rp$ is the limit of operators $\{M_{\varphi}T_{n}M_{\varphi}\}_{n}$ in the norm of $\B(L_{2}(\Omega_{i}))$.
		Thus, $\Ext_{\Omega_{i}}(W_{i}T W_{i}\rp)\in\K(L_{2}(\GG))$.  By Theorem \ref{symbol}, one obtains $\Ext_{\Omega_{i}}(W_{i}T W_{i}\rp)\in\Pi$ and therefore $T\in\Pi_{i}$.
		
		The condition (5) in Definition \ref{localalg} is immediate.
		
		For condition (6) in Definition \ref{localalg}, let $\phi\in C_{c}(\U_{i})$ and suppose that $\{T_{k}\}_{k\geq1}\subset M_{\phi}\Pi_{i} M_{\phi}$ are such that $T_{k}\to T$ in the uniform norm of $\B(L_{2}(M,\nu))$. It is clear that $T$ is compactly supported in $\U_{i}$. 
		Observe that $$\Norm{W_{i}T_{k}W_{i}\rp-W_{i}TW_{i}\rp}_{\B(L_{2}(\Omega_{i},\nu\circ h_{i}\rp))}=\Norm{T_{k}-T}_{\B(L_{2}(M,\nu))}$$ 
		and that the Radon-Nikodym derivative of $\nu\circ h_{i}\rp$ is strictly positive and smooth on $\Omega_{i}$.
		Since $T_{k}$ and $T$ are compactly supported in $\supp(\phi)\subset\U_{i}$, one gets $$\Ext_{\Omega_{i}}(W_{i}T_{k}W_{i}\rp)\to\Ext_{\Omega_{i}}(W_{i}TW_{i}\rp),\quad k\to\infty,$$ in the uniform norm of $\B(\lp{2})$. The left factor belongs to $\Pi$, so is its limit. In other words, $T\in\Pi_{i}$.
		
		For condition (7) in Definition \ref{localalg}, let $\phi\in C_{c}(\U_{i})$ and $T\in\Pi_{i}$. Note that
		\begin{align*}
			\Ext_{\Omega_{i}}(W_{i}[T,M_{\phi}]W_{i}\rp)
			&=\Ext_{\Omega_{i}}([W_{i}TW_{i}\rp,W_{i}M_{\phi}W_{i}\rp])\\
			&=[\Ext_{\Omega_{i}}(W_{i}TW_{i}\rp),M_{\phi\circ h_{i}\rp}].
		\end{align*}
		Observing that $\Ext_{\Omega_{i}}(W_{i}TW_{i}\rp)\in\Pi$ and $\phi\circ h_{i}\rp\in C_{c}(\GG)$, by Lemma \ref{compactness}, the last factor belongs to $\K(L_{2}(\GG))$. Therefore, the left factor is compact on $\lp{2}$. Since $\phi$ and $T$ are compactly supported in $\U_{i}$, and the Radon-Nikodym derivative of $\nu\circ h_{i}\rp$ is strictly positive and smooth on $\Omega_{i}$, it implies that $[T,M_{\phi}]\in\K(L_{2}(M,\nu))$.

		\bigskip
		Proof of (2): the condition (1) in Definition \ref{localhomo} is claimed in Lemma \ref{subalgebra-homomorphism} and
		the condition (2) in Definition \ref{localhomo} is claimed in Lemma \ref{symi=symj}.

		For condition (3) in Definition \ref{localhomo}, let $T\in\Pi_{i}$ be compact. As in the previous, $\Ext_{\Omega_{i}}(W_{i}T W_{i}\rp)$ is compact on $L_{2}(\GG)$. Since $\ker(\sym)=\K(L_{2}(\GG))$, it follows that 
		\begin{align*}
			\sym_{i}(T)=\Theta_{i}(\sym(\Ext_{\Omega_{i}}(W_{i}T W_{i}\rp))\circ h_{i})=\Theta_{i}(0)=0.
		\end{align*}
		Conversely, let $T\in\Pi_{i}$ satisfy $\sym_{i}(T)=0$. Observing that $\Theta_{i}:C_{c}(\U_{i},\M)\to C_{b}(E_{hom})$ is an embedding, one has $$\sym(\Ext_{\Omega_{i}}(W_{i}T W_{i}\rp))=0.$$
		Moreover, by $\ker(\sym)=\K(L_{2}(\GG))$, $$\Ext_{\Omega_{i}}(W_{i}T W_{i}\rp)\in\K(L_{2}(\GG)).$$
		Since $T$ is compactly supported in $\U_{i}$ and the Radon-Nikodym derivative of $\nu\circ h_{i}\rp$ is strictly positive and smooth on $\Omega_{i}$,
		this implies that $T\in\K(L_{2}(M,v))$. 
		
		The condition (4) in Definition \ref{localhomo} is immediate if we take ${\rm Hom}$ to be the natural embedding $C_{0}(M)\to C_{b}(E_{hom})$, which sends each $f\in C_{0}(M)$ to a family $F_{f}=\set{F_{i}\in C(\U_{i},\A_{2})}_{i\in\II}$ given by $F_{i}=f|_{\U_{i}}\cdot\one$, where $\one$ is the unit in $\A_{2}$.
		This completes the proof. 
	\end{proof}

	\subsection{Existence of principal symbol on manifolds }

	\begin{proof}[Proof of Theorem \ref{symbol on manifold}]
		Under the assumptions in Definition \ref{Pi-X} \eqref{Pi tranfer to group}, it is clear that $M_{\phi}TM_{\phi}$ is compactly supported in $\U_{i}$ and $\Ext_{\Omega_{i}}(W_{i}M_{\phi}TM_{\phi}W_{i}\rp)\in\Pi$, i.e. $M_{\phi}TM_{\phi}\in\Pi_{i}$. By Theorem \ref{Pii-symi satisfy local} \eqref{Pii satisfy local} and Theorem \ref{globalisation} \eqref{Cstar algebra}, $\Pi_{M}$ is a $C^{*}$-algebra which contains $\Pi_i$ for every $i\in\II$ and $\K(L_2(M,\nu))$. 
		
		Let the principal symbol mapping $\sym_{M}:\Pi_{M}\rightarrow C_{b}(E_{hom})$ be the $*$-homomorphism constructed in Theorem \ref{globalisation} from the collection $\{\sym_{i}\}_{i\in\II}$ defined in Definition \ref{local principal symbol}.
		The Theorem \ref{Pii-symi satisfy local} \eqref{symi satisfy local} and Theorem \ref{globalisation} \eqref{star homomorphism} imply the existence of the $*$-homomorphism $\sym_{M}:\Pi_{M}\to C_{b}(E_{hom})$ and $$\ker(\sym_{M})=\K(L_{2}(M,\nu)).$$
		It suffices to show that it is surjective.
		
		Let $F\in C_{b}(E_{hom})$ and let $\{\phi_{k}\}_{k\in\II}$ be a locally finite partition of unity subordinate to the atlas $\{\U_{k}\}_{k\in\II}$. Assuming that $\phi_{k}\in C_{c}(\U_{i_{k}})$, one has the $i_{k}$-th component $$q_{k}:=(F\phi_{k})_{i_{k}}\in C_{c}(\U_{i_{k}},\A_{2}).$$
		Since $\sym$ is a surjection, there exists a $T_{k}\in\Pi$ such that $\sym(T_{k})=q_{k}\circ h_{i_{k}}\rp$. Let $\psi_{k}\in C_{c}(\Omega_{i_{k}})$ be such that $\phi_{k}\circ h_{i_{k}}\rp=(\phi_{k}\circ h_{i_{k}}\rp)\psi_{k}$. Then $M_{\psi_{k}}T_{k}M_{\psi_{k}}\in\Pi$ and $$\sym(M_{\psi_{k}}T_{k}M_{\psi_{k}})=(q_{k}\circ h_{i_{k}}\rp)\psi_{k}^{2}=q_{k}\circ h_{i_{k}}\rp.$$
		Note that $M_{\psi_{k}}T_{k}M_{\psi_{k}}$ is bounded and compactly supported in $\Omega_{i_{k}}$. It follows that $T_{F,k}:=W_{i_{k}}\rp \Rest_{\Omega_{i_{k}}}(T_{k}) W_{i_{k}}$ is bounded and compactly supported in $\U_{i_{k}}$. Since $\Ext_{\Omega_{i_{k}}}(W_{i_{k}}T_{F,k}W_{i_{k}}\rp)=T_{k}$, it yields that $$T_{F,k}\in\Pi_{i_{k}}\subset\Pi_{M}$$ and
		\begin{align*}
			\sym_{M}(T_{F,k})=\sym_{i_{k}}(T_{F,k})
			&=\Theta_{i_{k}}(\sym(\Ext_{\Omega_{i_{k}}}(W_{i_{k}}T_{F,k}W_{i_{k}}\rp))\circ h_{i_{k}})\\
			&=\Theta_{i_{k}}(\sym(T_{k})\circ h_{i_{k}})\\
			&=\Theta_{i_{k}}(q_{k})\\
			&=\Theta_{i_{k}}((F\phi_{k})_{i_{k}})=F\phi_{k}.
		\end{align*} 
		Therefore, $\displaystyle T_{F}:=\sum_{k\in\II}T_{F,k}\in\Pi_{M}$ and $$\sym_{M}(T_{F})=\sum_{k\in\II}\sym(T_{F,k})=\sum_{k\in\II}F\phi_{k}=F.$$
		This implies surjectivity of $\sym_{M}$, which completes the proof.
	\end{proof}

	\bigskip
	
	{\bf Acknowledgments:} D. Farrell, F. Sukochev  and D. Zanin are supported by ARC DP230100434. 
 F. L. Yang is supported by National Natural Science Foundation of China No. 12371138 and he is grateful for the support of his supervisor professor Xiao Xiong.	
	


\end{document}